\documentclass[12pt,a4paper]{amsart}
\usepackage{graphicx}
\usepackage{epstopdf}
\usepackage{amssymb}
\usepackage{amsthm}
\usepackage{amsmath}
\usepackage{mathrsfs}
\usepackage{amscd}
\usepackage{enumerate}
\usepackage[all]{xy}
\newtheorem{thm}{Theorem}
\newtheorem{assum}[thm]{Assumption}
\newtheorem{cor}[thm]{Corollary}
\newtheorem{prop}[thm]{Proposition}
\newtheorem{lem}[thm]{Lemma}
\newtheorem{rem}[thm]{Remark}

\theoremstyle{definition}
\newtheorem{defn}[thm]{Definition}
\newtheorem{q}[thm]{Question}

\newtheorem{prop-def}[thm]{Proposition-Definition}
\newtheorem{example}[thm]{Example}
\newtheorem{rem-defn}[thm]{Remark-Definition}

\theoremstyle{remark}

\newcommand{\be}{\begin{equation}}
\newcommand{\bc}{\begin{cor}}
\newcommand{\bt}{\begin{thm}}
\newcommand{\bl}{\begin{lem}}
\newcommand{\bpr}{\begin{prop}}
\newcommand{\br}{\begin{rem}}
\newcommand{\bd}{\begin{defn}}
\newcommand{\ee}{\end{equation}}
\newcommand{\et}{\end{thm}}
\newcommand{\el}{\end{lem}}
\newcommand{\epr}{\end{prop}}
\newcommand{\er}{\end{rem}}
\newcommand{\ed}{\end{defn}}
\newcommand{\ec}{\end{cor}}

\newcommand{\R}{\Bbb{R}}

\begin{document}
\date{} 
\title{Graphs and obstruction theory for algebraic curves}
\author{Takeo Nishinou}
\thanks{email : nishinou@rikkyo.ac.jp}
\address{Department of Mathematics, Rikkyo University,
  Nishi-Ikebukuro, Toshima, Tokyo, Japan} 
\begin{abstract}
In this paper we study a construction of algebraic curves from 
 combinatorial data.
In the study of algebraic curves through degeneration, graphs usually appear as
 the dual intersection graph of the central fiber.
Properties of such graphs can be encoded in so-called tropical curves.
Our main concern is the relation between algebraic curves and tropical curves
 where the deformation problem is obstructed.
Particular emphasis is put on the role of higher valent vertices of 
 tropical curves, which has not been developed well so far in spite of its importance
 in this area of study.
We will give a general formula describing the obstruction,  
 a new criterion for the vanishing of the obstruction, 
 and a relation between the number of algebraic
 curves and the number of integral points in certain polytopes.
We also prove the optimal version of the correspondence between 
 tropical curves and algebraic curves when the tropical curves are regular, 
 generalizing \cite{CFPU, N3, Tyo}.
\end{abstract}
\maketitle
\section{Introduction}
The appearance of graphs in the study of algebraic curves is quite common.
In this context, 
 usually a graph appears as the dual intersection graph of the central fiber of
 a suitable degenerating family of algebraic curves.
Although any abstract finite graph can be realized as the dual intersection graph 
 associated to
 a degenerating family of algebraic curves (see \cite[Appendix B]{B}), 
 the situation becomes more interesting when one considers curves embedded in suitable
 ambient spaces.
Now the graphs have to satisfy balancing conditions at the vertices, which is the combinatorial 
 counterpart of the harmonicity, and there are subtler conditions coming from 
 superabundancy, which correspond to the obstructions to deform the degenerate
 algebraic curves.

The nature of the graphs associated to such degenerate algebraic curves
 can be encoded in so-called tropical curves.
G.Mikhalkin started this area of study with his pioneering work
 \cite{M}, in which he proved the correspondence between 
 imbedded tropical curves
 of any genus in $\R^{2}$, and holomorphic curves
 in toric surfaces specified by the combinatorial data of the tropical curves.
Recently ideas from non-archimedean geometry are pushing this field
 forward (see for example \cite{BPR}).

Nevertheless, the nature of degenerate algebraic curves is not yet well understood.
The main obstacle is the presence of obstructions
 to deform degenerate curves to smooth (or more generally irreducible) curves. 
The combinatorial counterpart of these curves are
 superabundant tropical curves.
In this paper, we study fundamental aspects of superabundant tropical curves
 and their implications to algebraic curves.
In particular, we investigate the role of
 higher valent vertices, which has not been studied 
 much so far in spite of its importance in the study of superabundant tropical curves.\\

The most fundamental problem concerning tropical curves is their relationship to
 algebraic curves.
This problem is essentially the same as the study of algebraic curves in a given variety
 though degeneration.
The basic strategy to handle this problem is two fold:
\begin{enumerate}
\item[Step 1:] First, construct a degenerate algebraic curve in a suitable ambient space, which is the central
 fiber of a degenerating family of varieties.
\item[Step 2:] Then, deform the curve to a general fiber of the family.
\end{enumerate}
When the tropical curve is not superabundant 
 (see Definition \ref{def:nonsuperabundant2} for the precise condition), 
 the relation to algebraic curves is quite nice in respect of the both steps
 above, and we have faithful correspondence between these curves
 (it is conditionally proved in \cite{CFPU, N3, Tyo}, and we prove it in the optimal
 form in this paper, see Theorem \ref{thm:regsm}).

On the other hand, when the tropical curve is superabundant, 
 several new phenomena emerge.
In this regard, Step 2 has been relatively well studied
 (see for example \cite{BPR, Katz, N3, S}).
However, it turns out that Step 1 is as crucial as Step 2,
 though this point does not seem to have been studied seriously.
The importance of this point becomes particularly clear when one considers
 tropical curves with higher valent vertices.
In this paper, we develop a general formalism to deal with higher valent vertices,
 and study the first and second steps using it.\\

The bridge which connected Step 1 to Step 2 is the calculation of the 
 cohomology group in which the obstruction to deform the degenerate algebraic curve lies.
On the other hand, the combinatorial counterpart of the obstruction is 
 the 
 superabundancy of tropical curves. 
Therefore, describing the superabundancy effectively is an important
 problem both in the theory of algebraic and tropical curves.
Namely, given a tropical curve, following problems immediately arise:
\begin{itemize}
\item Determine whether the tropical curve is superabundant or not.
\item When the tropical curve is superabundant, calculate
 the number of parameters of its deformation.
\end{itemize}
Both problems are not easily solved when one looks only at
 the tropical curve itself.
In Sections \ref{sec:dual obstruction} and \ref{sec:higherval}, 
 we will give a general answer to these problems by reducing them to a calculation of 
 certain sheaf cohomology group of the 
 algebraic counterpart of 
 tropical curves (Theorems \ref{thm:obstruction} and \ref{thm:obstruction2}).
When the tropical curve is 3-valent and is an  immersion, then this cohomology group is completely 
 determined by the combinatorial data of the tropical curve.
On the other hand, when there are higher valent vertices in the tropical curve, 
 this cohomology group is given as the solution space of certain system of linear equations,
 which is not in general determined by the combinatorial data.
With this description, we can study correspondence theorems for 
 superabundant tropical curves \cite{N3}.
It also gives basis for studies in 
 various other situations
 \cite{N4, N5, NY}.
An interesting point is that while usually tropical curves are used 
 to combinatorially understand algebraic curves, 
 here we use algebraic geometry 
 to understand combinatorics of tropical curves.

The fundamental set up for the study of higher valent vertices is 
 given in Section \ref{sec:higherval}.
The key point is that although algebraic geometric counter part of 
 a general higher valent vertex is not necessarily easy to handle, 
 such a vertex can be obtained from a standard higher valent vertex, 
 whose algebraic geometric counter part is a line in a projective space.
Using this, we can reduce various calculation on degenerate curves to 
 calculation on a line in a projective space, which can be done very 
 explicitly.
In particular, this enables us to calculate the obstruction 
 cohomology group in various situations.

In Sections \ref{sec:5} and \ref{sec:6}, we apply our formulation
 to concrete study of curves.
In Section \ref{sec:4}, we prepare some techniques useful for these studies.
Meanwhile they are also used to prove the optimal correspondence theorem 
 for regular tropical curves (Theorem \ref{thm:regsm}).
In Section \ref{sec:5}, we construct examples which exhibit subtleties 
 present in the first step and in the calculation of the obstruction
 cohomology groups.
The only known sufficient condition 
 which guarantees that a given superabundant tropical curve corresponds to an
 algebraic curve (such a tropical curve is called \emph{smoothable}
 in this paper)
 is so-called \emph{well-spacedness condition}
 for curves of genus one,
 first considered by Speyer \cite{S}. 
The situation becomes rather different in the presence of higher genera and 
 higher valent vertices.
In particular, we show the following new sufficient condition
 for smoothability of tropical curves of genus one (see Theorem \ref{thm:loop} for
 the precise statement).
\begin{thm}\label{thm:1}
Let $h\colon \Gamma\to\Bbb R^n$ be a tropical curve of genus one.
Assume that the directions of the edges emanating from the vertices on the loop
 span $\Bbb R^n$.
Then $(\Gamma, h)$ is smoothable if there is a degenerate algebraic curve
 corresponding to $(\Gamma, h)$.
\end{thm}
This is strikingly different from the well-spacedness condition.
For example, for an immersed 3-valent tropical curve of genus one to be well-spaced, 
 at least two vertices must exist from which edges
 not contained 
 in the minimal affine subspace containing the loop emanate. 
However, this is not the case in the situation of the above theorem.
Note that if $(\Gamma, h)$ is a superabundant tropical curve of genus one,
 then there must be higher valent vertices to fulfill the condition of the theorem.
This theorem can be extended to cover many situations where curves have higher genera
 (see Definition \ref{def:abundancysupport}).

Also note that the above theorem claims nothing about the existence of 
 a degenerate algebraic curve corresponding to $(\Gamma, h)$.
This point is also different from the cases to which the known well-spacedness condition
 applies.
Namely, when a tropical curve of genus one is an immersion, 
 one can show that there always exists a degenerate algebraic curve corresponding to it
 (that is, the first step can be always solved, see Corollary \ref{cor:genusonepre-logexist}).
However, as examples in Section \ref{sec:5} show, it is a subtle problem
 when there are higher valent vertices.
We give a general answer to this problem 
 in Section \ref{sec:6}.
The result is the following (see Theorem \ref{thm:pre-logexist} 
  for the precise statement)
\begin{thm}\label{thm:2}
Let $(\Gamma, h)$ be a superabundant tropical curve of genus one
 satisfying the condition of Theorem \ref{thm:1} (in particular, there are higher valent
 vertices on the loop).
There is a lattice parallelotope $\overline P$ in some $\Bbb R^N$ 
 and a finite set of hyperplanes $H_i$ constructed from the combinatorial data of
 $(\Gamma, h)$ such that a degenerate algebraic curve of type $(\Gamma, h)$
 exists if and only if
 $P\setminus (\cup_i H_i)$ contains a lattice point.
Here $P$ is the interior of $\overline P$.
\end{thm}
Also, it turns out that the number of families of degenerate 
 algebraic curves corresponding to $(\Gamma, h)$ is obtained by
 counting lattice points in $P\setminus(\cup_i H_i)$. 
Combining Theorems \ref{thm:1} and \ref{thm:2}, 
 we obtain the complete characterization of superabundant tropical curves 
 of genus one which correspond to algebraic curves, when there are no vertices
 outside the loop.
When there are vertices outside the loop, we need the calculation of 
 the obstructions contributed from them, and it is done in \cite{N3}.
Thus, the problem of the correspondence between algebraic curves and tropical curves
 of genus one has been solved (of course, there are still important problems
 of purely combinatorial nature, such as calculating the number of tropical curves
 of a given combinatorial type satisfying suitable constraint conditions).

The results so far apply to the study of algebraic curves on toric varieties through degeneration.
For the study of curves on more general varieties, we need to 
 know the effect of the singularities of the total space of the degeneration
 to the obstruction of curves.
Such calculation was done in \cite{N4} for the generic case.
Therefore, given any variety, 
 once we have a reasonable degeneration so that we can study degenerate 
 curves on it, we will be able to deduce large amount of information on the original variety from it, 
 at least about curves of genus at most one.

On the other hand, tropical varieties themselves are of much interest and 
 subject to active study.
In relation to our results, 
 Carolin Torchiani defined a nice counting number
 of tropical curves of genus one with desirable invariance properties
 in her thesis \cite{Tor}.
Her definition is based on the intersection theory on the moduli space of 
 tropical curves of genus one (which itself is a tropical variety).
In her definition, all the tropical curves satisfying the condition of 
 Theorem \ref{thm:1} contribute to the counting number positively.
However, Theorem \ref{thm:2} shows that the moduli space of 
 algebraic curves has more refined structure than the combinatorial one,
 reflecting the subtlety of Step 1 above.\\

\noindent
{\bf Acknowledgments.}
The author's
 study of tropical geometry began from the joint work \cite{NS} with
 B.Siebert, and many ideas from it appear in this paper, too.
It is a great pleasure for the author to express his gratitude to him.
The author would also like to thank C.Torchiani for sending him
 her thesis \cite{Tor}.
The author was supported by JSPS KAKENHI Grant Number
 26400061.
%$\bold{Notation.}$ 
%$N$ is a free abelian group of rank greater than or
% equal to 2.
%$N_{\Q} = N \otimes_{\Z} \Q$ and $N_{\R} = N \otimes_{\Z} \R$.
%$M = Hom_{\Z}(N, \Z)$.
%A toric variety is always regarded as a complex variety,
% sometimes with a symplectic structure induced from a dual polytope
% (this symplectic structure can be singular
% along lower dimensional toric strata,
% but this does not affect our argument).

\section{Preliminary}\label{sec:pre}
In this section, we recall and define some notations and notions
 which are used in this paper.

\subsection{Tropical curves}\label{subsec:pre}
First we recall some definitions about tropical curves, see \cite{M, NS}
for more information.
Let $\overline \Gamma$ be a weighted, connected finite graph.
Its sets of vertices and edges are denoted by $\overline \Gamma^{[0]}$,
$\overline \Gamma^{[1]}$, respectively.
We write by $w_{\overline \Gamma} \colon 
\overline \Gamma^{[1]} \to \Bbb N \setminus \{ 0 \}$
the weight function.
An edge $E \in \overline \Gamma^{[1]}$ has adjacent vertices
$\partial E = \{ V_1, V_2 \}$.
Let $\overline \Gamma^{[0]}_{\infty} \subset \overline \Gamma^{[0]}$
be the set of all 1-valent vertices.
We write $\Gamma = \overline \Gamma \setminus \overline\Gamma^{[0]}_{\infty}$.
Noncompact edges of $\Gamma$ are called \emph{unbounded edges}.
Let $\Gamma^{[1]}_{\infty}$ be the set of all unbounded edges.
Let $\Gamma^{[0]}, \Gamma^{[1]}, w_{\Gamma}$
be the sets of vertices and edges of $\Gamma$ and the weight function
of $\Gamma$ (induced from $w_{\overline\Gamma}$ in 
the obvious way),
respectively.
Let $N$ be a free abelian group of rank $n\geq 1$
and we write $N_{\Bbb K} = N\otimes_{\Bbb Z}\Bbb K$, where 
$\Bbb K = \Bbb Q, \Bbb R, \Bbb C$.

\begin{defn}
	We call a continuous map $h\colon \Gamma\to N_{\Bbb R}$ 
	a \emph{semi-affine map} if the following conditions hold.
	\begin{enumerate}
		\item $h$ is a proper map.
		In particular, the closure of the image of an unbounded edge is non-compact.
		\item For every edge $E \in \Gamma^{[1]}$, the restriction $h \big|_E$
		is either an embedding with the image $h(E)$ 
		contained in an affine line, or 
		a contraction so that $h(E)$ is a point.
	\end{enumerate}
\end{defn}

Let $h$ be a semi-affine map.
Let $\mathcal S\subset h(\Gamma)$ be the set of points
with the property that $p\in \mathcal S$ if and only if
for any neighborhood $O_p$ of $p$ in $N_{\Bbb R}$, 
the intersection $O_p\cap h(\Gamma)$ is not
homeomorphic to an open interval.
The following is easy to see.
\begin{lem}\label{lem:savert}
	Let $h$ be a semi-affine map as above.
	Then the following statements hold.
	\begin{enumerate}
		\item $\mathcal S$ is a finite set.
		\item By adding finite number of vertices to $\Gamma$, 
		we can assume that the inverse image $h^{-1}(p)$ 
		of each $p\in\mathcal S$
		consists of closed subgraphs of $\Gamma$.\qed
	\end{enumerate}
\end{lem} 

Under the condition of Lemma \ref{lem:savert} (2), the image $h(\Gamma)$
has a natural structure of a graph as follows.

\begin{defn}\label{def:imagevert}
	Let $h\colon \Gamma\to N_{\Bbb R}$ be a semi-affine map satisfying 
	the condition of Lemma \ref{lem:savert} (2).
	Then a vertex of the image $h(\Gamma)$ is the image of some vertex of $\Gamma$.
	Similarly, an edge of $h(\Gamma)$ is the image of some edge of $\Gamma$
	which is not contracted.
\end{defn}

Note that under the condition of Lemma \ref{lem:savert} (2), 
the inverse image $h^{-1}(v)$ of a vertex of $h(\Gamma)$ is a union of 
closed subgraphs.
Similarly, the inverse image $h^{-1}({\mathfrak E}^{\circ})$ of the open part
 ${\mathfrak E}^{\circ}$ of an edge 
$\mathfrak E$ of $h(\Gamma)$ is the union of the open part of the edges of $\Gamma$
each of which is mapped to $\mathfrak E$ homeomorphically.

Henceforth, we always assume 
Lemma \ref{lem:savert} (2) holds.
Let $p\in\mathcal S$ and $\Gamma_1$ be one of connected components of
$h^{-1}(p)$.
Then $\Gamma_1$ contains several 1-valent vertices $q_1, \dots, q_a$.
Let $E_1, \dots, E_b$ be the edges of $\Gamma\setminus\Gamma_1$
emanating from some of 
$q_1, \dots, q_a$.

\begin{defn}[{\cite[Definition 2.2]{M}}]\label{def:param-trop}
	A \emph{parametrized tropical curve} in $N_{\Bbb R}$ is a semi-affine map
	$h \colon \Gamma \to N_{\Bbb R}$ satisfying the following conditions.
	\begin{enumerate}
		\item[(i)] For every edge $E \in \Gamma^{[1]}$, 
		the image $h(E)$ is either contained
		in an affine line with a rational slope, or 
		a point.
		\item[(ii)] For every vertex $V \in \Gamma^{[0]}$, $h(V)\in N_{\Bbb Q}$.
		\item[(iii)] 
		The following \emph{balancing
			condition} holds.
		Namely, for each $p\in\mathcal S$ and a connected component
		$\Gamma_1$ of $h^{-1}(p)$, 
		the equality
		\begin{equation}\label{eq:balancingcond}
		\sum_{j=1}^b w_{\Gamma}(E_j)u_j = 0
		\end{equation}
		holds, using the notation in the above paragraph.
		Here $u_j$ is the primitive integral vector of $N$ 
		in the direction of the edge $h(E_j)$ emanating from $p$.
	\end{enumerate}
\end{defn}
\begin{rem}
In \cite{NS}, $h|_E$ is assumed to be an embedding
 $\text{(see \cite[Definition 1.1]{NS})}$ for every edge $E$.
The reason why we adopt the above definition 
 is that these cases appear naturally
 when we consider superabundant tropical curves.
\end{rem}
An isomorphism between parametrized tropical curves 
 $h \colon \Gamma \to N_{\Bbb R}$ and 
 $h' \colon \Gamma' \to N_{\Bbb R}$ 
 is a homeomorphism $\Phi \colon \Gamma \to \Gamma'$
 respecting the weights such that $h = h' \circ \Phi$.
\begin{defn}\label{def:tropical curve}
A \emph{tropical curve} is an isomorphism class of parametrized 
 tropical curves.
A tropical curve is \emph{3-valent} if any vertex of $\Gamma$ is at  most
 3-valent.
The \emph{genus} of a tropical curve is the first Betti number of $\Gamma$.
The set of \emph{flags} of $\Gamma$ is 
\[
F\Gamma = \{(V, E) \big|
     V \in \partial E \}.
     \]
\end{defn}
Note that in our definition, 
 there can be 2-valent vertices in the abstract graph $\Gamma$
 of a 3-valent tropical curve.
By (i) of Definition \ref{def:param-trop}, we have a map
 $u \colon F\Gamma \to N$ sending a flag $(V, E)$
 to the primitive integral vector $u_{(V, E)} \in N$
 emanating from $h(V)$ in the direction of $h(E)$ or to the zero vector. 
\begin{rem}\label{rem:0edge}
We allow the case 
 that even when $u_{(V, E)}\neq 0$, the edge $E$ is contracted
 by the map $h$
 $($but in such a case the direction $u_{(V, E)}$ does not contribute to the 
 balancing condition $\text{(\ref{eq:balancingcond}))}$.
We think of such a tropical curve as the limit of a
 family of tropical curves $(\Gamma, h_s)$, 
 $s\in [0, 1)\cap \Bbb Q$ where the image $h_s(E)$ has the direction $u_{(V, E)}$
 for each $s$ 
 and the length goes to zero as $s\to 1$
 (see Definition \ref{def:moduli}).
\end{rem}
\begin{defn}\label{def:type}
The (unmarked) 
 \emph{combinatorial type} (or simply the \emph{type}) of a tropical curve $(\Gamma, h)$
 is the graph $\Gamma$ 
 together with the map $u \colon F\Gamma \to N$.
We write this by the pair $(\Gamma, u)$.
\end{defn}

\begin{defn}\label{def:degree}
The \emph{degree} of a combinatorial type $(\Gamma, u)$
 is the function $\Delta(\Gamma, u) = \Delta\colon N \setminus \{ 0 \}
  \to \Bbb N$
 with finite support defined by 
\begin{equation*}
 \Delta(\Gamma, u)(v):= \sharp \{ (V, E) \in F\Gamma |
    E \in \Gamma^{[1]}_{\infty}, w(E)u_{(V, E)} = v \}. 
\end{equation*}
Let $e=|\Delta| = \sum_{v\in N\setminus\{0\}}\Delta(v)$.
This is the same as the number of unbounded edges 
 of $\Gamma$ (not necessarily of $h(\Gamma)$
 since some of the edges may have the same image).
\end{defn}
\begin{defn}\label{immersive}
We call a tropical curve $(\Gamma, h)$ \emph{immersive} 
 if for any $E\in\Gamma^{[1]}$, the restriction of 
 $h$ to $E$ is an embedding.
Note that even if $(\Gamma, h)$ is immersive, 
 some of the edges of $\Gamma$ can have the same image. 
\end{defn}
\begin{prop}[{\cite[Proposition 2.13]{M}}]\label{prop:trop_moduli}
The space parameterizing immersive 3-valent tropical
 curves of a given combinatorial type is, if it is
 non-empty, an open
 convex polyhedral domain in the real affine $k$-dimensional space,
 where $k \geq e+(n-3)(1-g)$.
Here $e$ is the number of unbounded edges of $\Gamma$ as in Definition \ref{def:degree},
 $n$ is the dimension of the target space $N_{\Bbb R}$, and $g$ is the genus of $\Gamma$. \qed
\end{prop}
\begin{defn}\label{def:moduli}
Fix a combinatorial type of 3-valent tropical curves whose parameter space of 
 immersive curves is non-empty.
Then we define
 the {parameter space} of \emph{all} 3-valent tropical
 curves of the given combinatorial type as the closure
 of the parameter space of 
 immersive curves.
\end{defn}
\begin{rem}\label{rem:type}
In other words, an element of the parameter space of 
 tropical curves of a given combinatorial type is a tropical curve
 which can be deformed into an immersive tropical curve
 of that combinatorial type.
If there is no immersive tropical curve of the given combinatorial type, 
 then the corresponding parameter space is defined to be empty
 in this paper (see Example \ref{ex:singular}).
See also Proposition \ref{prop:immexist} and Question \ref{q:1}.
\end{rem}

%

%The parameter space in the statement of
% Proposition \ref{prop:trop_moduli} is defined by a set of 
% affine linear inequalities which we write symbolically by $f_i>0$.
%If we allow some edges of $\Gamma$ contracted by the map $h$, 
% some of the inequalities $f_i>0$ become not strict: $f_i\geq 0$.
%In particular, we have the following.
The following is an immediate consequence of Proposition \ref{prop:trop_moduli}.
\begin{cor}\label{cor:trop_moduli}
Let $(\Gamma, h)$ be a 3-valent tropical curve.
Then the parameter space of 3-valent tropical curves
 of the given combinatorial type is, if it is not empty, a closed
 convex polyhedral domain in a real affine $k$-dimensional space,
 where $k \geq e+(n-3)(1-g)$.
\qed
\end{cor}
\begin{rem}
The term 
 'closed' in the statement of the corollary does not mean 'compact'.
Since we can parallel transport any tropical curve, 
 the parameter space is always non-compact.
\end{rem}
Before stating Assumption A below, we prepare some terminologies.
Let $\Gamma$ be a non-compact finite graph as above.
\begin{defn}\label{def:loops}
\begin{enumerate}[(i)]
\item An edge $E\in\Gamma^{[1]}$ is said to be a 
 \emph{part of a loop} of $\Gamma$
 if the graph given by $\Gamma\setminus E^{\circ}$ has the smaller first Betti number
 than $\Gamma$.
Here $E^{\circ}$ is the interior of $E$ (that is, 
 $E^{\circ} = E\setminus\partial E$).
\item The \emph{loop part} of $\Gamma$ 
 is the subgraph of $\Gamma$ composed of
 the union of all parts of the loops of $\Gamma$.
\item A \emph{bouquet} of $\Gamma$ is 
 a connected component of the loop part of $\Gamma$.
A subset of a bouquet which is homeomorphic to
 a circle is called a \emph{loop}.
\end{enumerate}
\end{defn}
In particular, a bouquet or a loop does not contain unbounded edges.

\begin{example}\label{ex:singular}
Proposition \ref{prop:trop_moduli} fails to hold for 
 non-immersive tropical curves.
For example, consider an abstract 3-valent graph $\Gamma$ 
 which has three unbounded edges $E_1, E_2, E_3$ of weight one.
Assume that the set $\Gamma\setminus 
 \{E^{\circ}_1, E^{\circ}_2, E^{\circ}_3\}$
 is a bouquet.
The following map $h\colon \Gamma\to \Bbb R^2$
 gives a tropical curve.
\begin{itemize}
\item $h$ maps the ends of $E_1, E_2, E_3$ 
 to the origin $(0, 0)\in\Bbb R^2$.
\item $h$ maps the edges $E_1, E_2, E_3$ onto the half lines
 $\{(x, 0)\;|\; x\geq 0\}, \{(0, y)\;|\; y\geq 0\}, \{(x, y)\;|\; x = y\leq 0\}$, 
 respectively.
\item $h$ contracts the other part of $\Gamma$ to $(0, 0)\in\Bbb R^2$.
\end{itemize}
Then it is easy to see that  
 there is no deformation of $(\Gamma, h)$ other than 
 parallel transports.
Therefore, if the genus of $\Gamma$ is positive and $h$ is not immersive,
 Proposition \ref{prop:trop_moduli}
 does not always hold.
This tropical curve is not contained in the parameter space
 in the sense of Definition \ref{def:moduli}
 for any combinatorial type.
 
Also, even when $\Gamma$ is tree, $(\Gamma, h)$ need not
 deform into an immersive tropical curve, see Figure \ref{fig:proj}.

\begin{figure}[h]
	\includegraphics{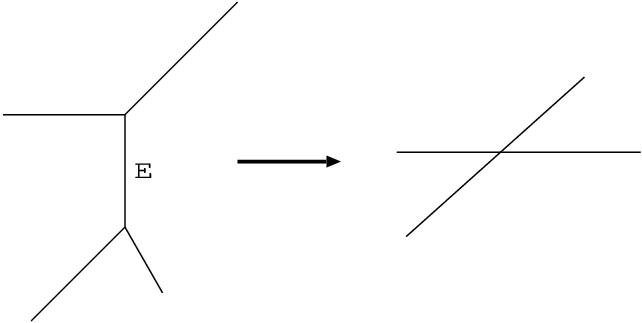}
	\caption{A tropical curve in $\Bbb R^2$ which cannot be
		deformed into an immersive curve, even though the domain
		abstract graph is a tree.
		The edge $E$ in the abstract graph
		is contracted to the unique vertex of the image.
		The two upper edges of the abstract graph are mapped into 
		the horizontal line in $\Bbb R^2$.
		Similarly, the two lower edges are mapped into the other line.
	}\label{fig:proj}
\end{figure}
\end{example}
As this example shows, when the map $h$ contracts loops of $\Gamma$, 
 it becomes difficult to give a unified treatment of tropical curves.
Therefore we introduce the following assumption. \\

\noindent
{\bf Assumption A.}
\begin{enumerate}[(i)]
\item The abstract graph $\Gamma$ is always 3-valent.
Therefore, $(\Gamma, h)$ is always a 3-valent tropical curve in view of
 Definition \ref{def:tropical curve},
 although some vertices of the image $h(\Gamma)$ may not be 2- or 3-valent.
\item The map $h$ may contract some
 of the bounded edges of $\Gamma$.
However, $h$ does not contract a loop to a point.
\item $(\Gamma, h)$ can be deformed into an immersive tropical curve.
Thus, the parameter space containing $(\Gamma, h)$ is not empty.
\item The inverse image of any vertex of $h(\Gamma)$ is a 
disjoint union of closed subgraphs of $\Gamma$ (by (ii), 
each component must be a tree or a vertex).
\end{enumerate}
%{\bf (2)の例の絵を入れる。loopのedgeだけでなく, loopに隣接するedgeも
%contractされない。}
\begin{rem}\label{rem:assump}
\begin{enumerate}
\item The essential part of Assumption A is the part (iii). 
The part (ii) is also related to the part (ii) (see Example \ref{ex:singular}).
The part (iv) is (2) of Lemma \ref{lem:savert}, and can always be realized by taking
 a suitable refinement of  the graph
$\Gamma$.
\item In particular, immersive 3-valent tropical curves satisfy Assumption A
 possibly after a suitable refinement of the graph $\Gamma$.
\item When $(\Gamma, h)$ satisfies Assumption A, we
 define a part of a loop, the loop part, a bouquet and a loop 
 of the image 
 $h(\Gamma)$ as the images of them of $\Gamma$.
\item By (iii) of Assumption A, 
 when a tropical curve $(\Gamma, h)$ satisfies Assumption A,
 the parameter space of the given type is a non-empty
 closed 
 convex polyhedral domain in a real affine $k$-dimensional space,
 where 
 $k \geq e+(n-3)(1-g)$, by Corollary \ref{cor:trop_moduli}.
The subset of tropical curves satisfying Assumption A
 is neither closed nor open in general,
 since some of the edges can be contracted while edges in a loop
 cannot be contracted simultaneously.
\item A general member of the parameter space in 
 (4) above is immersive.
\end{enumerate}
\end{rem}

Now we define the superabundancy of tropical curves.
In view of Example \ref{ex:singular}, it is reasonable to define it
 only for curves satisfying Assumption A(iii).
\begin{defn}[{\cite[Definition 2.22]{M}}]\label{def:superabundancy}
A 3-valent tropical curve satisfying Assumption A(iii)
 is called \emph{superabundant}
 if the parameter space is of dimension larger than 
 $e+(n-3)(1-g)$.
Otherwise it is called \emph{non-superabundant}.
See also Subsection \ref{subsec:compare}.
\end{defn}

To see whether a tropical curve 
 satisfying Assumption A
 of a given combinatorial type is superabundant or not, 
 it is enough to check it for an immersive tropical curve obtained by 
 deforming the original curve. 
On the other hand, we will see below that 
 the superabundancy of an 
 immersive tropical curve can be 
 effectively calculated via algebraic geometry.
Now we recall some notions from algebraic geometry 
 relevant to our purpose.

\subsection{Toric varieties associated to tropical curves and pre-log curves on them}\label{subsec:pre-logdef}
\begin{defn}\label{toric}
A toric variety $X$ defined by a fan $\Sigma$ 
 is called \emph{associated to a tropical curve $(\Gamma, h)$}
 if the set of the 
 rays of $\Sigma$ contains the set of the rays spanned by the vectors in
 $N$ which are in the 
 support of the degree map $\Delta\colon N\setminus\{0\}\to \Bbb N$
 of $(\Gamma, h)$.

If $\mathfrak E$
 is an unbounded edge of $h(\Gamma)$, there is an obvious unique divisor of $X$
 corresponding to it.
We write it by $D_{\mathfrak E}$ and call it the \emph{divisor associated to the edge 
 $\mathfrak E$}.
\end{defn}
\begin{defn}\label{def:degeneration}
Given a tropical curve $(\Gamma, h)$ in $N_{\Bbb R}$
 defined over $\Bbb Q$ (that is, the vertices have rational coordinates),
 we can construct a polyhedral
 decomposition $\mathscr P$ of $N_{\Bbb R}$ defined over $\Bbb Q$
 such that $h(\Gamma)$ is contained in the 1-skeleton of $\mathscr P$
 (\cite[Proposition 3.9]{NS}).
Given such $\mathscr P$, 
 we construct a degenerating family $\mathfrak X\to \Bbb C$
 of
 a toric variety $X$ associated to $(\Gamma, h)$ (\cite[Section 3]{NS}).
We call such a family a \emph{degeneration of $X$ defined respecting $(\Gamma, h)$}.
Let $X_0$ be the central fiber.
It is a union $X_0 = \cup_{v\in\mathscr P^{[0]}}X_{0,v}$
 of toric varieties intersecting along toric strata.
Here $\mathscr P^{[0]}$ is the set of the vertices of $\mathscr P$.
\end{defn}

\begin{rem}
To define $\mathfrak X$, in general we need to multiply $(\Gamma, h)$ by some constant
 so that it becomes defined over $\Bbb Z$.
There are choices of this multiplication factor, which correspond to 
 base changes on the algebraic geometry side.
One requirement to this choice is that the integral length of each bounded edge is a multiple of
 its weight (see \cite[Proposition 7.1]{NS}).
%If this condition is satisfied, the difference of the multiplication factor has essentially no
% effect (see Remark \ref{rem:2-val}).
In the following arguments,
 we just talk about degenerations defined 
 respecting $(\Gamma, h)$ which is defined over $\Bbb Q$, 
 assuming a suitable choice of the multiplication factor is made.
It is convenient to do this, since we sometimes refine the graph by adding 2-valent vertices,
 and mentioning this base change each time will be cumbersome.
\end{rem}

\begin{defn}[{\cite[Definition 4.1]{NS}}]\label{torically transverse}
Let $X$ be a toric variety.
A holomorphic curve $C\subset X$ is \emph{torically transverse}
 if it is disjoint from all toric strata of codimension greater than one.
A stable map $\phi\colon C\to X$ 
 is torically transverse if $\phi^{-1}(int X)\subset C$
 is dense and $\phi(C)\subset X$ is a torically transverse curve. 
Here $int X$ is the complement of the union of toric divisors.
\end{defn}
\begin{defn}\label{def:pre-log}
Let $C_0$ be a prestable curve.
A \emph{pre-log curve} on $X_0$ is a stable map 
 $\varphi_0\colon C_0\to X_0$
 with the following properties.
\begin{enumerate}
\item[(i)] For any $v\in\mathcal P^{[0]}$,
 the restriction $C_0\times_{X_0}X_{0,v}\to X_{0,v}$
 is a torically transverse stable map.
\item[(ii)] Let $P\in C_0$ be a point which maps to the singular locus of $X_0$.
Then $C_0$ has a node at $P$, and $\varphi_0$ maps the two branches
 $(C_0', P), (C_0'', P)$ of $C_0$ at $P$ to different irreducible components 
 $X_{0, v'}, X_{0, v''}\subset X_0$.
Moreover, if $w'$ is the intersection index 
 of the restriction $(C_0', P)\to (X_{0, v'}, D')$ with the toric divisor
 $D'\subset X_{0, v'}$, 
 and $w''$ accordingly for $(C_0'', P)\to (X_{0, v''}, D'')$,
 then $w' = w''$.
\end{enumerate}
\end{defn}
Let $X$ be a toric variety and $D$ be the union of toric divisors.
In \cite[Definition 5.2]{NS}, a non-constant
 torically transverse map $\phi\colon \Bbb P^1\to X$
 is called a \emph{line}
 if 
\[
\sharp\phi^{-1}(D)\leq 3.
\]
In this case, the image of $\phi$ is contained in the closure of 
 the orbit of a subtorus of dimension at most two 
 of the big torus acting on $X$ (\cite[Lemma 5.2]{NS}).
Because we consider 
 more general tropical curves, we have to extend this notion.

Let $\Gamma$ be a weighted 3-valent tree 
 and $h\colon \Gamma\to N_{\Bbb R}$
 be a map which gives $\Gamma$ a structure of a tropical curve,
 and assume that the image $h(\Gamma)$ has only one vertex $v$.
Let $\mathfrak E_1, \dots, \mathfrak E_s$ be the edges of $h(\Gamma)$
 (some of 
 these can coincide, and
 these correspond to the unbounded edges of $\Gamma$ in a natural way).
% and 
% $W_1, \dots, W_s$ be the weights (note that 
% in general $W_i$ is a sequence of integers).
Let $X$ be a toric variety associated to $(\Gamma, h)$.
\begin{defn}\label{def:typev}
A non-constant torically transverse map
 $\phi\colon \Bbb P^1\to X$ is called 
 \emph{of type $(\Gamma, h)$}, or  \emph{of type $v$}
 when $h$ is clear from the context,
 if $\phi$ satisfies the following property:
\begin{itemize}
\item Let $\mathfrak E_{i}$ be an edge of $h(\Gamma)$ and
 let $w_i$ be the weight of $\mathfrak E_{i}$.
Then $\phi(\Bbb P^1)$ has an intersection
 with the divisor $D_{\mathfrak E_i}$
 with 
  intersection multiplicity 
 $w_i$, and there is no intersection between $\phi(\Bbb P^1)$
 and toric divisors other than these.
\end{itemize} 
\end{defn}
Note that when some edges of $\mathfrak E_1, \dots, \mathfrak E_s$
 coincide, then there are several intersections between $\phi(\Bbb P^1)$
 and the corresponding toric divisor.
 
\begin{figure}[h]
	\includegraphics[height=4cm]{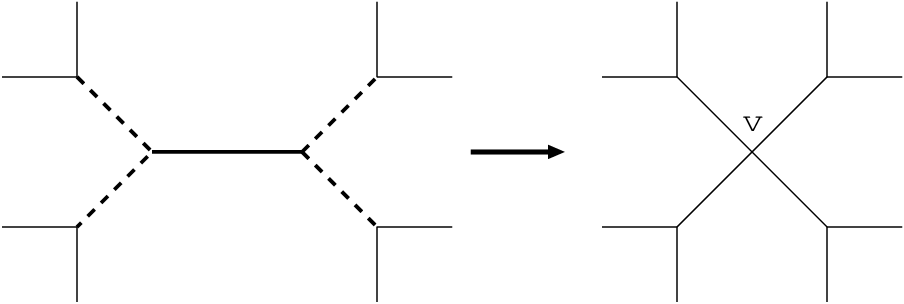}
	\caption{The abstract graph $\Gamma$ (the picture on the left)
		is mapped to a graph in $\Bbb R^n$ (the picture on the right).
		The bold line segment in $\Gamma$ is the inverse image of the vertex $V$
		and the union of the bold line segment and the dotted line segments
		in $\Gamma$ is the subgraph $\Gamma_V$.
		In this example, the graph $\Gamma_{\varphi_0}$
		is the same as the picture on the right as an abstract graph.}\label{fig:invim}
\end{figure}

Let $(\Gamma, h)$ be a tropical curve satisfying Assumption A.
Let $X$ be a toric variety associated to $(\Gamma, h)$ and  
 $\mathfrak X\to\Bbb C$ be a degeneration of $X$ defined respecting $(\Gamma, h)$.
Let $X_0$ be the central fiber.
\begin{defn}\label{def:oftype}
A pre-log curve $\varphi_0\colon C_0\to X_0$ is called 
 \emph{of type $(\Gamma, h)$}
 if for any $v\in h(\Gamma)^{[0]}$,  
 the restriction of $C_0\times_{X_0}X_{0, v}\to 
  X_{0, v}$ to each of the components $\{C_{0, v(i)}\}$
  of $C_0$
  mapped to $X_{0, v}$
  is a rational torically transverse curve of type 
  $(\Gamma_{v(i)}, h|_{\Gamma_{v(i)}})$.
 Here
  $\{\Gamma_{v(i)}\}$ is the set of open
  subgraphs of $\Gamma$ each element of which is 
  the union of 
 \begin{itemize}
 \item a connected component of $h^{-1}(v)$, and 
 \item the open
  part of the edges emanating from the vertices of it.
 \end{itemize}
See Figure \ref{fig:invim} for an example.
\end{defn}

\begin{defn}\label{def:intergraph}
When $\varphi_0\colon C_0\to X_0$ is a pre-log curve of type
 $(\Gamma, h)$, let $\Gamma_{\varphi_0}$ be the graph obtained from 
 $\Gamma$ by contracting 
 all the bounded edges of each $\Gamma_{v(i)}$
 in the notation of Definition \ref{def:oftype}.
The graph $\Gamma_{\varphi_0}$ is the dual intersection graph of $C_0$
 (containing the set of unbounded edges, which is
 the information of intersection of $\varphi_0(C_0)$
 and some of the toric divisors of $X_0$ seen as marked points of $C_0$)
 and it is clear that the map $h\colon \Gamma\to \Bbb R^n$
 factors through $\Gamma_{\varphi_0}$ (see Figure \ref{fig:invim}).
\end{defn}

The components $\{C_{0, v(i)}\}$ 
 are disjoint, since otherwise the graphs $\Gamma_{v(i)}$
 and $\Gamma_{v(i')}$ have a non-contracting common (open) edge $E$
 for some $i$ and $i'$, 
 but since both of the ends of $E$ are mapped to the vertex $v$, 
 the edge $E$ must be contracted.
Note that each $C_{0, v(i)}$ is an irreducible rational curve.

\begin{rem}[{\cite[Definition 5.6]{NS}}]\label{def:maximally degenerate}
A pre-log curve $\varphi_0\colon C_0\to X_0$ 
 is \emph{maximally degenerate} if
 it is of type $(\Gamma, h)$ for some \emph{immersive}
 tropical curve $(\Gamma, h)$. 
\end{rem}

\begin{defn}\label{def:smoothable}
A 3-valent 
 tropical curve $(\Gamma, h)$  is \emph{smoothable}
 if there is a pre-log curve $\varphi_0\colon C_0\to X_0$ 
 of type $(\Gamma, h)$ with the following property.
Namely, there exists 
 a family of stable maps over $\Bbb C$
\[
\Phi\colon \mathfrak C/\Bbb C\to \mathfrak X/\Bbb C
\]
 such that $\mathfrak C/\Bbb C$ is a flat family of pre-stable 
 curves whose fiber over $0$ is isomorphic to $C_0$,
 and the restriction of $\Phi$ to $C_0$ is a stable map
 equivalent to $\varphi_0$.
We also call such a pre-log curve \emph{smoothable}.
\end{defn}
\begin{rem}
The smoothability of a tropical curve does not depend on the choice of
 a toric variety
 $X$ associated to it or a degeneration of $X$ defined respecting the tropical curve.
\end{rem}
See \cite[Section 5]{NS}, for more information about lines and maximally
  degenerate pre-log curves.
Given an immersive 3-valent tropical curve
 $(\Gamma, h)$ which is non-superabundant,
 we can construct maximally degenerate 
 pre-log
 curves of type $(\Gamma, h)$
 (\cite[Proposition 5.7]{NS}), and vice versa (\cite[Construction 4.4]{NS}).
%The arguments there together with some additional arguments
% extend this result to not necessarily 
% immersive non-superabundant
% 3-valent tropical curves 
% and pre-log curves of type $(\Gamma, h)$.
%In fact, if $(\Gamma, h)$ satisfies
% Assumption A, 
% there is a pre-log curve of type $(\Gamma, h)$ (see Corollary \ref{cor:genexist}).
More generally, if $(\Gamma, h)$ is regular in the sense of Definition \ref{def:nonsuperabundant2},
 then we will show that there is a pre-log curve of type $(\Gamma, h)$ (see Proposition \ref{prop:exist}).
Also, it is not difficult to see that immersive non-superabundant
 tropical curves are regular.
However, when the tropical curve is superabundant, 
 a degenerate curve of the given type might not exist, see 
 Examples \ref{ex:1} and \ref{ex:1.5}.

%We put marked points on these maximally
% degenerate 
% pre-log
% curves, corresponding to the edges of the tropical curves.
%Namely, a marked point is a node when it corresponds to a bounded edge of
% the tropical curve, while a marked point is a non-singular point when it corresponds
%  to an unbounded edge.
The smoothing of these degenerate curves 
 of type $(\Gamma, h)$ is studied in the following sections
 through log-smooth deformation theory \cite{KF, KK}.
For informations about 
 log structures relevant to our situation, see \cite[Section 7]{NS}.
We do not repeat it here, because nothing new about log structures is required here, 
 other than those given in \cite{NS}. 
 
%\begin{rem}\label{rem:2-val}
%	In general the graph $\Gamma$ must contain 2-valent vertices
%	to satisfy (iv) of Assumption A.
%	However, since the components of pre-log curves corresponding to 
%	these vertices play rather minor role (see the proof of \cite[]{NS}, 
%	for example), we ignore these vertices
%	in many cases to simplify the exposition.
%\end{rem}
% 

\section{Combinatorial description of the
 dual space of obstructions for immersive tropical curves}\label{sec:dual obstruction}
Let $(\Gamma, h)$ be an immersive
	tropical curve and $X$ be a toric variety associated
	to it.
Let $\mathfrak X\to \Bbb C$ be a degeneration of $X$ 
	defined respecting $(\Gamma, h)$.
Let $\mathscr P$ be a polyhedral decomposition of $N_{\Bbb R}$ defining $\mathfrak X$.
Assume there is a pre-log curve $\varphi_0\colon C_0\to X_0$
 of type $(\Gamma, h)$ (In fact, the result in this section makes sense
 even if such a curve does not exist, see Remark \ref{rem:thm40add}).
\begin{rem}\label{rem:2val}
It is in general necessary to allow 2-valent vertices
  to $\Gamma$ in order to assure the property
	$h^{-1}(\mathscr P^{[0]}) = \Gamma^{[0]}$.
Then
	there are components of $C_0$ corresponding to these 2-valent vertices.
However, as in the proof of \cite[Proposition 7.1]{NS}, 
 these components do not play an essential role in the argument below. 
Therefore, we usually neglect them and regard $\Gamma$ as 
	if it has only 3-valent
	vertices for the simplicity of presentation. 
\end{rem}

We can give log structures to $C_0$ and $\mathfrak X$
 as in \cite[Section 7]{NS}.
There are log tangent sheaves associated to these
 log structures.
The tangent space and the obstruction space of the 
 deformation of $\varphi_0$ are calculated in terms of these sheaves.

%Proposition 8.2 of \cite{N} extends in a straightforward manner to our general case.
Suppose that a lift
$\varphi_{k-1}\colon C_{k-1}/O_{k-1}\to\mathfrak X$ of $\varphi_0$ is constructed.
Here $O_{k-1} = \Bbb C[\epsilon]/\epsilon^k$.
Then as in the proof of \cite[Lemma 7.2]{NS},  an extension $C_{k}/O_{k}$
of $C_{k-1}/O_{k-1}$ exists and such extensions are
parametrized by the space $H^1(C_0, \Theta_{C_0/O_0})$, here 
$\Theta_{C_0/O_0}$ is the log tangent sheaf.

On the other hand, 
the obstruction to lift the map $\varphi_{k-1}$ to the next order
 is given by a cohomology class in 
$H^{1}(C_{0}, \varphi^*_{0}\Theta_{\mathfrak X/\Bbb C})$
(note that even if we consider the deformation of 
$\varphi_{k-1}\colon C_{k-1}\to \mathfrak X$, the obstruction lies in 
the cohomology on $C_0$),
here $\Theta_{\mathfrak X/\Bbb C}$ is the log tangent sheaf 
relative to the base.
As in usual deformation theory of smooth varieties, there is a following 
standard result in the log smooth deformation theory \cite{KK}.
\begin{prop}\label{prop:nonobst}
If $H^1(C_0, \varphi_0^*\Theta_{\mathfrak X/\Bbb C})$ vanishes, 
	then the pre-log curve $\varphi_0$ is smoothable.\qed
\end{prop}

On the other hand,
 the sheaf $\varphi^*_{0}\Theta_{\mathfrak X/\Bbb C}$ fits in the exact
sequence
\begin{equation}\label{eq:deform}
	0 \to \Theta_{C_{0}/O_{0}} \to \varphi^*_{0}\Theta_{\mathfrak X/\Bbb C}
	\to \varphi^*_{0}\Theta_{\mathfrak X/\Bbb C}/\Theta_{C_{0}/O_{0}} \to 0.
\end{equation}
Note that
we have the natural isomorphism
\[
\Theta_{\mathfrak X/\Bbb C} \simeq N\otimes_{\Bbb Z}
\mathcal O_{\mathfrak X},
\]
here $N$ is the free abelian group such that
the fan defining a general fiber of $\mathfrak X$ lies in 
$N_{\Bbb R} = N\otimes\Bbb R$.

We have a map between cohomology groups
\[
H^1(C_{0}, \Theta_{C_{0}/O_{0}}) \to
H^1(C_{0}, \varphi^*_{0}\Theta_{\mathfrak X/\Bbb C}).
\]
The group  
 $H^1(C_{0}, \Theta_{C_{0}/O_{0}})$
 is the tangent space of the moduli space of deformations
 of $C_{0}$, and the obstruction classes in 
 $H^1(C_{0}, \varphi^*_{0}\Theta_{\mathfrak X/\Bbb C})$
 which are in the image of the above map
 can be cancelled when we deform 
 the moduli of the domain of the stable maps.
Namely, we have the following standard fact.
\begin{prop}
If the map $H^1(C_{0}, \Theta_{C_{0}/O_{0}}) \to
 H^1(C_{0}, \varphi^*_{0}\Theta_{\mathfrak X/\Bbb C})$
 is a surjection, then the map $\varphi_0$ is smoothable.
\end{prop}
\proof
Given the map $\varphi_0$, the obstruction to deform it
is defined, and it determines a class 
$\alpha\in H^1(C_{0}, \varphi^*_{0}\Theta_{\mathfrak X/\Bbb C})$.
This class can be calculated in the following way.
Take a suitable open covering $\{U_i\}$ of $C_0$ so that 
there is a first order
lift $\varphi_1|_{U_i}\colon \tilde U_i\to \mathfrak X$
of the restriction of $\varphi_0$ to each of the open set
$U_i$
(the existence of such a covering follows from the general theory
of log smooth deformations).
The set of lifts on $U_i$ forms a torsor over the abelian group of sections 
$\Gamma(U_i, \varphi^*_{0}\Theta_{\mathfrak X/\Bbb C})$, 
and the differences of the lifts on the intersections $U_i\cap U_j$
determine a $\varphi^*_{0}\Theta_{\mathfrak X/\Bbb C}$-valued
\v{C}ech 1-cocycle, which represents the class $\alpha$.

On the other hand, by assumption, the class $\alpha$ is mapped to 
$0\in H^1(C_0, \varphi^*_{0}\Theta_{\mathfrak X/\Bbb C}
/\Theta_{C_{0}/O_{0}})$.
This implies that we can perturb each lift $\varphi_1|_{U_i}$ by a section in
$\Gamma(U_i, \varphi^*_{0}\Theta_{\mathfrak X/\Bbb C})$
so that the images of the
differences of the lifts in $\varphi^*_{0}\Theta_{\mathfrak X/\Bbb C}
/\Theta_{C_{0}/O_{0}}$ are zero not only cohomologically but also
in the level of cocycles.
Then the differences of the lifts give $\Theta_{C_{0}/O_{0}}$-valued
\v{C}ech 1-cocycle, 
where we see $\Theta_{C_{0}/O_{0}}$ as a subsheaf of 
$\varphi^*_{0}\Theta_{\mathfrak X/\Bbb C}$ in the natural way.

In particular, the class $\alpha$ can be seen as a class 
of $H^1(C_0, \Theta_{C_0/O_0})$.
This implies that the domain $\tilde U_i$ of the local lift
$\varphi_1|_{U_i}$ glues into a global lift $C_1$ of $C_0$, 
so that the natural map $C_1\to \Bbb C[t]/t^2$ gives
a deformation of $C_0$ corresponding to the class $\alpha$, 
and the maps $\varphi_1|_{U_i}$ can be also glued into a global map
$C_1\to \mathfrak X$.
The existence of higher order lifts can be proved similarly. \qed\\

In other words, the obstruction to smooth $\varphi_0$ in fact
 lies in the cohomology group
 $H^1(C_{0}, \varphi_{0}^*\Theta_{\mathfrak X/\Bbb C}/
   \Theta_{C_{0}/O_{0}})$.
The purpose of this section is to calculate this group 
 (precisely speaking, its dual) when the tropical curve 
 $(\Gamma, h)$ is an immersion.
More general situation will be treated in Section \ref{sec:higherval}.

\subsection{Basis of 1-forms with logarithmic poles}\label{subsec:basis}
Let $z$ be an affine coordinate on $\Bbb P^1$ and
 let $a_1, \dots, a_s$ be distinct points on $\Bbb P^1$.
For later use, we fix a basis of the space of meromorphic 1-forms with logarithmic poles allowed
 at $a_1, \dots, a_s$
 as follows.

Assume for simplicity that none of $\{a_i\}$ is $\infty$.
Let $\widetilde{\omega}$ be the sheaf of meromorphic 1-forms allowing 
 logarithmic poles at $\{a_i\}$. 

Then the space of sections $\Gamma(\widetilde{\omega})$ is 
 an $(s-1)$-dimensional vector space spanned by
\[
\tau_1=\frac{dz}{(z-a_s)(z-a_1)},\;\; \cdots,\;\; 
 \tau_{s-1}=\frac{dz}{(z-a_s)(z-a_{s-1})}.
\]
Taking 
\[
\frac{dz}{z-a_1}, \;\; \cdots,\;\;  \frac{dz}{z-a_s}
\]
 as frames of 
 $\widetilde{\omega}$ at neighborhoods of $z = a_1, \dots, a_s$, respectively, 
 the section
 $\tau_j$ takes values 
\[
0,\;\; \dots,\;\; 0,\;\; \frac{1}{a_j-a_s},\;\; 0, \;\;\dots,\;\; 0,\;\; \frac{1}{a_s-a_j},
\]
 at $a_1, \dots, a_s$, respectively.
Here non-zero values are taken at $a_j$ and $a_s$.

In other words, 
 the space of sections $\Gamma(\widetilde{\omega})$
 is identified with the subspace of
\[
\Bbb C^s=\{(v_1, \dots, v_s)\;|\; v_1, \dots, v_s\in \Bbb C\}
\]
 defined by 
\[
\sum_{i=1}^s v_i = 0.
\]

Similarly, the space of sections of 
 $\tilde\omega\otimes\mathcal O_{\Bbb P^1}^{\oplus r}$
 is identified with the subspace of $(\Bbb C^r)^{\oplus s}$
 defined by
 \[
 \sum_{i=1}^s w_i = 0,\;\; w_i\in \Bbb C^r.
 \]
 
%Thus, we can lift $\varphi_{0}$ to $\varphi_k$ for any $k$ in this situation.
%Applying the same type of combinatorics introduced in the proof of 
% Proposition \ref{prop:canonical},
% we can analyze $H^1(C_{0}, \varphi_{0}^*\Theta_{\mathfrak X/\Bbb C}/
%  \Theta_{C_{0}/O_{0}})$.
%As we saw, 
% if $H^1(C_{0}, \varphi_{0}^*\Theta_{\mathfrak X/\Bbb C}/
%  \Theta_{C_{0}/O_{0}})$
%  vanishes, we know that the pre-log curves corresponding 
% to the tropical curve can be smoothed.
%In this section, we give an effective method to calculate 
% $H^1(C_{0}, \varphi_{0}^*\Theta_{\mathfrak X/\Bbb C}/
%  \Theta_{C_{0}/O_{0}})$ (Theorem \ref{thm:obstruction}).
  
\subsection {Calculation of superabundancy 
 for immersive tropical curves via
  algebraic geometry}\label{subsec:main}

%Since, under Assumption A, superabundancy only occurs 
% when the dimension $n$ of the ambient space
% is at least 3, we assume $n\geq 3$ hereafter.
In this subsection, 
 we assume $(\Gamma, h)$
 is an immersive tropical curve in $\Bbb R^n$, $n\geq1$.
Thus, the graphs $\Gamma$ and $\Gamma_{\varphi_0}$
 (see Definition \ref{def:intergraph}) are identical.
In particular, the dual intersection graph $\Gamma_{\varphi_0}$
 of $C_0$ is 3-valent.
We also assume there is no divalent vertex for notational simplicity.

\begin{rem}
	Since the superabundancy only occurs when $n\geq 2$, we may as well
	assume $n\geq 2$.
	Note that although all embedded plane tropical curves are regular in the sense of
	Definition \ref{def:nonsuperabundant2}
	as shown in \cite{M}, when we allow some edges to have the same image, 
	then even immersive 3-valent plane tropical curves can be superabundant.
	See Example \ref{ex:super1}.
	See also \cite[Remark 2.25]{M} for an example with higher valent vertices.
\end{rem}

By Serre duality for nodal curves, we have
\[
H^1(C_{0}, \varphi_{0}^*\Theta_{\mathfrak X/\Bbb C}/
  \Theta_{C_{0}/O_{0}})
  \cong H^0(C_{0}, (\varphi_{0}^*\Theta_{\mathfrak X/\Bbb C}/
  \Theta_{C_{0}/O_{0}})^{\vee}\otimes \omega_{C_0})^{\vee}.
\]
Here $\omega_{C_0}$ is the dualizing sheaf of the nodal curve $C_0$, 
 which is the sheaf of meromorphic 1-forms with logarithmic poles allowed at the nodes.
It is known that $\omega_{C_0}$ is an invertible sheaf.

Let $v$ be a vertex of $\Gamma$ and $C_{0, v}$ be the component
 of $C_0$ corresponding to $v$.
Then since each component of $C_0$ has three special points
 (they are the intersections with the toric divisors of the components of $X_0$.
 Some of them are nodes of $C_0$), 
 $\Theta_{C_{0}/O_{0}}|_{C_{0, v}}\cong \mathcal O(-1)$
  and $\omega_{C_0}|_{C_{0, v}} \cong \mathcal O(-2+s)$, 
  here $s$ is the number of nodes of the component.

From this, it is easy to see that when $s = 1$,  
\[
\Gamma(C_{0, v}, (\varphi_{0}^*\Theta_{\mathfrak X/\Bbb C}/
  \Theta_{C_{0}/O_{0}})^{\vee}\otimes \omega_{C_0})=0
 \] 
 holds.
  
When $s = 2$, we have
\[
 \Gamma(C_{0, v}, (\varphi_{0}^*\Theta_{\mathfrak X/\Bbb C}/
  \Theta_{C_{0}/O_{0}})^{\vee}\otimes \omega_{C_0})\cong
  \Gamma(C_{0, v}, (\varphi_{0}^*\Theta_{\mathfrak X/\Bbb C}/
  \Theta_{C_{0}/O_{0}})^{\vee})
  \]
 on the corresponding component.
Note that there is a following inclusion:
\[
(\varphi_{0}^*\Theta_{\mathfrak X/\Bbb C}/
  \Theta_{C_{0}/O_{0}})^{\vee}\subset 
  (\varphi_{0}^*\Theta_{\mathfrak X/\Bbb C})^{\vee}\cong N^{\vee}_{\Bbb C}
   \otimes \mathcal O_{C_{0}}.
\]

The edges emanating from $v$ 
 span one or two dimensional subspace $V_v$ of $N_{\Bbb C}$.
Then it is clear that 
 $\Gamma(C_{0, v}, (\varphi_{0}^*\Theta_{\mathfrak X/\Bbb C}/
  \Theta_{C_{0}/O_{0}})^{\vee}\otimes\omega_{C_0})$ is given by the subspace
  $V_v^{\perp}\subset N_{\Bbb C}^{\vee}$.
Namely, under the convention in Subsection \ref{subsec:basis},
 an element of 
 $\Gamma(C_{0, v}, (\varphi_{0}^*\Theta_{\mathfrak X/\Bbb C}/
  \Theta_{C_{0}/O_{0}})^{\vee})$ is given by the following data:
\begin{enumerate}
\item[(a)] Give the zero vector of $N_{\Bbb C}^{\vee}$
 to the flag $(v, E_0)$, where $E_0$ is the unique unbounded
 edge emanating from $v$.
\item[(b)] Give $\pm \alpha$, where $\alpha\in V_v^{\perp}$,  to the remaining
 flags associated to $v$.
\end{enumerate}

Let us consider the case $s = 3$.
For simplicity, first we consider the case $n=2$.
Let $C_{0, v}$ and $v$ be as above.
In this case, 
\[
(\varphi_{0}^*\Theta_{\mathfrak X/\Bbb C}|_{C_{0, v}}/
 \Theta_{C_{0}/O_0}|_{C_{0, v}})^{\vee}
 \cong \mathcal O(1)^{\vee}
 \cong \mathcal O(-1)
\]
 holds when $\dim V_v = 2$.
When $\dim V_v = 1$, then it becomes
\[
(\varphi_{0}^*\Theta_{\mathfrak X/\Bbb C_{0}}|_{C_{0, v}}/
 \Theta_{C_{0}/O_0}|_{C_{0, v}})^{\vee}
 \cong \mathcal O.
\]
Therefore, $\Gamma(C_{0, v}, (\varphi_{0}^*\Theta_{\mathfrak X/\Bbb C}/
  \Theta_{C_{0}/O_{0}})^{\vee}\otimes \omega_{C_{0}})$
 is isomorphic to $\Bbb C$ when $\dim V_v = 2$ and
 to $\Bbb C^2$ when $\dim V_v = 1$.
On the other hand, 
\[\begin{array}{ll}
\Gamma(C_{0, v}, (\varphi_{0}^*\Theta_{\mathfrak X/\Bbb C}/
  \Theta_{C_{0}/O_{0}})^{\vee}\otimes \omega_{C_{0}})& \subset 
  \Gamma(C_{0, v}, (\varphi_{0}^*\Theta_{\mathfrak X/\Bbb C})^{\vee}
  \otimes\omega_{C_{0}}) \\
  & \cong N^{\vee}_{\Bbb C}\otimes\Bbb C\langle \tau_1,
   \tau_2\rangle
  \end{array}\]
 on this component.
Here $\tau_1, \tau_2$ are basis vectors of 
 the space of meromorphic 1-forms on 
 $C_{0, v}$ allowing logarithmic poles
 at the three nodes of $C_0$ contained in $C_{0, v}$,
 which we introduced in Subsection \ref{subsec:basis}.

Let 
\[
(la, lb),\;\; (mc, md),\;\; (-la-mc, -lb-md)
\]
 be
 the slopes of the edges of the image $h(\Gamma)$ emanating from $h(v)$.
Here $l, m\in\Bbb Z_{>0}$ are the weights of the relevant edges
 and $(a, b), (c, d)$ are the primitive integral vectors.
Recall that these edges correspond to the intersections 
 between $C_{0, v}$ and the toric divisors
 of the toric surface defined by a fan 
 whose 1-skeleton is given by the tropical curve with 
 the unique vertex $h(v)$ 
 (see \cite[Definition 5.1]{NS}).
\begin{rem}
When $\dim V_v = 1$, then the vectors $(la, lb)$,
 $(mc, md)$ and $(-la-mc, -lb-md)$ are all parallel (up to 
 the orientation), and we can assume $(a, b) = (c, d)$.
Also, in this case the fan defining the toric surface consists of 
 two rays emanating from the origin of $\Bbb R^2$
 (and no two dimensional cone). 
\end{rem}

We take an inhomogeneous  coordinate $z$ on $C_{0, v}$ so that
 the values of $z$ at the points corresponding to the edges of slopes 
  $(la, lb), (mc, md), (-la-mc, -lb-md)$ 
  are 0, 1 and $\infty$, respectively.
As in Subsection \ref{subsec:basis}, 
 we can take $\tau_1, \tau_2$
 and local frames of the sheaf $\omega_{C_0}|_{C_{0, v}}$ at the nodes
 so that 
\[
\tau_1(0) = -\tau_1(\infty) = 1, \;\;
  \tau_1(1) = 0
 \]
  and 
\[
\tau_2(0) = 0,\;\; \tau_2(1)=-\tau_2(\infty) = 1.
\]
\begin{lem}\label{lem:dual1}
Assume $\dim V_v = 2$.
Let $u_1, u_2$ be the generators of $(\Bbb C\cdot(a, b))^{\perp}$,
 $(\Bbb C\cdot(c, d))^{\perp}$
 in $N_{\Bbb C}^{\vee}$ such that
 $u_1((c, d)) = u_2((a, b)) = 1$.
Then the space of sections 
 $\Gamma(C_{0, v}, (\varphi_{0}^*\Theta_{\mathfrak X/\Bbb C}/
  \Theta_{C_{0}/O_{0}})^{\vee}\otimes \omega_{C_{0}})$ 
  is given by the one dimensional subspace
\[
\Bbb C\cdot\langle 
   lu_1\tau_1-mu_2\tau_2\rangle\subset 
   N^{\vee}_{\Bbb C}\otimes \Bbb C\langle \tau_1, \tau_2\rangle.
 \]
In other words, the space of sections 
 $\Gamma(C_{0, v}, (\varphi_{0}^*\Theta_{\mathfrak X/\Bbb C}/
  \Theta_{C_{0}/O_{0}})^{\vee}\otimes \omega_{C_{0}})$ 
  is given by the space
\[
\{(\alpha_1, \alpha_2, \alpha_3)\in \Bbb C\cdot (a, b)^{\perp}\times \Bbb C\cdot (c, d)^{\perp}
 \times \Bbb C\cdot (la+mc, lb+md)^{\perp}\;|\; \alpha_1+\alpha_2+\alpha_3 = 0 \in N^{\vee}_{\Bbb C}\} 
\]
\end{lem}
\proof
The stalks of $\Theta_{C_{0}/O_{0}}$ at 0, 1, $\infty\in C_{0, v}$ are
 spanned by
 $(la, lb), (mc, md)$, and $(-la-mc, -lb-md)$, respectively, considered as 
 subsets of $N_{\Bbb C}\otimes \mathcal O_{C_{0,v}}$.
Sections of $\Gamma(C_{0,v}, (\varphi_{0}^*\Theta_{\mathfrak X/\Bbb C}/
  \Theta_{C_{0}/O_{0}})^{\vee}\otimes \omega_{C_{0}})$
  must annihilate these vectors, and this condition determines the
  mentioned subspace in the statement.\qed\\

Similarly, when $\dim V_v = 1$, we have the following.
\begin{lem}\label{lem:dual2}
Assume $\dim V_v = 1$ and $(a, b) = (c, d)$.
Let $u$ be a generator of $(\Bbb C\cdot (a, b))^{\perp}$ in $N_{\Bbb C}^{\vee}$.
Then the space of sections 
 $\Gamma(C_{0, v}, (\varphi_{0}^*\Theta_{\mathfrak X/\Bbb C}/
  \Theta_{C_{0}/O_{0}})^{\vee}\otimes \omega_{C_{0}})$ 
  is given by the two dimensional subspace
\[
\{(a\tau_1+b\tau_2)u\;|\; a, b\in \Bbb C\}\subset 
 N_{\Bbb C}^{\vee}\otimes \Bbb C\langle\tau_1, \tau_2\rangle.
 \]
Note that again this space is presented as
\[
\{(\alpha_1, \alpha_2, \alpha_3)\in \Bbb C\cdot (a, b)^{\perp}\times \Bbb C\cdot (c, d)^{\perp}
 \times \Bbb C\cdot (la+mc, lb+md)^{\perp}\;|\; \alpha_1+\alpha_2+\alpha_3 = 0 \in N^{\vee}_{\Bbb C}\}. 
\]
\qed
\end{lem}

From this, one sees that in the general case where 
 $n$ is an integer with $n\geq 2$,
 the space $\Gamma(C_{0, v}, (\varphi_{0}^*\Theta_{\mathfrak X/\Bbb C}/
  \Theta_{C_{0}/O_{0}})^{\vee}\otimes \omega_{C_{0}})$
  is described as follows.
\begin{lem}
The space $\Gamma(C_{0, v}, (\varphi_{0}^*\Theta_{\mathfrak X/\Bbb C}/
  \Theta_{C_{0}/O_{0}})^{\vee} \otimes \omega_{C_{0}})$
  is naturally identified with the subspace
\begin{equation}\label{s=3}
 \{(\alpha_1, \alpha_2, \alpha_3)\in \Bbb C\cdot (a, b)^{\perp}\times \Bbb C\cdot (c, d)^{\perp}
  \times \Bbb C\cdot (la+mc, lb+md)^{\perp}\;|\; \alpha_1+\alpha_2+\alpha_3 = 0 \in N^{\vee}_{\Bbb C}\}. 
\end{equation}
 of $N_{\Bbb C}^{\vee}\otimes \Bbb C\langle\tau_1, \tau_2\rangle$.
\end{lem}
\proof
This follows from Lemmas \ref{lem:dual1} and  \ref{lem:dual2}
 and the residue theorem.\qed\\

Using these results, we can combinatorially describe the dual space of obstructions.
Consider a vertex $v$ of the graph $\Gamma$ and let $s$ be the number of bounded edges 
 emanating from it as above.
We attach vectors in $N_{\Bbb C}^{\vee}$ to the flags of $\Gamma$ as follows.
\begin{enumerate}
\item When $s = 1$, then give the zero vector to all the flags whose vertex is $v$.
\item When $s = 2$, then we give the zero vector to the unbounded edge and
 give values $\alpha$ and $-\alpha$ to the remaining flags.
Here $\alpha$ is a vector which belongs to the annihilator subspace of the 
 one or two dimensional plane
 spanned by the direction vectors of the edges emanating from $v$.
\item When $s = 3$, we give values $\alpha_1, \alpha_2, \alpha_3$ to the flags so that
 they satisfy $\alpha_1+\alpha_2+\alpha_3 = 0$ in $N_{\Bbb C}^{\vee}$.
Here $\alpha_i$ belongs to the annihilator subspace of the direction of the
 corresponding edge.
\end{enumerate} 
\begin{defn}\label{def:compatiblelabel}
We say that the set of these values attached to the flags is a \emph{compatible labelling}
 when the sum of the values of the two flags associated to each bounded edge is zero.
\end{defn} 
This 
 reflects the relation of the frames
\[
\frac{dz_1}{z_1}+\frac{dz_2}{z_2} = 0
\] 
 at the corresponding node of $C_0$, 
 here $z_1, z_2$ are coordinates of the two branches at the node.

By the argument so far, we have the following.
\begin{prop}
The set of elements of the dual obstruction group 
 $H^0(C_{0}, (\varphi_{0}^*\Theta_{\mathfrak X/\Bbb C}/
 \Theta_{C_{0}/O_{0}})^{\vee}\otimes \omega_{C_0})$
 has a natural
 one to one correspondence with the set of compatible labellings of the flags of $\Gamma$.\qed
\end{prop}

Now we study the set of compatible labellings.
Let $L = \cup_i L_i$ be the loop part
  of $\Gamma$ (Definition \ref{def:loops}), where $L_i$ 
 are bouquets.
This is a closed subgraph of $\Gamma$.
Let $\Gamma_T = \Gamma\setminus L$.
The closure of a connected component of $\Gamma_T$ is a tree.
There are two types of these trees, namely:
\begin{itemize}
\item[(U)] It contains only one
 flag whose vertex is contained in a loop. 
\item[(B)] Otherwise.
\end{itemize}
By inductive argument, it is easy to see that all the flags in a component
 of type $(U)$
 must have the value zero, including the unique 
 flag whose vertex is contained in a loop.
For the type $(B)$ too, we have the following result.
\begin{lem}\label{lem:bridge}
All the flags of a component of type $(B)$, including the flags whose
 vertices are contained in the loops, must have the value 
 $0\in N_{\Bbb C}^{\vee}$. 
\end{lem}
\proof
Note that $\Gamma$ can be written in the form 
 as in the following figure (Figure \ref{fig:disks}).

\begin{figure}[h]
\includegraphics{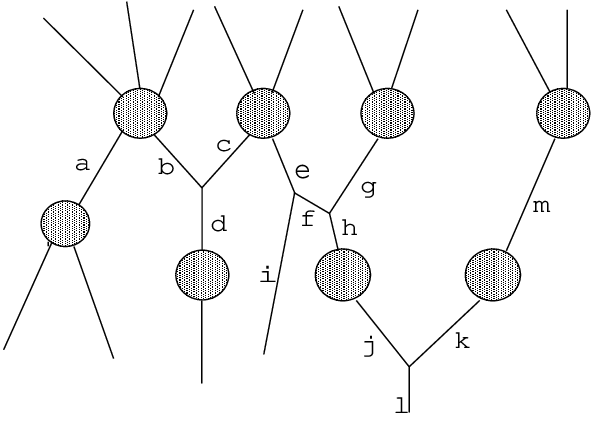}
\caption{}\label{fig:disks}
\end{figure}

Here, each shaded disk corresponds to some component $L_i$ of 
 the loop part of $\Gamma$.
By definition of $\{L_i\}$, if we regard these disks as vertices, we 
 obtain another tree $\Gamma'$.

Recall the above remark that all the edges contained in the components of type 
 $(U)$ have the value zero.
In Figure \ref{fig:disks}, this means that all the edges (outside the shaded disks)
 except the ones labeled by $a, b, c, \cdots, m$ have the value zero.
We call the edges outside the shaded disks the \emph{bridges}.

Now, by the fact that $\Gamma'$ is a tree,
 it is easy to see that there is a shaded disk such that the bridges emanating from
 it have the value zero except one bridge.
Let us call this bridge $r$ and call the remaining bridges 
 emanating from the shaded disc $a_1, \dots, a_s$.
By the condition that the sum of the
 values of the three edges emanating from each vertex 
 is
 zero, we see that the sum of the values attached to $a_1, \dots, a_s, r$
 is zero.
Since $a_1, \dots, a_s$ have the
 value zero, it follows that $r$ has also the value zero. 
By induction on the number of shaded disks, we see that all the bridges,
 and so all the edges of the components of  
 type $(B)$ also have the value zero.\qed\\

According to this lemma, 
 we only need to consider the flags contained in some bouquet.
Let $L_i$ be a bouquet.
This is a graph with 2-valent and 3-valent vertices.
Every 3-valent vertex $v$ of $L_i$ determines a one or two dimensional 
 subspace of $N_{\Bbb C}$ spanned by the edges emanating from it.
We write it by $V_v$ as before.
Also, every edge $E$ of $L_i$ determines a one dimensional subspace of 
 $N_{\Bbb C}$.
%We write it by $W_e$.

Let us describe the space
 $H = H^1(C_{0}, \varphi_{0}^*\Theta_{\mathfrak X/\Bbb C}/
  \Theta_{C_{0}/O_{0}})^{\vee}$.
Let $\{ v_i\}$ be the set of 3-valent vertices of $L$.
Cutting $L$ at each $v_i$, we obtain
 a set of piecewise linear segments $\{l_m\}$ in $h(\Gamma)$. 
Let $U_m$ be the linear subspace of $N_{\Bbb C}$
 spanned by the direction vectors of the segments of $l_m$.
The following theorem follows from the argument so far.
\begin{thm}\label{thm:obstruction}
Let $(\Gamma, h)$ be a 3-valent immersive 
 tropical curve.
Elements of the space $H$ are described by the following data.
\begin{enumerate}[(I)]
\item Give the zero vector of $(N_{\Bbb C})^{\vee}$ to all the flags not contained in $L$.
\item Give a value
 in $(U_m)^{\perp}\subset (N_{\Bbb C})^{\vee}$
 to each of the flags  associated to the edges of $l_m$.
\item The data in $(I)$ and $(II)$ give an element of $H$ if and only if the following 
 conditions are satisfied.
\begin{enumerate}
\item At each vertex
 $v$ of $\Gamma$, 
\[
u_1+u_2+u_3=0
\]
 holds as an element of $(N_{\Bbb C})^{\vee}$.
Here $u_1, u_2, u_3$ are the data attached to the three
 flags in $\Gamma$
 which have $v$ as the vertex.
\item The data in $(II)$
 is compatible on each edge of $l_m$, 
 in the sense that the sum of the values attached to the two flags
 of any edge of $l_m$ is zero. \qed
\end{enumerate}
\end{enumerate} 
\end{thm}
\begin{thm}\label{thm:dim}
Let $(\Gamma, h)$ be a 3-valent immersive tropical curve.
Then the number of parameters to deform it is given by
\[
(n-3)(1-g)+e + \dim H,
\]
 here $H$ is the space in Theorem \ref{thm:obstruction}.
\end{thm}
\proof
Let us write the dimension of  
$H^0(C_{0}, \varphi^*_{0}\Theta_{\mathfrak X/\Bbb C}/\Theta_{C_{0}/O_{0}})$
by $d_1$ and
the dimension of 
$H^1(C_{0}, \varphi^*_{0}\Theta_{\mathfrak X/\Bbb C}/\Theta_{C_{0}/O_{0}})$
by $d_2$.
In the proof of \cite[Proposition 5.7]{NS} (see also Subsection \ref{subsec:corrthm}
 and Remark \ref{rem:degdeform}), 
it is shown that when the tropical curve $(\Gamma, h)$ is immersive, 
 $d_1$ 
 is the same as the dimension of the parameter
 space of the corresponding
 tropical curve.
Also, it is known that the dimension of the cohomology group
 $H^0(C_0, \omega_0)$ equals $g$ (see \cite[Exercise 3-4]{HM}).
Therefore, the dimension of the cohomology group
 $H^1(C_0, \varphi_0^*\Theta_{X_0})\cong H^0(C_0, \omega_{C_0}^{\oplus n})^{\vee}$
 equals $ng$.

Then by the long exact sequence associated to the sequence (\ref{eq:deform})
 after Proposition \ref{prop:nonobst}, 
 we have the equality
\[
d_1-d_2 = (n-3)(1-g)+e.
\] 
Note that $\dim H = d_2$.
These observations prove the theorem. \qed 

\begin{rem}\label{rem:thm40add}
An important remark is that 
 the sheaves $\varphi_0^*\Theta_{\mathfrak X/\Bbb C}$
 and $\varphi_0^*\Theta_{\mathfrak X/\Bbb C}/\Theta_{C_0/O_0}$
 exist even when the map $\varphi_0$ does not exist globally.
They can be constructed by gluing local pieces on each component, 
 and these local pieces are determined by the combinatorial data 
 of the tropical curve $(\Gamma, h)$ around each vertex corresponding to 
 the components of $C_0$.
In particular, the space $H$ also makes sense and
 Theorem \ref{thm:dim} gives us the number of freedom to deform the given
 immersive  
 superabundant tropical curve even if a pre-log curve
 of type $(\Gamma, h)$ does not exist.
\end{rem}
\begin{rem}
We see that it almost suffices to check the conditions only at the 3-valent vertices of
 $L$.
The conditions (I) and (III) ($a$) imply that at a 2-valent vertex of $L$, 
 the values $u, u'$ associated to the relevant two flags satisfy $u+u' = 0$.
Together with the condition (III) (b), we see that
  on each $l_m$, the values associated to the flags are unique up to sign.
\end{rem}
The following is immediate from this,
 because when the genus of $\Gamma$ is one,
 there is no 3-valent vertex in $L$.
See also \cite{S}.
\begin{cor}
When $(\Gamma, h)$ is an immersive 3-valent tropical curve of genus one, 
 then $H \cong U^{\perp}$,
 here $U$ is the linear subspace of $N_{\Bbb C}$
 spanned by the direction vectors of the segments of the loop of $\Gamma$.\qed
\end{cor}

\begin{cor}\label{cor:genusonepre-logexist}
When the tropical curve $(\Gamma, h)$ is of genus one and immersive, 
 a pre-log curve of type $(\Gamma, h)$ exists even if it is 
 superabundant.
\end{cor}
\proof	
Let $A$ be the minimal dimensional affine subplane of $N_{\Bbb R}$
 containing the loop part of $h(\Gamma)$.
Then the loop part of $h(\Gamma)$ defines a tropical curve $(\bar\Gamma, \bar h)$
 of genus one in $A$
 in a natural way, and it is non-superabundant 
 (or even regular in the sense of Definition \ref{def:nonsuperabundant2}).
Then there is a pre-log curve of type $(\bar\Gamma, \bar h)$
 by Proposition \ref{prop:exist} (see also the last part of Section \ref{sec:pre}).
Since $(\Gamma, h)$ is constructed from $(\bar\Gamma, \bar h)$ by 
 grafting trees, it is easy to see that we can extend a pre-log curve of type $(\bar\Gamma, \bar h)$
 to a pre-log curve of type $(\Gamma, h)$.\qed\\

\noindent
\subsection{{\bf Examples}}\label{subsec:ex}
\begin{example}
Let us consider immersive tropical curves $\Gamma_1$ and $\Gamma_2$
 of genus two
 in $\Bbb R^3$ given in Figure \ref{fig:supex}.

\begin{figure}[h]
\includegraphics{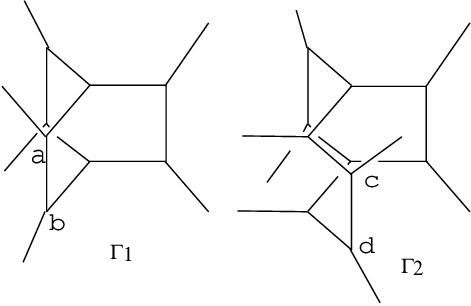}
\caption{}\label{fig:supex}
\end{figure}

The curve $\Gamma_1$ has 6 unbounded edges of the directions
\[
(1, 0, 1),\;\; (1, 0, -1),\;\; (-1, -1, 1),\;\; (-1, -1, -1),\;\; (0, 1, 1),\;\; (0, 1, -1).
\]
The bounded edges are:
\begin{itemize}
\item Three parallel vertical edges of the direction $(0, 0, 1)$.
\item Three pairs of parallel edges of the directions
\[
(1, 0, 0),\;\; (-1, -1, 0), \;\; (0, 1, 0),
\]
 respectively.
\end{itemize}
The curve $\Gamma_2$ is a modification of $\Gamma_1$
 at the vertices $a$ and $b$.
Namely:
\begin{enumerate}
\item Delete the edge $\overline{ab}$ as well as the neighboring unbounded edges.
\item Add a pair of parallel unbounded edges of the direction $(-1, 0, 0)$,
 and a pair of parallel bounded edges of the direction 
 $(1, 1, 0)$ of the same length.
\item Connect the end points $c, d$ of the bounded edges added in (2) by
 a segment of the direction $(0, 0, 1)$.
\item Add unbounded edges of the directions $(1, 1, 1), (1, 1, -1)$ at
 the vertices $c, d$, respectively.
\end{enumerate}
Using Theorem \ref{thm:obstruction}, it is easy to see that $\Gamma_1$
 is superabundant, while $\Gamma_2$ is non-superabundant.
Namely, the sets of piecewise linear segments $\{l_m\}$ of
 these tropical curves are given by
 the following three components, respectively (Figure \ref{fig:decomp}).
 
\begin{figure}[h]
\includegraphics{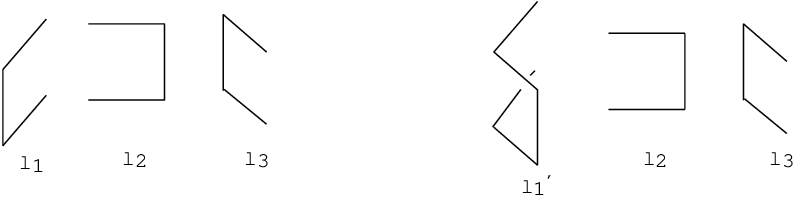}
\caption{}\label{fig:decomp}
\end{figure}

We write the corresponding linear subspaces of $N_{\Bbb C}\cong \Bbb C^3$
 by $U_{l_1}, U_{l_2}$, etc..
Then, using standard nondegenerate 
 quadratic form on $\Bbb C^3$ to identify it with its dual, 
\[
(U_{l_1})^{\perp}\cong \Bbb C\cdot (1, 0, 0),\;\;
(U_{l_2})^{\perp}\cong \Bbb C\cdot (0, 1, 0),\;\;
(U_{l_3})^{\perp}\cong \Bbb C\cdot (1, -1, 0).
\]
Then it is easy to see that the space $H$ 
 for $\Gamma_1$ is a one dimensional vector space.
Thus, $\Gamma_1$ is superabundant.

On the other hand, since $U_{l'_1}\cong \Bbb C^3$, 
 $(U_{l'_1})^{\perp} = \{0\}$.
From this, it is easy to see that the space $H$ for $\Gamma_2$ is $\{0\}$. 
Therefore, $\Gamma_2$ is non-superabundant.
\end{example}

\begin{example}\label{ex:super1}
Plane tropical curves were studied by Mikhalkin \cite{M} in great detail.
There it was shown that any conventional (that is, immersive and
 without multiple edges) plane tropical curve is non-superabundant
 and smoothable.
On the other hand, if we only assume Assumption A, 
 then even 3-valent immersive (in the sense of Definition \ref{immersive}) 
 plane tropical curves
 can be superabundant and non-smoothable.
The simplest example is given by the one in Figure \ref{fig:0G}.

\begin{figure}[h]
\includegraphics{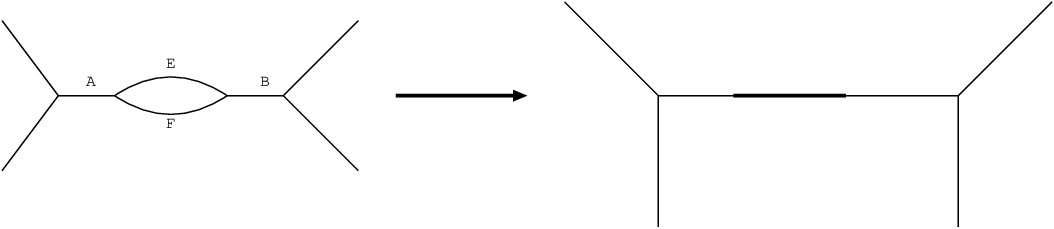}
\caption{The abstract graph (the picture on the left) is mapped into $\Bbb R^2$.
The edges $E$ and $F$ have the same image
 (the bold line in the picture on the right).}\label{fig:0G}
\end{figure}

In Figure \ref{fig:0G}, all the edges except $E$ and $F$ have weight two, 
 while the edges $E$ and $F$ have weight one.
The loop part is the union of $E$ and $F$, and the space $H$ is one dimensional.
According to the study in \cite{N3}, this tropical curve is smoothable 
 if and only if the lengths of the images of the edges $A$ and $B$ are the same.

\end{example}

\section{Algebraic curves corresponding to higher valent vertices}
\label{sec:higherval}
Theorem \ref{thm:obstruction} describes the superabundancy of 
 immersive tropical curves in a satisfactory way.
However, from the point of view of the general study of 
 superabundant tropical curves, it is yet not satisfactory.
The reason is that in the study of superabundant tropical curves
 (and associated algebraic curves), non-immersive tropical curves
 appear in typical situations.

In this section,
 we generalize Theorem \ref{thm:obstruction}
 to describe the dual obstruction cohomology group
 associated to algebraic curves corresponding to 
 a tropical curve satisfying Assumption A (in particular, it need not
 be immersive), see Theorem \ref{thm:obstruction2}.

Contrary to the immersive case (Theorem \ref{thm:obstruction}), the description is not purely 
 combinatorial reflecting the fact that higher valent vertices correspond to rational curves with
 $k (\geq 4)$ marked points, which have their own moduli.
Nevertheless, we will see later that the combinatorics of the tropical curve still have
 crucial data of the corresponding holomorphic curves.
This study owes a lot to the degeneration technique.
Namely, given a curve we want to study, the structure of its degeneration is not necessarily simple. 
However, in the degenerate situation,
 we can reduce various calculation to the standard case 
 (typically curves of degree one in projective spaces).
Thus, we can obtain important quantities such as obstruction classes rather explicitly.

This study reveals the importance of
 higher valent vertices in the study of superabundant 
 tropical curves in various ways.
In Subsection \ref{subsec:highval}, we prove a new criterion for the smoothability
 of superabundant tropical curves of genus one (Theorem \ref{thm:loop})
 which is totally different from the known criterion (well-spacedness condition
 \cite{S, N3}).
Also, we study a superabundant tropical curve of genus two in which
 a higher valent vertex plays a crucial role in the computation of the 
 obstruction (Example \ref{ex:3}).

\subsection{Computation of the dual obstruction cohomology groups associated to  higher valent vertices}\label{subsec:localobst}
Let $\Gamma_0$ be a 3-valent tree graph with unbounded edges 
 (see the beginning of Subsection \ref{subsec:pre}).
Let $r$ be the number of vertices of $\Gamma_0$.
Then $\Gamma_0$ has $r+2$ unbounded edges.
We assume all the edge weights of $\Gamma_0$ are 1.
We write by
\[
E_1, \dots, E_{r+2}
\]
 the unbounded edges.
Let 
\[
h_0\colon \Gamma_0\to \Bbb R^{r+1}
\]
 be the map which contracts all the bounded edges 
 and maps the edges $E_i$ onto the 
 rays as follows:
\[
E_i\to \Bbb R_{\geq 0}\cdot (0, \dots, \stackrel{i}{\breve{1}}, \dots, 0),\;\;
 i = 1, \dots, r+1,
\]
\[
E_{r+2}\to \Bbb R_{\geq 0}\cdot (-1, \dots, -1).
\]
In particular, $h_0$ maps all the vertices of $\Gamma_0$ to the origin.
Clearly this satisfies the balancing condition and gives $\Gamma_0$
 a structure of a tropical curve.

Let $(\Gamma, h)$, $h\colon \Gamma\to \Bbb R^{n}$,
 be a tropical curve which contains a subgraph 
 such that the restriction of $h$ to it is isomorphic to 
 $(\Gamma_0, h_0)$ restricted to a suitable open subset containing all the vertices.
We regard this subset of $h_0(\Gamma_0)$ as a subgraph of $h(\Gamma)$.
We assume $(\Gamma, h)$ satisfies Assumption A.
Let $X$ be a toric variety associated to $(\Gamma, h)$
 and $\mathfrak X$ be a toric degeneration of $X$ defined respecting 
 $(\Gamma, h)$.
Let $\varphi_0\colon C_0\to X_0$ be a pre-log curve of type $(\Gamma, h)$.
We put a log structure on $C_0$ so that the 
 map $\varphi_0$ extends to a morphism between log schemes
 (with the log structure on $\mathfrak X$ coming from the toric structure)
 and the composition of $\varphi_0$ and the projection 
 $\mathfrak X\to\Bbb C$
 is log smooth.
\begin{rem}
As we mentioned several times,
 in general, there exist tropical curves such that there are 
 no degenerate holomorphic
 curves of those types, even if we restrict our attention to embedded
 tropical curves.
Here we assume there is a 
 degenerate holomorphic curve of type $(\Gamma, h)$.
\end{rem}
Our purpose is to calculate the obstruction cohomology group
 $H^1(C_0, \mathcal N_{\varphi_0})$ or its dual 
 $H^0(C_0, \mathcal N^{\vee}_{\varphi_0}\otimes\omega_{C_0})$
 as in the previous section.
Here $\mathcal N_{\varphi_0} = \varphi_0^*\Theta_{\mathfrak X/\Bbb C}/\Theta_{C_0/O_0}$
 is the log normal sheaf.
As before, the sheaf $\mathcal N_{\varphi_0}$
 is locally free, and elements of 
 $H^0(C_0, \mathcal N^{\vee}_{\varphi_0}\otimes\omega_{C_0})$ can be
 described by gluing the sections of 
 $\mathcal N^{\vee}_{\varphi_0}
 \otimes\omega_{C_0}$ restricted to components of $C_0$.
Therefore, we first concentrate on the study of the restriction of 
 $\mathcal N_{\varphi_0}^{\vee}\otimes\omega_{C_0}$ to the component
 $C_{0, v}$
 of $C_0$ corresponding to the vertex $v$ of $h(\Gamma)$
 which is modelled on the unique vertex of $h_0(\Gamma_0)$.
It is easy to see the following.
\begin{lem}\label{lem:highnormal}
The restriction of $\mathcal N_{\varphi_0}$ to $C_{0, v}$  is isomorphic
 to 
\[
\mathcal O_{\Bbb P^1}(1)^{r}\oplus \mathcal O_{\Bbb P^1}^{n-r-1}. 
\]\qed
\end{lem}
Let $s$ be the number of nodes of $C_0$ contained in $C_{0, v}$.
This is the same as the number of the 
 edges among $h_0(E_1), \dots, h_0(E_{r+2})$
 which are the restrictions of the bounded edges
 of $h(\Gamma)$.
\begin{lem}\label{lem:highdim}
\[
\dim H^0(C_{0, v}, \mathcal N^{\vee}_{\varphi_0}
 \otimes\omega_{C_0}|_{C_{0, v}})
 = \begin{cases}
  0,\;\; (s = 0, 1)\\
 r(s-2)+ (n-r-1)(s-1),\;\; (s\geq 2)
\end{cases}
\]
\end{lem}
\proof
Note that the restriction  
 $\omega_{C_{0, v}} = \omega_{C_0}|_{C_{0, v}}$ is isomorphic to 
 $\mathcal O_{\Bbb P^1}(-2+s)$.
Thus, 
\[
\begin{array}{ll}
H^0(C_{0, v}, \mathcal N^{\vee}_{\varphi_0}\otimes
 \omega_{C_0}|_{C_{0, v}})
 & \cong 
 H^0(C_{0, v}, 
 (\mathcal O_{\Bbb P^1}(-1)^{r}\oplus \mathcal O_{\Bbb P^1}^{n-r-1})
  \otimes \mathcal O_{\Bbb P^1}(-2+s))\\
 & \cong
 H^0(C_{0, v}, \mathcal O_{\Bbb P^1}(-3+s)^r\oplus
  \mathcal O_{\Bbb P^1}(-2+s)^{n-r-1}).
\end{array}
\]
The result follows from this.\qed\\

First, we assume all the edges $E_1, \dots, E_{r+2}$ are
 the restrictions of the bounded edges of $\Gamma$, so that
 $s = r+2$ (see Remark \ref{rem:unbounded} for general cases).
The sheaf $\omega_{C_{0, v}}$ is the sheaf of meromorphic 
 1-forms such that they can have logarithmic poles at the points of 
 $C_{0, v}$ which are mapped to the toric divisors of the corresponding
 component $X_{0, v}$ of $X_0$.
By the residue theorem,
 the residues $a_1, \dots, a_{s}$ at these poles
 sum up to 0:
\[
a_1+\cdots + a_{s} = 0.
\]
Let us fix an affine coordinate $\zeta$ on $C_{0, v}$
 such that the coordinates of the points corresponding to the edges
 $E_1, \dots, E_{r+1}, E_{r+2}$ are
\[
p_1, \dots, p_{r+1}, \infty.
\]
Then we can take
\[
\sigma_i = \frac{d\zeta}{\zeta - p_i},\;\; i = 1, \dots, r+1
\]
 as the basis of the space of the sections of $\omega_{C_{0, v}}$.

Recall that the log normal sheaf $\mathcal N_{\varphi_0}$
 is the quotient
 $\varphi_0^*\Theta_{\mathfrak X/\Bbb C}/
 \Theta_{C_0/O_0}$ and the sheaf
 $\varphi_0^*\Theta_{\mathfrak X/\Bbb C}$ is naturally isomorphic to 
 the sheaf $N\otimes_{\Bbb Z}\mathcal O_{C_0}$.
Here $N$ is the free abelian group of rank $n$.
Then the sheaf  
 $\mathcal N^{\vee}_{\varphi_0}\otimes
 \omega_{C_0}|_{C_{0, v}}$
 can be seen as a subsheaf of
 the sheaf of $N_{\Bbb C}^{\vee}$-valued meromorphic 1-forms.
In particular, a section of the restriction of 
 $\mathcal N^{\vee}_{\varphi_0}\otimes
 \omega_{C_0}$
 to the component $C_{0, v}$ can be written in the form
\begin{equation}\label{eq:localobst}
\sum_{i, j = 1}^{r+1}a_{i, j}e_i^{\vee}\otimes\sigma_j
 + \sum_{i = r+2}^n\sum_{j = 1}^{r+1}c_{i, j}e_i^{\vee}\otimes\sigma_j,
\end{equation}
 where $a_{i, j}$ and $c_{i, j}$ are complex numbers and 
 $\{e_i\}$, $i = 1, \dots, n$ is a basis of $N$ such that 
 $\{e_i\}$, $i = 1, \dots, r+1$ is a basis of the sublattice of $N$
 obtained as the intersection of $N$ with the subspace parallel to the
 affine subspace
 of $N_{\Bbb R}$ spanned by the edges emanating from the vertex $v$.
Note that $\{e_i\}$, $i = 1, \dots, r+1$ can be taken so that $e_i$
 is the direction vectors of $h(E_i)$.
The set of vectors $\{e_i^{\vee}\}$ is the dual basis of $\{e_i\}$.
 
\begin{assum}
In this subsection, hereafter we assume $n = r+1$ for notational simplicity, since
 the directions $e_i^{\vee}$, $i\geq r+2$ play rather trivial role. 
\end{assum}

By the argument as in Section  \ref{sec:dual obstruction},
 it is easy to see that
 at the pole corresponding to the edge $E_i$, 
 the vector valued residue of the above section
 should take a value in the subspace $e_i^{\perp}$
 of $N_{\Bbb C}^{\vee}$.
This implies the following.
\begin{lem}\label{lem:cond1}
The coefficients of 
 a section $\sum_{i, j = 1}^{r+1}a_{i, j}e_i^{\vee}\otimes\sigma_j$
 of $\mathcal N^{\vee}_{\varphi_0}\otimes
 \omega_{C_0}|_{C_{0, v}}$
 satisfy
\[
a_{1,1} = a_{2, 2} = \cdots = a_{r+1, r+1} = 0.
\]\qed
\end{lem} 
A general torically transverse curve of degree one in $\Bbb P^{r+1}$ 
 is defined by the equations of the following form:
\[
b_2X_1+X_2+c_2  = 0,\;\; b_3X_1+X_3+c_3 = 0,\;\;
 \cdots, \;\; b_{r+1}X_1+X_{r+1}+c_{r+1} = 0,
\]
 where $X_i = \frac{x_i}{x_{r+2}}$
 are affine coordinates of $\Bbb P^{r+1}$
 corresponding to $e_i^{\vee}$,
 and $x_i$ are homogeneous coordinates of $\Bbb P^{r+1}$.

When $X_2X_3\cdots X_{r+1}\neq 0$,  
 a log tangent vector of this line in $\Bbb P^{r+1}$ is given by
\[
X_1\partial_{X_1} + \frac{X_2+c_2}{X_2}X_2\partial_{X_2}
 + \cdots + \frac{X_{r+1}+c_{r+1}}{X_{r+1}}X_{r+1}\partial_{X_{r+1}}.
\]
A section of $\mathcal N^{\vee}_{\varphi_0}\otimes
 \omega_{C_0}|_{C_{0, v}}$ should annihilate this log tangent vector.
Note that the vector $e_i^{\vee}$ is naturally identified with the log 
 cotangent vector $\frac{dX_i}{X_i}$.
 
Let $\zeta$ be the parameter of $C_{0, v}$ given by 
 $\varphi_0^*X_1$. 
Taking $\zeta$ in this way,
 the coordinate $p_i$ of the pole of $\sigma_i$
 is given by
\[
p_i = -\frac{c_i}{b_i},\;\; i = 1, \dots, r+1,
\]
 here we take $c_1 = 0, b_1 = -1$ (thus, $p_1 = 0$).
%Thus, the 1-form $\sigma_i$ is written as
%\[
% \sigma_i = \frac{d\zeta}{\zeta+\frac{c_i}{b_i}},\;\; i = 1, \dots, r+1.
%\]

Then the above annihilating condition implies
\[
\sum_{i, j = 1}^{r+1}a_{i, j}\frac{\zeta}{\zeta-p_i}\cdot
  \frac{1}{\zeta-p_j} 
  =
 \zeta \sum_{i, j = 1}^{r+1}a_{i, j}
   \frac{\prod_{l\neq i, j}(\zeta-p_l)}{\prod_{l=1}^{r+1}(\zeta-p_l)}= 0.
\]
Thus, we have the following.
Let us define a polynomial $P(\zeta)$ of $\zeta$ by
\[
P(\zeta) = \sum_{i, j = 1}^{r+1}a_{i, j}\prod_{l\neq i, j}(\zeta-p_l)
 = \sum_{k = 0}^{r-1}A_k(a_{i, j})\zeta^k,
\]
 where 
\[
A_k(a_{i, j}) = (-1)^{r-1-k}\sum_{i, j = 1}^{r+1}
 \sum_{\substack{J\subset I\setminus\{i, j\},\\ |J| = r-1-k}}
 p_{j_1}\cdots p_{j_{r-1-k}}a_{i, j}
\]
 is a linear polynomial of $a_{i,j}$. 
Here in the summation the set $J = \{j_1, \dots, j_{r-1-k}\}$
 runs over all subsets of $I\setminus\{i, j\}$ of cardinality $r-1-k$.
\begin{lem}\label{lem:cond2}
The coefficients $a_{i, j}$ of 
 a section $\sum_{i, j = 1}^{r+1}a_{i, j}e_i^{\vee}\otimes\sigma_j$
 of $\mathcal N^{\vee}_{\varphi_0}\otimes
 \omega_{C_0}|_{C_{0, v}}$
 satisfy
\[
A_k(a_{i, j}) = 0,\;\; k = 0, \dots, r-1.
\]\qed
\end{lem}
\begin{rem}
Note that Lemmas \ref{lem:cond1} and \ref{lem:cond2}
 give $2r+1$ linear conditions to the coefficients $\{a_{i,j}\}$.
These are all the conditions imposed on these coefficients,
 giving $(r+1)^2-(2r+1) = r^2$ freedom to 
 them.
This is compatible with the dimension calculated in 
 Lemma \ref{lem:highdim}
 (with $s = r+2, n = r+1$).
\end{rem} 
\begin{rem}\label{rem:unbounded}
When there are unbounded edges of $\Gamma_0$ which are also
 unbounded in $\Gamma$, 
 then the numbers $a_{i, j}$ associated to these edges should be set
 to zero.
Explicitly, when the edge of direction $e_j$ is also unbounded in 
 $\Gamma$, then $a_{i, j}$, $i  = 1, \dots, r+1$ should be all zero.
\end{rem}

\subsection{Dual obstruction spaces for general vertices}\label{subsec:gendual}
So far we considered tropical curves with a higher valent vertex whose
 image is locally isomorphic to the image of the tropical curve
 $(\Gamma_0, h_0)$ introduced at the beginning of Subsection 
 \ref{subsec:localobst}.
Such a vertex corresponds to a line in a projective space.
In this subsection we consider more general higher valent vertices.

Namely, let us take a 3-valent
 abstract tree $\Gamma_1$ with $r+2$ unbounded edges
 so that it is combinatorially isomorphic to $\Gamma_0$.
However, now we allow the edges of $\Gamma_1$ to have 
 nontrivial weights (both on bounded and unbounded edges).
Then consider a proper map 
\[
h_1\colon \Gamma_1 \to \Bbb R^n = N_1\otimes\Bbb R,
\]
 where $N_1$ is a free abelian group of rank $n(\leq r+1)$,  with the following properties:
\begin{itemize}
\item The map $h_1$ contracts all the bounded edges of $\Gamma_1$.
\item The map $h_1$ gives 
 a structure of a tropical curve to $\Gamma_1$.
\item The image of $h_1$ is not contained 
 in a proper affine hypersurface of $\Bbb R^n$.
\end{itemize}
\begin{rem}
For applications, it is important to note that the map $h_1$ may send
 some of the unbounded edges to the same image.
Although the image of these edges is a half line which corresponds to 
 a toric divisor through the construction in Subsection \ref{subsec:pre-logdef},
 a pre-log curve of type $(\Gamma_1, h_1)$
 should intersect this toric divisor at several different points
 of the domain curve (their images on the toric divisor can be the same),
 contrary to the cases of conventional tropical curves.
Also, the weight of the image is represented by a tuple of
 positive integers.
See \cite{N3}.
\end{rem}

Let us take a standard basis $\{e_1, \dots, e_{r+1}\}$ of 
 $\Bbb R^{r+1}$.
Recall that the directions of the edges of the image $h_0(\Gamma_0)$
 are given by $e_1, \dots, e_{r+1}$ and $f = -e_1-\cdots - e_{r+1}$.
We write by $E_1, \dots, E_{r+1}$ and $F$ the unbounded edges
 of $\Gamma_0$ which are mapped to the edges of these directions
 by $h_0$.

We fix the natural identification $\iota$
 of the abstract graphs $\Gamma_0$ and 
 $\Gamma_1$.
Let 
\[
n_1, \dots, n_{r+1}\in N_1
\]
 be the primitive integral vectors of 
 the directions of $h_1\circ\iota(E_1), \dots, h_1\circ\iota(E_{r+1})$.
Let $w_1, \dots, w_{r+1}$ be the weights of these edges.

Then define a linear map 
\[
\Xi\colon \Bbb R^{r+1}\to \Bbb R^n
\]
 by extending the map
\[
e_i\mapsto w_in_i,\;\; i = 1, \dots, r+1
\]
 linearly.
It is easy to see that $h_0(\Gamma_0)$ is mapped by $\Xi$
 onto $h_1(\Gamma_1)$.

The images of the tropical curves $(\Gamma_0, h_0)$
 and $(\Gamma_1, h_1)$ can be seen 
 as non-complete fans consisting of one dimensional cones.
Let $\Bbb P_{(\Gamma_0, h_0)}$ and
 $\Bbb P_{(\Gamma_1, h_1)}$
  be the toric varieties defined by these fans. 
The map $\Xi$ induces a map 
$f_{\Xi}$ between toric 
varieties from $\Bbb P_{(\Gamma_0, h_0)}$ 
to 
$\Bbb P_{(\Gamma_1, h_1)}$.
The map $f_{\Xi}$ sends pre-log curve of type $(\Gamma_0, h_0)$
 to pre-log curves of type $(\Gamma_1, h_1)$, 
 and all pre-log curves of type $(\Gamma_1, h_1)$
 can be obtained in this way.

Let 
\[
\psi: \Bbb P^1\to \Bbb P_{(\Gamma_1, h_1)}
\]
 be a curve of type $(\Gamma_1, h_1)$.
Assume that $\psi$ is obtained from 
 a line $\varphi\colon\Bbb P^1\to \Bbb P^{r+1}$
 by composing with $f_{\Xi}$ (note that $\Bbb P_{(\Gamma_0, h_0)}$
 is naturally an open subvariety of $\Bbb P^{r+1}$).
Its logarithmic conormal sheaf $\mathcal N^{\vee}_{\psi}$
 can be pulled back by the map $f_{\Xi}$.
Note that a section of the sheaf
 $\mathcal N_{\psi}^{\vee}\otimes\omega_{\Bbb P^1}$
 can be seen as a 1-form 
 on $\Bbb P^1$ with values in $(N_1\otimes \Bbb C)^{\vee}$.

The map $\Xi$ gives a linear 
 map $\Bbb Z^{r+1}\to N_1$, 
 and it in turn induces a map between the dual spaces
\[
\Xi^*\colon (N_1\otimes\Bbb C)^{\vee}\to (\Bbb C^{r+1})^{\vee}.
\] 
 
Then the pull back 
 $f_{\Xi}^{\ast}\mathcal N_{\psi}^{\vee}$
 can be seen as a sheaf on $\Bbb P^1$ with values in 
 $(\Bbb C^{r+1})^{\vee}$.
Since the sheaf 
 $f_{\Xi}^{\ast}\mathcal N_{\psi}^{\vee}$
 is naturally a subsheaf of $\mathcal N_{\varphi}^{\vee}$, 
 sections of $f_{\Xi}^{\ast}
 \mathcal N_{\psi}^{\vee}\otimes\omega_{\Bbb P^1}$
 can also be described by the numbers $\{a_{i, j}\}$
 given in Equation (\ref{eq:localobst}) in Subsection \ref{subsec:localobst}.

Let $y_i$ be the point of $\Bbb P^1$ which is mapped by $\varphi$
 to the toric divisor of $\Bbb P^{r+1}$ 
 corresponding to the edge $E_i$ of $\Gamma_0$ 
 ($i = 1, \dots, r+1$).
 
The fiber of $\mathcal N_{\varphi}^{\vee}$ at each $y_i$
 can be identified with the annihilator subspace $(E_i)^{\perp}$
 in 
 $(\Bbb C^{r+1})^{\vee}$ of the direction of the edge $E_i$.
\begin{defn}
Let 
 $F_i$ be the subspace of $(E_i)^{\perp}$ of codimension $r+1-n$
 which annihilates the kernel of the map
\[
\Xi_{\Bbb C}\colon \Bbb Z^{r+1}\otimes\Bbb C\to N_1\otimes\Bbb C. 
\]
\end{defn}
It is clear that the fiber of $f_{\Xi}^{\ast}\mathcal N_{\psi}^{\vee}$ at $y_i$
is canonically isomorphic to $F_i$. 
\begin{lem}\label{lem:genhigh}
A section of 
 $f_{\Xi}^{\ast}
  \mathcal N_{\psi}^{\vee}\otimes\omega_{\Bbb P^1}$
 is described by the set of numbers $\{a_{i, j}\}$
 with the following properties.
\begin{itemize}
\item The vector valued residue
\[
\sum_{i=1}^{r+1}a_{i, j}e_i^{\vee}
\]
 at $y_j$ is contained in $F_j$.
\item The conditions given in Lemmas \ref{lem:cond1} and 
 \ref{lem:cond2}, namely
\[
a_{1, 1} = \dots = a_{r+1, r+1} = 0
\]
 and
\[
A_k(a_{i, j}) = 0, \;\; k = 0, \dots, r-1
\]
 hold.
\end{itemize}
\end{lem} 
\proof
It is clear that the conditions for $\{a_{i, j}\}$
 in the statement are necessary to 
 define a section of 
 $f_{\Xi}^{\ast}\mathcal N_{\psi}^{\vee}\otimes
  \omega_{\Bbb P^1}$.

For a general point $z$ on $\Bbb P^1$, the fiber of 
 $\mathcal N_{\varphi}^{\vee}$ is the annihilator of the (log)
 tangent space of $\varphi(z)$, which can be identified
 with a subspace of $(\Bbb C^{r+1})^{\vee}$, 
 and the fiber of 
 $f_{\Xi}^{\ast}\mathcal N_{\psi}^{\vee}$
 is the subspace $F_z\subset (\Bbb C^{r+1})^{\vee}$
 consisting of the vectors which also annihilate
 the kernel of the above map $\Xi_{\Bbb C}$.
Therefore, we need to show that if the numbers $\{a_{i, j}\}$ satisfy
 the conditions in the statement, the value at $z$ of the section of 
 $\mathcal N_{\varphi}^{\vee}\otimes\omega_{\Bbb P^1}$
 determined by $\{a_{i, j}\}$ 
 is contained in $F_z$.

But this is easy to see.
Namely, given the vector valued residue
 $\alpha_j = \sum_{i=1}^{r+1}a_{i, j}e_i^{\vee}$ at $y_j$
 as in the statement, the corresponding section of 
 $\mathcal N_{\varphi}^{\vee}\otimes\omega_{\Bbb P^1}$ is given by
\[
\sum_{j = 1}^{r+1}\alpha_j\sigma_j
\]
 as in Equation (\ref{eq:localobst}).
Since the vectors $\alpha_j$ all annihilate the kernel of 
 $\Xi_{\Bbb C}$, the value at any point of $\Bbb P^1$
 also annihilates it.
This proves the lemma. \qed

\subsection{Description of the dual obstruction spaces for 
 general global pre-log curves}\label{subsec:gen}
The argument so far establishes the description of the 
 dual space of obstructions restricted to a component of a pre-log curve
 corresponding to a higher valent vertex of a tropical curve.
Now we consider the description of the dual spaces
 of obstructions of the original pre-log curves.

Let $(\Gamma, h)$ be a 
 tropical curve satisfying Assumption A.
Let $\varphi_0\colon C_0\to \mathfrak X$ be a pre-log curve of 
 type $(\Gamma, h)$ in a suitable toric degeneration defined respecting 
 $(\Gamma, h)$.
We would like to give a description of elements of 
 the space 
 $H = H^0(C_0, \mathcal N_{\varphi_0}^{\vee}\otimes\omega_{C_0})$.

As in the case of Theorem \ref{thm:obstruction}, which described the
 space $H$ when the tropical curve $(\Gamma, h)$ is immersive, 
 one can prove 
 that the space $H$ is determined by the contributions
 from each bouquet of $h(\Gamma)$.

Let $L$ be the loop part of 
 $h(\Gamma)\subset N_{\Bbb R}\cong\Bbb R^n$
  (see Definition \ref{def:loops}
 and Remark \ref{rem:assump}).
Then, as before, 
 cut $L$ at each 3- or higher valent vertex.
Then we obtain a set of chains of segments
 $\{l_m\}$.
Let $U_m$ be the linear subspace of $N_{\Bbb C}$
 spanned by the direction vectors of the segments of $l_m$.
 
\begin{rem}
Under Assumption A, it can happen that the inverse image of a vertex $v$ of $h(\Gamma)$
 has several connected components.
Let $\{v_1, \dots, v_k\}$ be the vertices of $\Gamma_{\varphi_0}$ mapped to $v$, 
 see Definition \ref{def:intergraph}.
Let $n_1, \dots, n_k$ be the valence of these vertices.
Then when we talk about the valence of $v$, we regard it as a disjoint union of 
 the vertices $\{v_1, \dots, v_k\}$.
Namely, we think that at the position of $v$, $k$ vertices with valences $n_1, \dots, n_k$
 are overwrapping.
\end{rem} 
 
Then we can generalize Theorem \ref{thm:obstruction}
 as follows.
\begin{thm}\label{thm:obstruction2}
Elements of the space $H$ are described by the following data.
\begin{enumerate}[(I)]
\item Give the value zero to all the flags not contained in $L$.
\item Give a value 
 in $(U_m)^{\perp}\subset (N_{\Bbb C})^{\vee}$
 to each of the flags  associated to the edges of $l_m$.
\item The data in (I) and (II)
 give an element of $H$ if and only if the following 
 conditions are satisfied.
\begin{enumerate}
\item At each 3-valent vertex
 $p$ of $h(\Gamma)$, 
\[
u_1+u_2+u_3=0
\]
 holds as an element of $(N_{\Bbb C})^{\vee}$.
Here $u_1, u_2, u_3$ are the data attached to the three flags in 
 $h(\Gamma)$
 which have $p$ as the vertex.
\item At each $r+2$-valent vertex ($r\geq 2$) $q$, the data in (II)
 attached to the flags whose vertex is $q$ satisfy the condition
 in Lemma \ref{lem:genhigh}.
\item The data in (II) is compatible  
 in the sense that the sum of the values attached to the two flags
 of any edge of $l_m$ 
 is zero. \qed
\end{enumerate}
\end{enumerate} 
\end{thm}
There are several remarks for the notation in this theorem.
\begin{rem}
\begin{enumerate}
\item Note that the numbers $\{a_{i, j}\}$ in
 Lemma \ref{lem:genhigh} are associated to edges of 
 the tropical curve $(\Gamma_0, h_0)$ in $\Bbb R^{r+1}$,
 and the data in $(II)$ determine these numbers through the
 projection argument in Subsection \ref{subsec:gendual}.
\item When the subspace spanned by the directions of the edges emanating
 from the higher valent vertex $q$ is strictly smaller than
 $\Bbb R^n$, there are also the $c_{i, j}$ parts
 of (\ref{eq:localobst}) in Subsection \ref{subsec:localobst}
 (these $c_{i, j}$ do not have local constraints at each vertex,
 contrary to $\{a_{i, j}\}$).
\item Also note that the numbers $\{a_{i, j}\}$ and $\{c_{i, j}\}$
 are associated to $r+1$ edges among the $r+2$ edges emanating
 from $q$.
The residue at the point corresponding to the remaining 
 edge is determined from them by the residue theorem.
For some choice of $r+1$ edges out of the total $r+2$ edges, the numbers $\{a_{i, j}\}$
 constructed from the data in $(I)$, $(II)$ and this choice of the edges satisfy 
 the condition $(III)(b)$ above if and only if the numbers $\{a_{i, j}\}$ satisfy it
 for any choice of $r+1$ edges.
\end{enumerate}
\end{rem}

The proof of the theorem 
 is the same as that of Theorem \ref{thm:obstruction}.
This is the most general form of the description of the space $H$
 for pre-log curves corresponding to tropical curves satisfying 
 Assumption A.

 In view of Theorems \ref{thm:obstruction}
  and \ref{thm:obstruction2}, we give the following definition.
\begin{defn}\label{def:abundancysupport}
Let $(\Gamma, h)$ be a tropical curve satisfying Assumption A.	
Then the \emph{support of superabundancy} of $(\Gamma, h)$ is
 the closed subgraph $\Gamma_{ss}$ of $h(\Gamma)$ such that 
 for any edge $\mathfrak E$ of $\Gamma_{ss}$, there is an element of 
 $H$ such that the value of it on the flags associated to $\mathfrak E$
 is not zero.
\end{defn}
Large part of results which are valid for superabundant curves of genus one
 can be extended to those tropical curves whose support of superabundancy
 is a loop.

\section{Regular tropical curves and smoothability}\label{sec:4}
While our main concern is superabundant tropical curves, 
 in this section we deal with those tropical curves which are not superabundant.
We include this section since the technique used here is also useful in the later 
 study of superabundant tropical curves.
On the other hand, the same technique allows us to deduce the optimal correspondence theorem
 for regular tropical curves (see Theorem \ref{thm:regsm}), generalizing 
 results in \cite{CFPU, N3, Tyo}.

\subsection{Comparison of two definitions of superabundancy}\label{subsec:compare}
We defined the superabundancy of a tropical curve
 in Definition \ref{def:superabundancy}.
On the other hand, there is another known definition of superabundancy.

Let $\Gamma$ be a weighted abstract graph as in the beginning of 
 Subsection \ref{subsec:pre}, but now we allow
 that $\Gamma$ is not necessarily 3-valent.
Let $h\colon \Gamma\to N_{\Bbb R}$ be an embedding giving $\Gamma$
 a structure of a tropical curve.
We identify the graph $\Gamma$ with its image.
Let $L$ be the loop part of $\Gamma$ 
 and let $L^{[1]} = \{E_1, \dots, E_l\}$ be the set of edges in $L$
 seen as 1-chains of the simplicial homology group
 (in particular, we choose an orientation for each $E_i$).
Let $u_i$ be the primitive integral vector in the direction of $E_i$.
\begin{defn}[{\cite[Definition 4.1]{CFPU}, see also \cite[Section 2.6]{M}
 and \cite[Section 1]{Katz}}]\label{def:nonsuperabundant2}
The tropical curve $(\Gamma, h)$ is \emph{regular}
 if the following abundancy map
\[
\Phi_{(\Gamma, h)}\colon
 \Bbb R^{L^{[1]}}\to Hom(H_1(\Gamma), N_{\Bbb R}),
\]
\[
(\ell_{E_i})\mapsto \left(\sum_i a_i[E_i]\mapsto \sum_i\ell_{E_i}a_iu_i
 \right)
\]
 is surjective.
\end{defn}
\begin{rem}
In \cite{CFPU}, a tropical curve satisfying this condition is called non-superabundant.
We call it differently, since some people already call these tropical curves regular, 
 and we have used the word non-superabundant in Definition \ref{def:superabundancy}.
Also, from the point of view of Theorem \ref{thm:regsm} below, the use of the
 word regular for these tropical curve
 seems reasonable.
\end{rem}

First we remark that this definition is the same as the following.
Let us write 
\[
rank H_1(\Gamma) = r
\]
 and choose $r$ edges 
 $F_1, \dots, F_r$ from $L^{[1]}$ so that 
 the graph $\Gamma\setminus\{F_1^{\circ}, \dots, F_r^{\circ}\}$
 is a tree.
Let 
\[
N_i = N_{\Bbb R}/\Bbb R\cdot v_i
\]
 be the quotient space,
 where $v_i$ is the primitive integral vector in the direction of 
 $F_i$.
Let $\partial_1F_i$ and $\partial_2 F_i$ be the ends of $F_i$.

Note that since $\Gamma\setminus\{F_1^{\circ}, \dots, F_r^{\circ}\}$
 is a tree, for every 
 $(\ell_{E_i})\in \Bbb R^{L^{[1]}\setminus\{F_1, \dots, F_r\}}$,
 there is a corresponding map from 
 $\Gamma\setminus\{F_1^{\circ}, \dots, F_r^{\circ}\}$ to $N_{\Bbb R}$
 with the length of $E_i$ given by $\ell_{E_i}$
 (the minus length means the direction inverse to the given direction in 
 the original map $h$). 
There are many such maps due to parallel transports, but the vectors
 $\partial_1F_i-\partial_2F_i$, $i = 1, \dots, r$ do not depend on the choice.

\begin{lem}\label{lem:nonsuperabundant}
The tropical curve $(\Gamma, h)$ is regular if and only if
 the map
\[
\Bbb R^{L^{[1]}\setminus\{F_1, \dots, F_r\}}
 \to \oplus_{i = 1}^rN_i,
\]
\[
(\ell_{E_i})\mapsto (\overline{\partial_1F_i-\partial_2F_i})
\]
 is surjective.
Here $(\overline{\partial_1F_i-\partial_2F_i})$ is the image of  
 $\partial_1F_i-\partial_2F_i\in N_{\Bbb R}$ in $N_i$.  
 \qed
\end{lem}
From this observation, we can see the following.
\begin{prop}\label{prop:immexist}
Let $(\Gamma, h)$ be a regular tropical curve.
Then it can be deformed into an immersive (in fact, embedded) 3-valent tropical curve
 of the same genus
 in the sense of Definitions \ref{def:tropical curve} and \ref{immersive}.
Moreover, the deformed immersive tropical curve is also 
 regular.
\end{prop}
\proof
We identify $\Gamma$ and $h(\Gamma)$ as noted above.
Let $V\in \Gamma$ be a higher valent vertex.
It suffices to show that the star of any such $V$ can be deformed into 
 a 3-valent graph while keeping the 
 balancing condition and the homotopy type of
 the entire graph.

As in the above argument, take edges 
 $F_1, \dots, F_r$ of $\Gamma$ so that
 the graph $\Gamma\setminus\{F_1^{\circ}, \dots, F_r^{\circ}\}$
 is a tree.
The interior $F_i^{\circ}$ of 
 these $F_i$ may be either contained or not contained in the star of
 $V$.

Then perturb the star of $V$ so that it becomes
 a trivalent tree graph $\Gamma_V$
 satisfying the balancing condition.
Since $\Gamma\setminus\{F_1^{\circ}, \dots, F_r^{\circ}\}$
 is a tree, we can replace the intersection of the star of $V$ and
 $\Gamma\setminus\{F_1^{\circ}, \dots, F_r^{\circ}\}$ by 
 the graph $\Gamma_V\setminus\{F_1^{\circ}, \dots, F_r^{\circ}\}$
 without changing the homotopy type
 and the slopes of the other edges.
 On the other hand, the difference 
 $\overline{\partial_1F_i-\partial_2F_i}$ (which is 
 clearly 0 in the original 
 $\Gamma$)
 may have 
 slightly changed through this process.

However, due to Lemma \ref{lem:nonsuperabundant}, 
 we can perturb the resulting graph so that 
 the difference $\overline{\partial_1F_i-\partial_2F_i}$
 becomes 0.
Note that since the domain of the map of 
 Lemma \ref{lem:nonsuperabundant} consists of the lengths of the 
 part of the edges of $\Gamma$ (and does not contain
 those of the newly introduced edges in $\Gamma_V$), 
 this perturbation can be done without changing $\Gamma_V$.
 
Then the resulting graph can be completed by adding the edges 
 $F_i$ to produce a tropical curve which has the same homotopy
 type as $\Gamma$.
Note that the resulting tropical curve is still 
 regular
 (since introducing the new edges makes the domain 
 of the abundancy map larger, the surjectivity of the map is not affected).
Repeating this process at all the higher valent vertices, we obtain
 an embedded 3-valent tropical curve which is a deformation of 
 $(\Gamma, h)$.\qed\\

In particular, we see the following.
\begin{cor}\label{cor:nonsuperabund}
The set of regular tropical curves 
 is contained in the set of non-superabundant tropical curves.
\end{cor}
\proof
It suffices to prove that 
 a deformation of a given regular tropical curve
 into an immersive tropical curve
 (which exists by Proposition \ref{prop:immexist}) 
 is non-superabundant in the sense of 
 Definition \ref{def:superabundancy}.
This is equivalent to show that 
 the dual obstruction space $H$ vanishes for 
 pre-log curves of the type given by the deformed immersive tropical curve.
The existence of such a curve can be proved by the argument in 
 \cite[Proposition 5.7]{NS} and the surjectivity of the abundancy map, 
 see the argument before Proposition \ref{prop:abundancymap} below.
In fact, as we noted in Remark \ref{rem:thm40add}, 
 the space $H$ makes sense even without the existence of a
 pre-log curve of the corresponding type, and the following argument
 is also valid in this sense.
 
Now \cite[Proposition 4.2]{CFPU} 
 implies that, given a tropical curve, if the abundancy map is surjective, 
 then
 the space $H$ vanishes for 
 pre-log curves of the type given by the tropical curve.
Thus, tropical curves which are regular is non-superabundant.\qed

\begin{rem}
%\item In the proof of the corollary, we need to deform the 
% original tropical curve to an immersed one since to calculate the space
% $H$, we need to construct a pre-log curve of the given type, while
% we have not proved a general existence of pre-log curves 
% corresponding to a tropical curve which is not immersive.
%When the tropical curve is immersive and 
% non-superabundant in the sense of 
% Definition \ref{def:nonsuperabundant2}, this existence is proved below.
The statement of the corollary is a little imprecise in that
 regular tropical curves in the sense of 
 Definition \ref{def:nonsuperabundant2} do not have the data of 
 3-valent domain graphs which we always assume to be
 associated to tropical curves.
Therefore, the inclusion makes sense when we specify this data to each 
 regular tropical curve
 (see also Corollary \ref{cor:pre-loggen} below).
\end{rem}

The proof of the proposition actually shows the following stronger statement.
\begin{cor}\label{cor:pre-loggen}
Let $(\Gamma, h)$ be a regular tropical curve.
Let $\{V_i\}$ be the set of higher valent vertices.
The star of each $V_i$ can be seen as an embedded tropical curve
 with one vertex.
Let $\Gamma_{V_i}$ be any deformation of the star of $V_i$
 into an embedded 3-valent tropical curve.
Let $\tilde{\Gamma}$ be the abstract graph obtained from $\Gamma$
 by replacing the star of each $V_i$ by $\Gamma_{V_i}$.
The embeddings of $\Gamma$ and $\Gamma_{V_i}$ into $N_{\Bbb R}$
 determine the combinatorial type of the graph $\tilde{\Gamma}$.
Then $(\Gamma, h)$ can be deformed into an embedded 3-valent
 tropical curve whose domain is $\tilde{\Gamma}$ with this combinatorial
 type.\qed
\end{cor}

\begin{rem}
	Note that the graph $\Gamma_{V_i}$ need not be a tree.
	In this case, of course the resulting tropical curve may not be non-superabundant
	in general.
\end{rem}

When the tropical curve is not regular, it may not deform into an immersive tropical curve, 
 see Example \ref{ex:singular}.
However, smoothable tropical curves known so far behave well.
Therefore, we can pose the following question.
\begin{q}\label{q:1}
Is there a smoothable tropical curve which does not deform into an immersive tropical curve?
\end{q}

%Here we only fix the genus and we can vary the domain curve $\Gamma$ and the combinatorial type.
In \cite{N3}, it is proved that for curves of genus one, there is no such a tropical curve.
Note that it can happen that a given tropical curve is smoothable, 
and it deforms into immersive tropical curves, but any deformed immersive curve is not smoothable
(see Example \ref{ex:2})

\subsection{General correspondence theorem for regular tropical curves}\label{subsec:corrthm}

As we discussed in \cite[Proposition 7.3]{NS}, 
 the deformation of tropical curves
 of a fixed type is closely related to the construction of 
 a pre-log curve of the corresponding type.
It is easy to see that if $(\Gamma, h)$ is a tropical curve
 satisfying Assumption A
 where $\Gamma$ is a tree, then a pre-log curve of the corresponding
 type always exists, regardless of whether there are higher valent
 vertices or not.

\begin{figure}[h]
\includegraphics{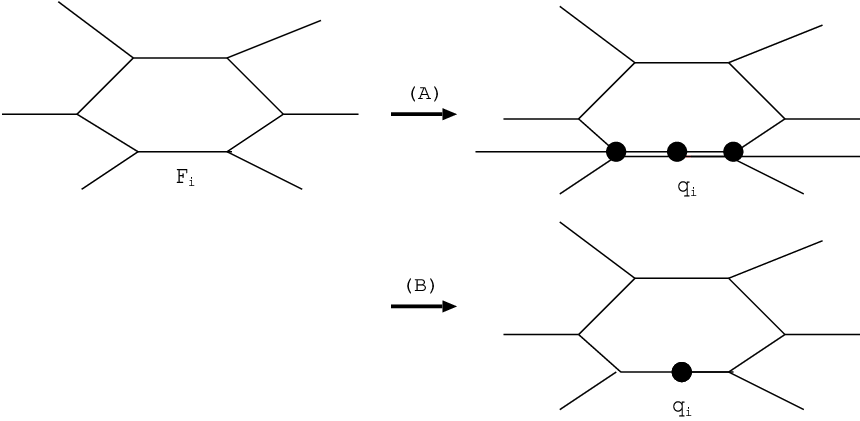}
\caption{First cut the edge $F_i$ at the midpoint, extend the cut edges
 to infinity, and add suitable 2-valent vertices (process (A)).
We obtain a tropical curve $(\Gamma', h')$.
Then consider pre-log curves of type $(\Gamma', h')$. 
If the two components corresponding to the vertex $q_i$
 coincide, then we can discard the other components corresponding to the 2-valent
 vertices on the
 new unbounded edges, and obtain another pre-log curve
 which is of type $(\Gamma, h)$, with a 2-valent vertex $q_i$ added
 on the edge $F_i$ (process (B)).
}\label{fig:cut_and_glue}
\end{figure}

Suppose we have an embedded regular tropical curve $(\Gamma, h)$,
 and choose the edges $\{F_1, \dots, F_r\}$ as before
 ($\Gamma$ is not necessarily 3-valent).
Then cut each $F_i$ at the midpoint and extend the cut edges to infinity.
Also, add 2-valent vertices to the extended edges at the points
 where there were the midpoints and where the extended edges intersect 
 other vertices
 (see Figure \ref{fig:cut_and_glue}).
The result can naturally be regarded as an immersed tropical curve 
 of genus zero.
The midpoint $q_i$
 of $F_i$ 
 now becomes the common image of 2-valent vertices on 
 two newly produced unbounded edges. 
Let us write this tropical curve by $(\Gamma', h')$.
By the remark in the preceding paragraph, there is a pre-log curve 
 $\varphi_0'\colon C_0'\to X_0$
 of type 
 $(\Gamma', h')$ in the central fiber of a suitable toric degeneration.

Let us write by $C'_{0, i_1}$ and $C'_{0, i_2}$ the components of 
 $C'_0$ corresponding to the vertex $q_i$
 (as we noted above, since $q_i$ is the
 image of vertices on two different edges, 
 there are two components of $C'_0$ corresponding to it).
The images of these components lie in the same component of $X_0$ and
 are the closures of orbits of the 
 one dimensional subtorus corresponding to the direction of $F_i$. 
However, these orbits are different in general.
When these images coincide, 
 then by discarding other components of $C_0'$ corresponding to 
 the 2-valent vertices of the 
 new unbounded edges, 
 one obtains a pre-log curve of positive genus 
 (see Figure \ref{fig:cut_and_glue}).

Let us write by 
\[
\Bbb G_i = N\otimes\Bbb C^*/\Bbb G_m\cdot u_i
\]
 the torus which is the quotient of the big torus acting on the 
 components of $X_0$ by the one dimensional sub-torus corresponding
 to the direction of $F_i$.
The difference of the positions of the above two orbits determines
 a point of $\Bbb G_i$ (well-defined up to inversion. We fix
 one choice).

The pre-log curve $\varphi_0'(C_0')$ 
 can be regarded as the result of gluing of the 
 components corresponding to the vertices of $\Gamma'$.
By the argument in \cite[Proposition 7.3]{NS}, the freedom of 
 changing the lengths of the edges of
 $L^{[1]}\setminus \{F_1, \dots, F_r\}$
% , after tensoring by $\Bbb C^*$
% (here $\Bbb R$ acts on $\Bbb C^*$ by $r\cdot \zeta =  e^r\zeta$),
 corresponds to transporting each component of $C_0'$ by 
 the torus action.
Now we briefly recall the construction.
% on the corresponding 
% component of $X_0$, keeping the homotopy type of the whole curve $\varphi_0'(C_0')$.

Namely, let $(\Gamma_1, h_1)$ be a tropical curve where $\Gamma_1$
 is a tree and the image $h_1(\Gamma_1)$ has only one vertex (which we assume
 to be the origin).
Let $\psi\colon C\to \Bbb P$ be a torically transverse curve of type
 $(\Gamma_1, h_1)$.
Here $\Bbb P$ is the toric variety defined by $h_1(\Gamma_1)$ seen as a fan.
Then $C$ is a nonsingular rational curve with $|\Delta_{\Gamma_1}|$
 marked points.
Here $|\Delta_{\Gamma_1}|$ is the number of unbounded edges of $\Gamma_1$.
Each unbounded edge of $h_1(\Gamma_1)$ determines a one dimensional 
 subtorus of the torus acting on $\Bbb P$.

Now the tropical curve $(\Gamma', h')$ is obtained by gluing pieces
 like $h_1(\Gamma_1)$.
Let $V$ be the number of vertices of $h'(\Gamma')$
 and let $\{v_1, \dots, v_V\}$ be the set of these vertices.

Each vertex has the freedom of parallel transport, which we write by $\Bbb R^n$.
In the product space $(\Bbb R^n)^V$, let $T$ be the set of the directions 
 tangent to the parameter space of tropical curves containing $(\Gamma', h')$.
This means that if $a = (a_1, \dots, a_V)\in T$, 
 and if $e$ is a bounded edge of $h'(\Gamma')$ whose ends are $v_i, v_j$, 
 then $a_i\equiv a_j$ in $\Bbb R^n/\Bbb R\cdot u_e$.
Here $u_e$ is the direction of the edge $e$.
Discarding the freedom of the overall parallel transport, we consider the 
 quotient space $T/\Bbb R^n$.
We may think this as fixing the position of one vertex (say, $v_1$)
 so that $a_1 = 0$.
Since $\Gamma'$ is a tree, each edge has a well-defined direction 
 pointing opposite to the vertex $v_1$.
 
The space $T/\Bbb R^n$ 
 has a natural integral structure, and its set of generators  is given
 by the integral generators of the bounded edges of $h'(\Gamma')$
 in the direction
 mentioned above.
Namely, since $\Gamma'$ is a tree, deformations of $h'(\Gamma')$ 
 (modulo parallel transports) are given by changing the lengths of the bounded edges.
Let $b_1, \dots, b_c$ be these set of generators.
Then, 
\[
T/\Bbb R^n\cong \oplus_{i=1}^c\Bbb Z b_i\otimes\Bbb R.
\]

Changing the length of the edge $e_i$ of the direction $b_i$ by unit integral length
 deforms the map $h'$ and displaces some vertices.
Namely, the edge $e_i$ divides $h'(\Gamma')$ into two parts, 
 one of them contains the vertex $v_1$. 
The above change transports the vertices on the part not containing $v_1$
 by $b_i$.
This attaches a vector in $(\Bbb R^n)^{V-1}$ to $b_i$, 
 which we write by $\tau(b_i)$ (the $-1$ corresponds to the cancelled freedom
 of the vertex $v_1$).
Therefore, we have a homomorphism 
\[
\tau\colon T/\Bbb R^n\to (\Bbb R^n)^{V-1}
\]
 extending the map
\[
b_i\mapsto \tau(b_i)
\]
 linearly.

Furthermore, recall that the graph $h'(\Gamma')$ has a pair of
 2-valent vertices associated to the edges $F_1, \dots, F_r$.
We write them by $v_{1, 1}, v_{1, 2}, v_{2, 1}, v_{2, 2}, \dots, v_{r, 1}, v_{r, 2}$.
These are a part of the set of vertices $\{v_1, \dots, v_V\}$ above.
Therefore, to each $b_i$ and $v_{j, k}$, we attach a vector
 $\tau_{j, k}(b_i)\in\Bbb R^n$.
Restricting the domain of the map $\tau$ to the subspace corresponding to the
 edges in $L^{[1]}\setminus \{F_1, \dots, F_r\}$, 
 we have another linear map
\[
\bar\tau\colon \Bbb R^{L^{[1]}\setminus \{F_1, \dots, F_r\}}\to
 \oplus_{i=1}^r(\Bbb R^n/\Bbb R\cdot u_i)
\]
 given by extending the map
\[
b_i\mapsto (\tau_{j, 1}(b_i)-\tau_{j, 2}(b_i))_{j=1, \dots, r}
\]
 linearly.
Here $u_i$ is the direction of the edge $F_i$.
The regularity of $(\Gamma, h)$ is equivalent to the claim that
 this map is surjective.

Now let $X_0'$ be the central fiber of a toric degeneration defined respecting $(\Gamma', h')$.
For each vertex $v_i$, there is the corresponding component $X_{0, i}'$.
Let $\varphi_0'\colon C_0'\to X_0'$ be a pre-log curve of type $(\Gamma', h')$, 
 which exists since $(\Gamma', h')$ is a tree.
For any component $X_{0, i}'$, it gives a
 pre-log curve $\varphi_{0, i}'\colon C_{0, i}'\to X_{0, i}'$
 of type $v_i$ in $X_{0, i}'$ in the sense of Definition \ref{def:typev}.
Here we regard the star of $v_i$ in $h'(\Gamma')$ as a (open subset of) tropical curve
 with one vertex.

Consider the group $(\Bbb C^*)^{L^{[1]}\setminus \{F_1, \dots, F_r\}}$.
Let $\zeta = (\zeta_1, \dots, \zeta_K)$ be an element of it, here 
 $K = \sharp(L^{[1]}\setminus \{F_1, \dots, F_r\})$.
It acts on each component $\varphi_{0, i}'(C_{0, i}')$ of $\varphi_0'(C_0')$
 by the action of $(\zeta_1^{\tau_i(b_1)_1}\cdots\zeta_K^{\tau_i(b_K)_1}, 
  \dots, \zeta_1^{\tau_i(b_1)_n}\cdots\zeta_K^{\tau_i(b_K)_n})\in N\otimes\Bbb C^*$.
Here $\tau_i(b_j) = (\tau_i(b_j)_1, \dots, \tau_i(b_j)_n)\in\Bbb Z^n$ is the component of
 $\tau(b_i)\in (\Bbb Z^n)^{V-1}$ corresponding to the vertex $v_i$.

By construction, the action of $(\Bbb C^*)^{L^{[1]}\setminus \{F_1, \dots, F_r\}}$
 on each component is compatible with $\varphi_0'$ in the sense that
 the resulting components on $X_{0, i}'$ glue again so that 
 it can be regarded as the image of another pre-log curve 
 $\zeta\cdot\varphi_{0}'\colon C_0'\to X_0'$
 of type $(\Gamma', h')$.
Note that the domain curve $C_0'$ does not change through the action.

In this way, we have an action of $(\Bbb C^*)^{L^{[1]}\setminus \{F_1, \dots, F_r\}}$
 on the set of pre-log curves of type $(\Gamma', h')$ with the fixed domain curve $C_0'$.
Recall that the curve $C_0'$ contains two components corresponding to 
 each of the edges $\{F_1, \dots, F_r\}$.
Let $D_{i, 1}, D_{i, 2}$ be these components corresponding to the edge $F_i$.
They are mapped by $\varphi_0'$ to the same component $X_{0, i}'$
 of $X_0'$, and their images are the closure of the orbit of the
 action of a one dimensional 
 subtorus of $\Bbb T$. 
The difference of these orbits gives an element of $\Bbb G_i = N\otimes\Bbb C^*/\Bbb G_m\cdot u_i$
 introduced above.

This determines the map 
\[
\Phi_{(\Gamma, h)}^{\Bbb C}\colon
 (\Bbb C^*)^{L^{[1]}\setminus \{F_1, \dots, F_r\}}\to \prod_{i=1}^r\Bbb G_i
\]
 which is the $\Bbb C^*$-tensored version of the map $\bar\tau$ above.

Therefore, the regularity of $(\Gamma, h)$ 
 implies the following statement for the pre-log curves.
\begin{prop}\label{prop:abundancymap}
The above map $\Phi_{(\Gamma, h)}^{\Bbb C}$ is surjective.\qed
\end{prop}
This results in the following proposition.
Let $(\Gamma, h)$
 be an embedded regular tropical curve (recall that in this case
 $\Gamma$ need not be 3-valent).
For each vertex $v$ of $h(\Gamma)$, take any torically transverse curve
 $C_{0, v}$ of type $v$ (see Definition \ref{def:typev}) 
　in the corresponding
 component $X_{0, v}$ of the central fiber $X_0$
 of a toric degeneration defined respecting $(\Gamma, h)$.
Note that $C_{0, v}$ is a rational curve with $r_v$ fixed marked points,
 where $r_v$ is the number of the edges emanating from $v$.
\begin{prop}\label{prop:exist}
Under the situation of the above paragraph, 
 there is a pre-log curve of type $(\Gamma, h)$
 whose component corresponding to $v$ is isomorphic to 
 the given $C_{0, v}$ as a rational curve with marked points.\qed
\end{prop}
In particular, regular tropical curves can be realized 
 at least as a degenerate holomorphic curve.

Now the main result of \cite{CFPU} states that once we have
 a pre-log curve whose type
 is given by a regular tropical curve, then it is smoothable.
Combining it with Proposition \ref{prop:exist}, we have the following 
 optimal form of the corresponding theorem for regular tropical curves.

\begin{thm}\label{thm:regsm}
Any regular tropical curve is smoothable.\qed
\end{thm}
\begin{rem}
The proof of this theorem 
 was originally given in the first version of the preprint \cite{N3} and also by Tyomkin \cite{Tyo}
 in immersive cases.
Recently it was proved 
 by Cheung, Fantini, Park and Ulirsch \cite{CFPU}
with broader applicability, but still under technical assumptions. 
\end{rem}

\begin{rem}\label{rem:degdeform}
Even when $(\Gamma, h)$ is not regular, once one knows the existence of 
 a pre-log curve of type $(\Gamma, h)$, then the above argument shows that
 the set of deformations of the pre-log curve with the fixed domain curve
 can be essentially identified with the set of deformations of $(\Gamma, h)$
 with $\Gamma$ fixed (in this case, 
 we assume $h$ is an immersion while we do not assume $\Gamma$ to be 3-valent). 
Precisely, the tangent space of the parameter space of the pre-log curves
 is canonically isomorphic to that of the parameter space of 
 the tropical curves with $\Bbb C^*$ tensored.
Generally speaking, the space of deformations of pre-log curves are the same as that of 
 tropical curves, and there are no obstruction to these deformations, 
 contrary to the deformation to the curves on general fibers. 
\end{rem}

\section{Higher valent vertices and smoothability}\label{sec:5}
In this section and the next, we study aspects of degenerate algebraic curves
 whose dual intersection complexes have higher valent vertices.
Usually, the study of algebraic curves through tropical curves depends on 
 the rather precise correspondence between the freedom of deforming these curves.
In particular, 3-valent immersive tropical curves constitute 
 a natural class to study (\cite{M, S, Tyo}),
 since a 3-valent vertex corresponds to a rational curve with three special points, 
 which has no moduli.

However, in the study of curves with positive genus, tropical curves with
 higher valent vertices appear all over the place. 
These vertices correspond to rational curves with more than three 
 special points, which have their own moduli.
Therefore, at first sight the tropical curves seem to miss lots of information
 about algebraic curves.
For example, even when one can resolve a higher valent vertex into 3-valent vertices
 by deforming a given tropical curve, such a freedom does not always correspond
 to a freedom to deform algebraic curves.  
Nevertheless, we will see that combinatorial information of tropical curves 
 quite effectively restore the information of algebraic curves.
This is particularly true when the curve is genus one, and in this case 
 we can actually read off all the information about the smoothability of 
 corresponding degenerate algebraic curves. 
 
The process in the study of correspondence between tropical curves
 and algebraic curves is two fold. 
The first is constructing degenerate algebraic curves of the type 
 specified by a given tropical curve.
The second is the calculation of the obstruction to deform the degenerate curves.
So far the second point has been studied by many authors.
This should be partly because the first point is not very intricate for regular tropical curves, 
 as we saw in the previous section.
However, the first point turns out to be as crucial as the second one 
 when one considers superabundant tropical curves.

We consider the second point for superabundant tropical curves with higher valent 
 vertices in this section.
In the next section, we solve the first point for tropical curves of genus one.
\subsection{Examples of tropical curves with higher valent vertices and obstruction}\label{subsec:highval}
In this subsection, we begin with several examples of tropical curves
 relevant to the theory we developed so far.
Particular emphasis is put on examples 
 which exhibit phenomena peculiar to higher valent vertices.

\begin{example}\label{ex:1}
First we consider the following tropical curve
 $(\Gamma, h)$ of genus one in $\Bbb R^3$,
 see Figure \ref{fig:0X}.

\begin{figure}[h]
\includegraphics[height=4cm]{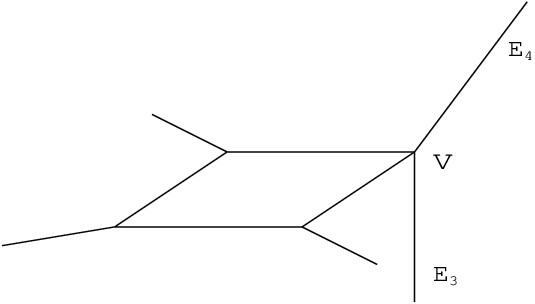}
\caption{}\label{fig:0X}
\end{figure}

Here the vertex $V$ is the standard 4-valent vertex (that is, 
 the directions of the edges emanating from it are
 $-e_1, -e_2, -e_3$ and $f = e_1+e_2+e_3$).
The direction of the edge $E_4$ is $f$ and that of $E_3$ is $-e_3$.
The complement of the edges $E_3$ and $E_4$ in $h(\Gamma)$,
 which we write by $L^{\circ}$, is contained in 
 an affine plane parallel to the plane spanned by $e_1$ and $e_2$.
The abstract 3-valent graph $\Gamma$ is isomorphic to the graph 
 which resolves the vertex $V$ into two 3-valent vertices by introducing
 an edge of direction $e_1+e_2$.

It is easy to see that this is a superabundant tropical curve.
Moreover, there is no pre-log curve of type 
 $(\Gamma, h)$ in any toric degeneration respecting $(\Gamma, h)$.
Namely, since the part $L^{\circ}$ is contained in an affine plane, 
 the function $z$ corresponding to $e_3^{\vee}$ (where
 $\{e_1^{\vee}, e_2^{\vee}, e_3^{\vee}\}$ is the dual basis of 
 $\{e_1, e_2, e_3\}$) restricted to
 the part of a pre-log curve corresponding to $L^{\circ}$ has the constant value
 (note that the vertex $V$ is not contained in $L^{\circ}$). 
Therefore the curve corresponding to the vertex $V$ should have the same value
 of $z$ at the two points corresponding to the 
 two bounded edges emanating 
 from it.
However, since the curve corresponding to $V$ is a line in $\Bbb P^3$, 
 this is impossible if we also require the torically transversality. \qed
\end{example}

\begin{example}\label{ex:1.5}
There are even embedded tropical curves to which a pre-log curve
 of the given type does not exist.
The following is one of such tropical curves (G. Mikhalkin also 
 found a similar example).

Consider the following tropical curve of genus two in 
 $\Bbb R^3 = \Bbb Z^3\otimes \Bbb R$, see Figure \ref{fig:16}.

\begin{figure}[h]
\includegraphics[height=5cm]{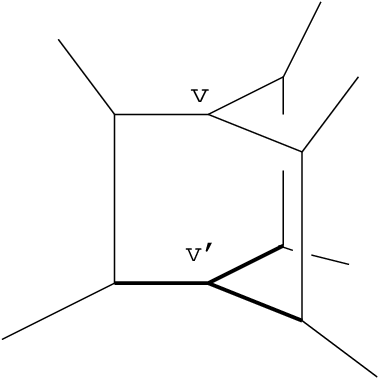}
\caption{}\label{fig:16}
\end{figure}

Here the directions and the lengths of the edges emanating from 
 the vertex $v'$ (drawn by bold lines in the figure)
  are exactly the same as those emanating from the vertex
 $v$, but all the edges emanating from $v'$ have weight two, 
 while all the other edges have weight one.
Assume that the direction vectors of the edges emanating from 
 $v$ (and also $v'$) generates the lattice given by the intersection of 
 $\Bbb Z^3$ and the affine plane spanned by these edges.
 
The loop part of this tropical curve is the following (see Figure \ref{fig:17}).

\begin{figure}[h]
\includegraphics[height=4cm]{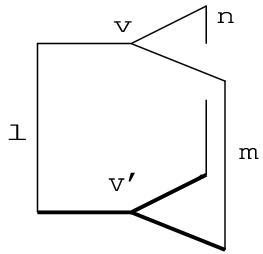}
\caption{}\label{fig:17}
\end{figure}

By cutting at the 3-valent vertices, 
 the loop part decomposes into three parts $l, m$ and $n$.
Let $L_l, L_m$ and $L_n$ be the two dimensional subspaces of 
 $\Bbb R^3$ spanned by the direction vectors of the edges
 contained in $l, m$ and $n$, respectively.
Let $L_l^{\perp}, L_m^{\perp}$ and $L_n^{\perp}$ be the 
 one dimensional subspaces of the dual space $(\Bbb R^3)^{\vee}$
 which annihilate vectors in $L_l, L_m$ and $L_n$ respectively.

By the above assumption, the integral generators of any two of 
 the spaces $L_l^{\perp}, L_m^{\perp}$, $L_n^{\perp}$
 generate a rank two saturated sublattice of $(\Bbb Z^3)^{\vee}$.

In particular, these generators can be thought of as affine coordinates of
 the projective plane $\Bbb P^2$,
 and the lines (in the sense of \cite[Definition 5.1]{NS}, see 
 Subsection \ref{subsec:pre-logdef})
 corresponding to the vertices $v$
 and $v'$ can be identified with lines in $\Bbb P^2$.
The line corresponding to the vertex $v$ is the line (that is, a curve of 
 degree one) in the usual sense.
However, the line corresponding to the 
 vertex $v'$ has degree two and intersects each toric divisor with
 multiplicity two.

Now let $x, y$ be affine coordinates of $\Bbb P^2$ mentioned above.
Then a torically transverse degree one curve can be described by 
 the defining equation
\[
px+qy+1 = 0,
\]
 where $p$ and $q$ are nonzero complex numbers.
Its intersection with the toric divisors have coordinates
\[
(x, y) = \left(-\frac{1}{p}, 0\right), \,\, \left(0, -\frac{1}{q}\right)
\]
and
\[
\left(\frac{1}{x}, \frac{y}{x}\right) = \left(0, -\frac{p}{q}\right).
\]
In particular, writing $\alpha = -\frac{1}{p}$, $\beta = -\frac{1}{q}$
 and $\gamma = -\frac{p}{q}$, we have
\[
\alpha\beta^{-1}\gamma = -1.
\]
On the other hand, one can see that the same product of the coordinates
 of the intersection between a line in $\Bbb P^2$
 corresponding to the vertex $v'$ and
 the toric divisors is one.
Thus, there is no pre-log curve of the type given by the above tropical curve.
\end{example}

\begin{example}\label{ex:2}
Let us consider a tropical curve
 which is slightly different from the one in Example \ref{ex:1},
 see Figure \ref{fig:0Y}.
This is also a superabundant tropical curve as in Example \ref{ex:1}.

\begin{figure}[h]
\includegraphics[height=4cm]{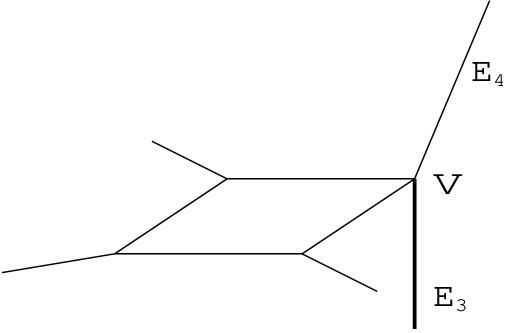}
\caption{}\label{fig:0Y}
\end{figure}

Now the vertex $V$ is the standard 4-valent vertex magnified
 in the $e_3$ direction by two times, and the edge $E_3$ has weight two. 
The part $L$ is not affected by this, and the condition that 
 the component of a pre-log curve corresponding to the vertex $V$
 should have the same value of the function $z$ at
 the two points associated to the edges $E_1$ and $E_2$ 
 can now be satisfied.
Here $z$ is the function corresponding to 
 $e_3^{\vee}$.
Also, $E_1$ and $E_2$ are the bounded edges emanating from $V$.
In particular, their directions are generated by $-e_1$ and $-e_2$, 
 respectively.

Explicitly, take any line (in the usual sense) in $\Bbb P^3$
 whose $z$-coordinates at the intersections with the toric divisors
 corresponding to $E_1$ and $E_2$ 
 are $\alpha$ and $-\alpha$, respectively. 
Here $\alpha$ is any nonzero complex number.
Then the image of such a line by the map between toric varieties 
 associated to the map magnifying the $e_3$ direction by two times
 gives the desired pre-log curve corresponding to the vertex $V$.
Thus, in this example there is a pre-log curve of type $(\Gamma, h)$.

Now let $\varphi_0\colon C_0\to X_0$ be such a pre-log curve and 
 let us calculate the dual obstruction based on Theorem \ref{thm:obstruction2}.
This calculation is done by solving the equations imposed on 
 the numbers $\{a_{i, j}\}$ associated to the vertex $V$, 
 see Subsection \ref{subsec:localobst}.
According to the Lemma \ref{lem:genhigh}, we attach to the 
 edges $E_1, E_2$ and $E_3$ the numbers
 $\{a_{2, 1}, a_{3, 1}\}$, $\{a_{1, 2}, a_{3, 2}\}$
 and $\{a_{1, 3}, a_{2, 3}\}$, respectively.
 
Since $E_3$ is an unbounded edge, the numbers  
 $a_{1, 3}, a_{2, 3}$ must be zero.
Let $\zeta$ be an affine
 coordinate of the component $C_{0, V}$ of $C_0$
 whose value is zero
 at the node corresponding to the edge $E_1$.
Let $p_2$ and $p_3$ be the values of $\zeta$  
 at the node corresponding to the edge $E_2$ and
 the point corresponding to 
 the unbounded edge $E_3$, respectively. 
Moreover, by condition (III)(c) of Theorem \ref{thm:obstruction2},
 $a_{1, 2}$ and $a_{2, 1}$ must also be zero.
Then  the polynomial $P(\zeta)$ introduced in Subsection 
 \ref{subsec:localobst} is
\[
a_{3, 1}(\zeta-p_2)+a_{3, 2}\zeta
 = (a_{3, 1}+a_{3, 2})\zeta
  -a_{3, 1}p_2.
\]

Thus, for the given configuration $\varphi_0$ of $C_0$, 
 an element of the dual obstruction space $H$ is given by 
 the data $\{a_{i, j}\}$ which satisfy
\[
a_{3, 1}+a_{3, 2} = 0
\]
 and
\[
a_{3, 1}p_2 = 0.
\]
Since $p_2$ is not zero, these conditions 
 imply that the space $H$ is $\{0\}$.
In other words, if there is a pre-log curve of the type given by
 the above tropical curve, then 
 its obstruction automatically vanishes, though the tropical curve
 itself is superabundant. \qed

\end{example}

\begin{rem}
The existence of pre-log curves of type $(\Gamma, h)$ of genus one and 
 the smoothability of those curves will be generalized much further 
 in Subsection \ref{subsec:genusonesm}.
\end{rem}

\subsubsection{Well-spacedness condition for higher genus curves with higher valent vertices}
Higher valent vertices are indispensable in the study of higher genus
 curves.
We give an example exhibiting a new phenomenon caused by a higher valent
 vertex which is important to the study of higher genus curves.

\begin{example}\label{ex:3}
 We consider 
 a tropical curve $(\Gamma, h)$ of genus two in $\Bbb R^3$,
 see Figure \ref{fig:00C}.

\begin{figure}[h]
\includegraphics[height = 6cm]{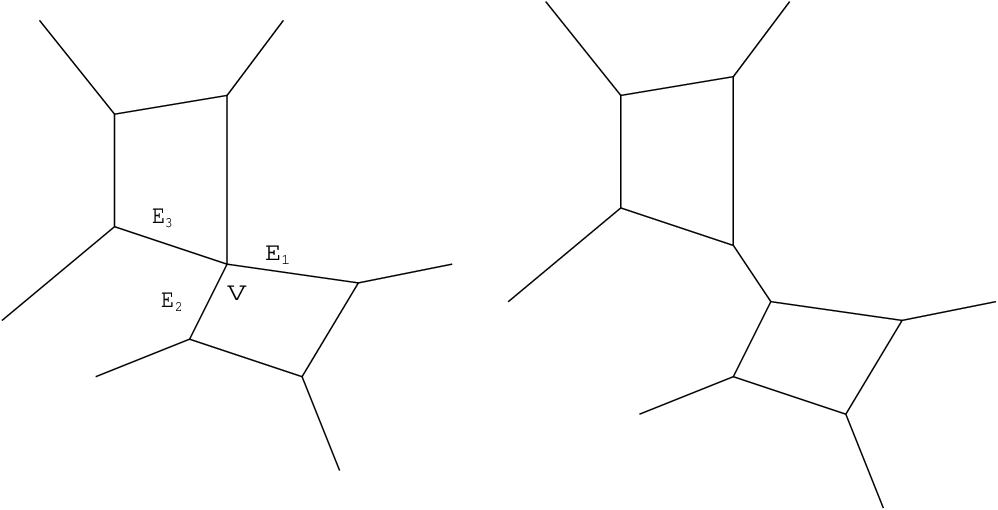}
\caption{A tropical curve of genus two with a 4-valent vertex 
 $(\Gamma, h)$ (the picture on the left), 
 and its deformation (the picture on the right).}\label{fig:00C}
\end{figure}

Here the two loops are contained in different planes.
Assuming the existence of a pre-log curve of type $(\Gamma, h)$
 in the central fiber of a suitable toric degeneration, 
 we calculate the dual obstruction space $H$.

If the 4-valent vertex $V$ is resolved into
 two 3-valent vertices (see Figure \ref{fig:00C}), 
 then by Theorem \ref{thm:obstruction}, the space $H$ is two dimensional.
We now see the effect of the existence of the 4-valent vertex
 to the space $H$.

As in the previous example, we associate the numbers $a_{i, j}$, 
 $i, j = 1, 2, 3$ to the edges $E_1, E_2, E_3$ emanating from $V$.
We define the numbers $a_{i, j}$ (see (4) in Subsection \ref{subsec:localobst})
 using the pull back of the component $C_{0, V}$ of $C_0$ to a linear curve in $\Bbb P^3$.
The equation 
 $P(\zeta) = 0$ in Subsection \ref{subsec:localobst} becomes 
\[
(a_{1, 2}+a_{2, 1})(\zeta-p_3)+
 (a_{2, 3}+a_{3, 2})(\zeta - p_1)+(a_{3, 1}+a_{1, 3})(\zeta - p_2) = 0.
\]
As before, we take the coordinate $\zeta$ on the component of the degenerate curve
 corresponding to the vertex $V$ so that 
 $p_1 = 0$.
Then the equation implies the conditions
\[
a_{1, 2}+a_{2, 1}+a_{2, 3}+a_{3, 2}+a_{3, 1}+a_{1, 3} = 0,\;\;
(a_{3, 1}+a_{1, 3})p_2 + (a_{1, 2}+a_{2, 1}) p_3 = 0.
\]
Theorem \ref{thm:obstruction2} implies further conditions
\[
a_{1, 2} = a_{2, 1} = 0,\;\; a_{1, 3}+a_{2, 3} = 0.
\]
Since $p_2$ is nonzero, these equations imply
\[
a_{1, 3} = a_{3, 2} = -a_{3, 1} = -a_{2, 3}
\]
 and the other $a_{i, j}$s are zero.
These numbers can be extended to the whole loop part.
Thus, in this case the space $H$ is one dimensional and 
 the generator of it has support on the whole loop part, contrary to the
 resolved curve (the picture on the right of Figure \ref{fig:00C}. 
 In that case the space $H$ is two dimensional and 
 the generators are supported on the individual loops).
Thus, the 4-valent vertex cancels part of the dual obstruction, and
 mixes up the rest.

Using the calculation in \cite{N3}, we can deduce 
 the criterion for the smoothability of 
 tropical curves of genus two which has $(\Gamma, h)$ as a subgraph.
Namely, extending the well-spacedness condition for genus one
 curves (see \cite{N3, S}), a curve in the following Figure \ref{fig:00D}
 is smoothable.

\begin{figure}[h]
\includegraphics[height = 6cm]{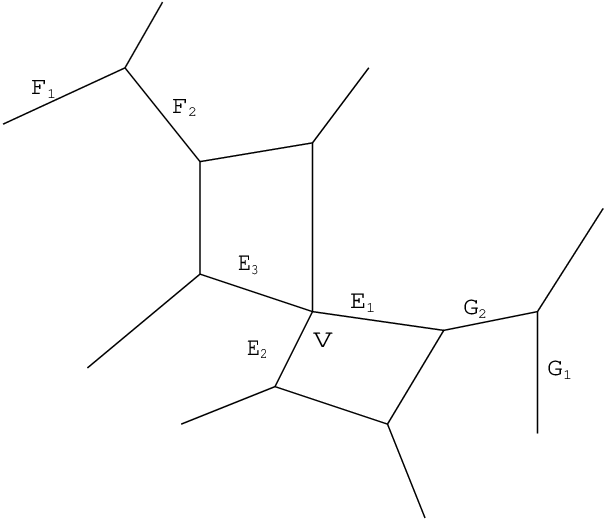}
\caption{}\label{fig:00D}
\end{figure}

Here the edge $F_1$ is not contained in the affine subplane spanned by 
 the edges in the loop containing $E_3$, and
 the edge $G_1$ is not contained in the affine subplane spanned by the
 edges in the loop containing $E_1, E_2$.
Moreover, the integral lengths of the edges 
 $F_2$ and $G_2$ should be the same. 
The case when the lengths of $F_2$ and $G_2$ are infinity
 (that is, the picture on the left of Figure \ref{fig:00C})
 is also smoothable.\qed

\end{example}

\subsection{Smoothability of tropical curves whose support of superabundancy has genus one}\label{subsec:genusonesm}
In this subsection 
 we prove a general result 
 concerning the smoothability of
 a pre-log curve corresponding to tropical curves of genus one
 which have 
 higher valent vertices at the loop.
More generally, it applies to tropical curves whose supports of superabundancy
 are of genus one (see Definition \ref{def:abundancysupport}).
This smoothability criterion has a very different form compared to the
known criterion (the well-spacedness condition, see \cite{N3, S}.
A hint of it can be found in \cite{Katz, Tor}).
This plays a fundamental role in the study of 
 algebraic curves through degeneration.
See for example \cite{N3}.

\begin{thm}\label{thm:loop}
	Let $h\colon \Gamma\to\Bbb R^n$ be a tropical curve satisfying Assumption A.
    Suppose its support of superabundancy $L$ is of genus one.
	We assume the following.
	\begin{itemize}
		\item Let $\{F_j\}$ be the set of flags of $(\Gamma, h)$
		whose vertices are on $L$.
		There is a natural direction $v_j$ of $\Bbb R^n$ 
		determined by each $F_j$.
		Then these directions span $\Bbb R^n$.
	\end{itemize}
	Then, if $\varphi_0\colon C_0\to X_0$ is a pre-log curve of type
	$(\Gamma, h)$ in the central fiber of a suitable toric degeneration, 
	the curve $\varphi_0$ is smoothable.
\end{thm} 
\begin{rem}
	Note that this theorem does not insist anything on the existence of 
	a pre-log curve of type $(\Gamma, h)$, and
	the existence of such a curve is an assumption. 
	As Examples \ref{ex:1} and \ref{ex:2} show, the existence of 
	such a curve is a nontrivial problem.
In the case of genus one, we have a solution to this problem,
 see the next section.
\end{rem}
\proof
We prove that under the assumption of the theorem, 
the space $H$ of the dual obstructions is $\{0\}$.
Then the obstruction to deform $\varphi_0$ automatically vanishes, 
and it follows that $\varphi_0$ is smoothable.

Let $V_1, \dots, V_l$ be the higher valent vertices on $L$
 and take one of them.
We write it by $V_a$. 
Let $E_{a, 1}, \dots, E_{a, k_a}$ be the edges of $h(\Gamma)$
 emanating from $V_a$, and let $E_{a, 1}$ and $E_{a, 2}$
 be the edges contained in $L$.
We write by $B_a$ the subspace of $\Bbb R^n$ 
 spanned by the directions
 of the edges  $E_{a, 1}, \dots, E_{a, k_a}$.
%Let $b_a$ be the dimension of $B_a$.
%Take a basis $\{e_{a, 1}, \dots, e_{a, n}\}$ of $\Bbb Z^n\subset \Bbb R^n$ so that 
%the set of the first $b_a$ vectors $\{e_{a, 1}, \dots, e_{a, b_a}\}$
%is a basis of $B_a$.

By the argument in Subsections \ref{subsec:localobst} and 
 \ref{subsec:gendual}, we attach numbers
 $\{a_{i, j}\}$ and $\{c_{k, l}\}$ to the vertex $V_a$ as in 
 the expression (\ref{eq:localobst}) in Subsection \ref{subsec:localobst}.
%Here the tropical curve with only one vertex $V_a$
%is obtained as the image of the standard tropical curve with one vertex
%by 
%composing a linear map
%\[
%\Xi\colon \Bbb R^{k_a-1}\to B_a\subset \Bbb R^n, 
%\]
%see Subsection  \ref{subsec:gendual}.
%Note that the vectors in the dual basis
%$\{e_{a, 1}^{\vee}, \dots, e_{a, n}^{\vee}\}$ can be seen as 
%vectors in the dual space of $\Bbb R^{k_a-1}$ 
%through the dual of this map (the vectors $e_{a, b_a+1}^{\vee}, \dots, 
%e_{a, n}^{\vee}$ are regarded as zero vectors, 
%which reflects the fact that there are no restrictions to the values of
%$\{c_{k, l}\}$).

By Theorem \ref{thm:obstruction2}, it is easy to see that
the numbers $a_{i, j}$ and $c_{k, l}$ attached to the edges
except $E_{a, 1}$ and $E_{a, 2}$ must be zero.
In particular, 
\[
a_{i, j} = 0,\;\; j\geq 3.
\]
The equation $P(\zeta)=0$ (see Subsection \ref{subsec:localobst})
 becomes into the form
\[
\sum_{j = 1, 2}\sum_{i=1}^{k_a-1}a_{i, j}\prod_{l\neq i, j}(\zeta-p_l)=0,
\]
where $\zeta$ is a coordinate on $C_{0, V_a}$ which is zero at the point
corresponding to the edge $E_{a, 1}$.
Substituting $\zeta = p_i$, $i = 3, \dots, k_a-1$, we have 
a series of equations
\[
(p_i-p_2)a_{i, 1}+ p_ia_{i, 2} = 0,\;\; i = 3, \dots, k_a-1.
\]
%Under these conditions, substituting $\zeta = p_1$ and $\zeta = p_2$
% gives the same condition
%\[
%\frac{a_{1, 2}+a_{2, 1}}{p_2} = -\frac{a_{3, 1}}{p_3}-\cdots - 
% \frac{a_{k_i-1, 1}}{p_{k_i-1}}.
%\]

On the other hand, by the residue theorem we also have the equations
\[
\sum_{j=1}^{k_a-1} a_{i, j} = 0,\;\; i = 1, \dots, k_a-1.
\]
Note that $\zeta = \infty$ corresponds to the edge $E_{a, k_a}$
 and the residue there is zero by the same reason as
 $a_{i, j} = 0$ for $j\geq 3$ above.

It follows that $a_{1, 2} = a_{2, 1} = 0$.
Also, since $p_i$, $i = 2, \dots, k_a-1,$ are nonzero, 
the equation $(p_i-p_2)a_{i, 1}+ p_ia_{i, 2} = 0$
cannot be proportional to $a_{i, 1}+a_{i, 2} = 0$, $i = 3, \dots, k_a-1$.
Thus, all the numbers $a_{i, j}$ must be zero.
This is imposed at all higher valent vertices on $L$.

Therefore, the only nonzero terms attached to higher valent vertices are 
$c_{i, j}$ associated to the directions in $N_{\Bbb R}^{\vee}$ which
annihilate the subspace $B_a$ of $N_{\Bbb R}$ defined above.
Then, by Theorem \ref{thm:obstruction2}, it follows that 
the space $H$ is naturally contained in the space $N_{\Bbb R}^{\vee}$
and a vector in $H$ must annihilate all the directions of the edges
attached to a vertex on $L$.
By assumption, these directions span $N_{\Bbb R}$, 
 and the space $H$ must be zero.\qed
\begin{rem}\label{rem:supple}
	If the directions of the edges emanating from the vertices on $L$ 
	do not span $N_{\Bbb R}$, the space $H$ is identified with the 
	annihilator subspace of the space spanned by these directions.
	The vanishing of the remaining obstruction is then argued by the method
	in \cite{N3}.
\end{rem}

The proof actually shows the following.
Namely, let $(\Gamma, h)$ be a tropical curve satisfying Assumption A.
Let $V$ be a vertex of $h(\Gamma)$ with valence $k\geq 4$.
We attach numbers $\{a_{i, j}\}$ and $\{c_{i, j}\}$
 to the edges emanating from $V$ as in 
 (\ref{eq:localobst}) in Subsection \ref{subsec:localobst}. 
\begin{cor}
If these numbers except the ones attached to 
 two of the edges vanish, then these numbers entirely vanish.\qed
\end{cor}
This is a statement which does not depend on the genus of the curve, and often reduces 
 the dimension of the (dual) obstruction space of curves of higher genus.

The same argument applies to the cases where 
 the numbers attached to more than two edges are (possibly) nonzero.
Let $E_1, E_2, E_3, \dots, E_k$ be the edges emanating from $V$
 and assume that all the numbers attached to the edges 
 except $E_1, E_2, E_3$ are zero (for example, assume that the edges $E_4, \dots, E_k$
 are all unbounded).
Then we have the following.
\begin{prop}
The numbers $\{a_{i, j}\}$ satisfy the relations
\[
\begin{array}{l}
a_{1, 2}+a_{1, 3} = 0,\\
a_{2, 1}+a_{2, 3} = 0,\\
a_{3, 1}+a_{3, 2} = 0,\\
a_{i, 1}+a_{i, 2}+a_{i, 3} = 0, \;\; (i = 4, 5, \dots, k-1)\\
\frac{a_{i, 1}}{p_i}+\frac{a_{i, 2}}{p_i-p_2}+\frac{a_{i, 3}}{p_i-p_3} = 0, \;\; (i = 4, 5, \dots, k-1),
\end{array}
\]
 here we used the same notation as in the proof of the theorem.
In this case, these cover all the local relations imposed on the numbers 
 $\{a_{i, j}\}$. \qed
\end{prop}
Local means that they are the consequence of Lemma \ref{lem:cond2}.
 In general, relations given by Theorem \ref{thm:obstruction2} give more constraints.

The relations of the same form also hold in the case where
 there are more than three edges
 on which the attached numbers can be nonzero. 
However, in such cases these relations do not cover all the constraints
 given by Lemma \ref{lem:cond2}, and there will be more complicated relations.
Nevertheless, as Example \ref{ex:3} shows, in many cases the relations
 from Theorem \ref{thm:obstruction2} reduce the freedom of the numbers
 and the result becomes simple.

\section{Pre-log curves of genus one and lattice point count}\label{sec:6}
Theorem \ref{thm:loop} shows that a part of the obstruction 
vanishes on degenerate algebraic curves corresponding
 to superabundant tropical curves of genus one
 when there are higher valent vertices on the loop.
On the other hand, as Examples \ref{ex:1} and \ref{ex:1.5} show, 
when the tropical curve is superabundant, in general there are cases
where corresponding pre-log curves do not exist and understanding
 when there are such curves is a rather complicated problem.
However, when the curve is of genus one, this problem can be
solved, as we will see below.

Let $(\Gamma, h)$, $h\colon \Gamma\to \Bbb R^n = N\otimes\Bbb R$
be a tropical curve of genus one satisfying Assumption A.
Let $\tilde L\subset h(\Gamma)$ be the union of all the open edges emanating from 
vertices on the loop and the loop itself (in other words, $\tilde L$ is the
'open star' of the loop part $L$ of $h(\Gamma)$).

The part $\tilde L$ can be thought of as (the image of) a tropical curve by extending the 
open edges to infinity.
We write it by $(\Gamma', h')$.
If there is a pre-log curve of type $(\Gamma', h')$, it is easy to 
see that there is also a pre-log curve of type $(\Gamma, h)$, 
since $\Gamma\setminus \Gamma'$ is a union of trees.
Thus, we assume there are no vertices outside $L$, so that
$h(\Gamma)$ equals $\tilde L$.
We also assume that the image $h(\Gamma)$ is not contained
in a proper subspace of $\Bbb R^n$.

Under these assumptions, if $h$ is an immersion, then 
it is easy to see that the tropical curve $(\Gamma, h)$ is
regular, and there is a pre-log curve of 
type $(\Gamma, h)$ by Proposition \ref{prop:exist}. 
Therefore, the problem is the case where higher valent vertices exist. 

Let $A$ be the minimal affine subspace of $N_{\Bbb R}\cong \Bbb R^n$ containing 
the loop of $h(\Gamma)$.
Let $\bar A$ be the linear subspace of $N_{\Bbb R}$ parallel to $A$.
Take a basis $\{e_1, \dots, e_n\}$ of $N$ so that 
$\{e_1, \dots, e_r\}$ is a basis of $\bar A$.
Let $\{e^{\vee}_1, \dots, e^{\vee}_n\}$ be the dual basis of $\{e_1, \dots, e_n\}$.

Let $\mathfrak X$ be a toric degeneration defined respecting $(\Gamma, h)$
(see Definition \ref{def:degeneration}) and $X_0$ be its central fiber.
A pre-log curve $\varphi_0\colon C_0\to X_0$, if it exists, can be obtained by
gluing irreducible components $\varphi_{0, v}\colon C_{0, v}\to X_{0,v}$, 
where $v$ is a vertex of $h(\Gamma)$. 

Recall that the degeneration $\mathfrak X$ is defined
by a fan in $(N\oplus\Bbb Z)\otimes\Bbb R$
having non-negative part with respect to the $\Bbb Z\otimes\Bbb R$ summand.
Let $u$ be the positive generator of the $\Bbb Z$-summand of $N\oplus \Bbb Z$ and
consider the dual basis $ 
e^{\vee}_1, \dots, e^{\vee}_{n}, u^{\vee}$ of the basis
$e_1, \dots, 
e_{n}, u$ 
of $N\oplus\Bbb Z$.
For any vertex $v$ of $h(\Gamma)$, 
and any vector $w$ in $N^{\vee}$,  
there is a function on $\mathfrak X$ corresponding to 
the vector $w+j u^{\vee}$ 
for a unique integer $j$ (which depends on $v$) such that
it is not constantly 0 or $\infty$ on the component $X_{0,v}$
of the central fiber $X_0$ of $\mathfrak X$ corresponding to the 
vertex $v$.
We call this the function corresponding to 
$w$ on $X_{0, v}$.

Since the loop of $h(\Gamma)$ is contained in the affine subspace $A$, 
it follows that there are functions $Z_{r+1}, \dots, Z_n$ on $\mathfrak X$
corresponding to 
$e^{\vee}_{r+1}, \dots, e^{\vee}_n$ with the property that
 any of which is not constantly zero or 
 diverges on the whole component $X_{0,v}$.
Here $v$ is any vertex on the loop of $h(\Gamma)$.

For such a vertex $v$, let $\varphi_{0, v}\colon C_{0, v}\to X_{0,v}$ be 
a pre-log curve of type $v$ in the sense of Definition \ref{def:typev}.
The star of $v$ in $h(\Gamma)$ has two edges contained in the loop 
of $h(\Gamma)$.
We write them by $E_{v, 1}$ and $E_{v, 2}$.
Let $p_{v, 1}$ and $p_{v, 2}$ be the points on $\varphi_{0, v}(C_{0, v})$
corresponding to these edges.
In particular, these points are contained in toric divisors of $X_{0, v}$.

Let $\{v_1, \dots, v_k\}$ be the vertices on the loop of $h(\Gamma)$
ordered cyclically.
The vertex $v_1$ is also referred to as $v_{k+1}$.
We assume the edges $E_{v_i, 2}$ and $E_{v_{i+1}, 1}$ coincide.

\begin{lem}\label{lem:loopcomp}
	If there are pre-log curves $\varphi_{0, v_i}\colon C_{0, v_i}\to X_{0,v_i}$
	such that 
	\[
	Z_j(p_{v_i, 2}) = Z_j(p_{v_{i+1}, 1}),\;\; \forall j\in \{r+1, \dots, n\},\;\;\forall i\in\{1, \dots, k\}, 
	\]
	then there is a pre-log curve of type $(\Gamma, h)$.
\end{lem}
\proof
Let us take the edge $E$ of the loop of $h(\Gamma)$ whose ends are
$v_1$ and $v_2$.
Cutting it at the midpoint, extending the cut edges to infinity, and
adding suitable vertices as in Figure \ref{fig:cut_and_glue}
of Subsection \ref{subsec:corrthm}, we obtain a tropical curve of genus zero.
Let $(\Gamma', h')$ be the resulting tropical curve. 

Since $\Gamma'$ is a tree, we can construct a pre-log curve 
$\varphi_0'\colon C_0'\to X_0$ of type $(\Gamma', h')$
In particular, we have two vertices corresponding to the midpoint of $E$, 
and these in turn correspond to two orbits of a one dimensional subtorus 
the torus $N\otimes \Bbb C^*$ acting on components of $X_0$.
Moreover, by the assumption of the lemma, we can take these orbits so that
their values of the functions $Z_{r+1}, \dots, Z_n$ are the same
(note that each $Z_i$, $i = r+1, \dots, n$, is constant on these orbits).
Therefore, the positions of these two orbits are different only in the direction
of $\bar A\otimes \Bbb C^*$.

On the other hand, 
by the argument as in Subsection \ref{subsec:corrthm},
we can change the relative position of these one dimensional orbits 
in any direction of $\bar A\otimes \Bbb C^*$ by 
moving the components of $\varphi_0'(C_0')$ by the action of 
the torus in the direction of $\bar A\otimes \Bbb C^*$, without changing the 
homotopy type of $\varphi_0'(C_0')$.
Eventually, we can take $\varphi_0'$ so that the orbits corresponding to the 
two vertices at the midpoint of $E$ coincide.
Then again as in the argument in Subsection \ref{subsec:corrthm},
we can discard some components of $C_0'$ so that 
we have a pre-log curve of type $(\Gamma, h)$ 
(see also Figure \ref{fig:cut_and_glue}).\qed\\

Now our problem is reduced to determining whether there are pre-log curves of 
type $v_i$ satisfying the assumption of Lemma \ref{lem:loopcomp}
for each $j$.
We will interpret this condition by looking at each curve 
 corresponding to $v_i$ closely to classify the possible values
 of the functions $Z_j$, $r+1\leq j\leq n$, 
 at the points $p_{v_i, 1}$ and $p_{v_i, 2}$.

Let $v$ be a vertex on the loop of $h(\Gamma)$.
The open star of $v$ in $h(\Gamma)$ can be thought of as
a tropical curve $(\Gamma_v, h_v)$ whose image has only one vertex.
Here $\Gamma_v$ is a subgraph of $\Gamma$ and $h_v$ 
is the map induced from $h$ in a natural way. 
As we discussed in Subsection \ref{subsec:gendual}, such a tropical curve
can be obtained from a standard tropical curve 
$h_0\colon\Gamma_0\to \Bbb R^{l+1}$ introduced
at the beginning of Subsection \ref{subsec:localobst} 
by composing $h_0$ with a suitable integral affine linear map 
$\Xi\colon\Bbb R^{l+1}\to \Bbb R^n$.
Here the valence of $h(\Gamma)$ at the vertex $v$ is $l+2$.
Let $\{\varepsilon_1, \dots, \varepsilon_{l+1}\}$ be the standard basis of
$\Bbb Z^{l+1}$.

Correspondingly, any pre-log curve of type 
$v$ (see Definition \ref{def:typev}) can be obtained from 
a pre-log curve $\psi_0\colon D_0\to \Bbb P^{l+1}$
of type $(\Gamma_0, h_0)$ by composing $\psi_0$ with
a map between toric varieties $f_{\Xi}\colon \Bbb P^{l+1}\to \Bbb P_v$ induced by $\Xi$.
Here $D_0$ is an irreducible nonsingular rational curve, and $\Bbb P_v$
 is the toric variety defined by $h_v(\Gamma_v)$ seen as a fan, 
 see Subsection \ref{subsec:gendual}. 

Now let $u_{v, 1}, \dots, u_{v, l+2}\in N$ be the primitive integral 
generators of the edges
emanating from $v$ in $h(\Gamma)$ (some of these may coincide).
We assume that $u_{v, 1}$ and $u_{v, 2}$ are the directions of the edges 
contained in the loop of $h(\Gamma)$.
Let $w_{v, 1}, \dots, w_{v, l+2}$ be the weights of these edges.
Then the linear part of the map $\Xi$ is given by the $n\times (l+1)$-matrix
\[
G = (w_{v, 1}u_{v, 1}\ w_{v, 2}u_{v, 2} \ \cdots \ w_{v, l+1}u_{v, l+1}).
\] 

The dual vector $e^{\vee}_j\in N^{\vee}$ introduced above is sent to the vector
\[
(w_{v, 1}u_{v, 1, j}, w_{v, 2}u_{v, 2, j}, \dots, w_{v, l+1}u_{v, l+1, j})\in (\Bbb R^{l+1})^{\vee}
\]
 by the adjoint $G^T$ of $G$,  
here $u_{v, i, j}$ 　is the $j$-th component of the vector
$u_{v, i}\in\Bbb Z^n$.
Note that since the part $\varepsilon_1, \varepsilon_2$ of the standard generators 
of $\Bbb Z^{l+1}$ are mapped to
$w_{v, 1}u_{v, 1}$ and $w_{v, 2}u_{v, 2}$, respectively, the equations
\[
u_{v, 1, j} = u_{v, 2, j} = 0,\;\; r+1\leq j\leq n
\]
hold.

Now the image of the curve $\psi_0$ is given by the defining equations
\[
x_1+a_2x_2+b_2 = 0,\;\; x_1+a_3x_3+b_3 = 0,\;\; \dots,\;\; x_1+a_{l+1}x_{l+1}+b_{l+1} = 0, 
\] 
here $x_i$ are inhomogeneous coordinates of $\Bbb P^{l+1}$ 
 corresponding to the dual basis
 $\{\varepsilon_1^{\vee}, \dots, \varepsilon_{l+1}^{\vee}\}$
 of $\{\varepsilon_1, \dots, \varepsilon_{l+1}\}$,
 and $a_i, b_i$ are nonzero complex numbers satisfying 
\[
b_i\neq b_j,\;\; \forall i\neq j.
\]

Let $S$ be the coordinate on $D_0$ given by $\psi_0^*x_1$.
Then the function $Z_j$ on $\Bbb P_v$, $j \in\{ r+1, \dots, n\}$, corresponding to the vector 
 $e^{\vee}_j$
is pulled back to 
\[
(\clubsuit) \;\;
\psi_0^*(x_1^{w_{v, 1}u_{v, 1, j}}
%x_2^{w_{v, 2}u_{v, 2, j}}
\cdots x_{l+1}^{w_{v, l+1}u_{v, l+1, j}})
= \left(-\frac{b_3+S}{a_3}\right)^{w_{v, 3}u_{v, 3, j}}
%\left(-\frac{b_4+S}{a_4}\right)^{w_{v, 4}u_{v, 4, j}}
\cdots  \left(-\frac{b_{l+1}+S}{a_{l+1}}\right)^{w_{v, l+1}u_{v, l+1, j}}
\] 
 by the map $f_{\Xi}\circ \psi_0$.
 
Let $p_{v, 1}$ and $p_{v, 2}$ be the points
 on $D_0$ corresponding to the edges of directions
 $\varepsilon_1, \varepsilon_2$ as before.
The values of the coordinate $S$ at these points are $0$ and $-b_2$, 
respectively.
Therefore, the values of the function $Z_j$ at $p_{v, 1}$ and $p_{v, 2}$ are given by
\[
Z_j(p_{v, 1}) = \left(-\frac{b_3}{a_3}\right)^{w_{v, 3}u_{v, 3, j}}
\left(-\frac{b_4}{a_4}\right)^{w_{v, 4}u_{v, 4, j}}
\cdots  \left(-\frac{b_{l+1}}{a_{l+1}}\right)^{w_{v, l+1}u_{v, l+1, j}}
\]
 and
\[
Z_j(p_{v, 2}) = \left(-\frac{b_3-b_2}{a_3}\right)^{w_{v, 3}u_{v, 3, j}}
\left(-\frac{b_4-b_2}{a_4}\right)^{w_{v, 4}u_{v, 4, j}}
\cdots  \left(-\frac{b_{l+1}-b_2}{a_{l+1}}\right)^{w_{v, l+1}u_{v, l+1, j}}, 
\]
 respectively.
The ratio of these values is
\[
\frac{Z_j(p_{v, 2})}{Z_j(p_{v, 1})} = 
 \left(\frac{b_3-b_2}{b_3}\right)^{w_{v, 3}u_{v, 3, j}}
 \left(\frac{b_4-b_2}{b_4}\right)^{w_{v, 4}u_{v, 4, j}}
 \cdots  \left(\frac{b_{l+1}-b_2}{b_{l+1}}\right)^{w_{v, l+1}u_{v, l+1, j}}.
\]

Recall that $\{v_1, \dots, v_k\}$ are the vertices on the loop of $h(\Gamma)$
 ordered cyclically.
We take pre-log curves 
\[
\varphi_{0, v_i}\colon D_{0, v_i}\to X_{0, v_i}
\]
 of type $v_i$ for each $i = 1, \dots, k,$ and also take points $p_{v_i, 1}$
 and $p_{v_i, 2}$ corresponding to the edges $E_{v_i, 1}$ and $E_{v_i, 2}$
 (recall that $E_{v_i, 2} = E_{v_{i+1}, 1}$. It follows that the edges
 $E_{v_k, 2} = E_{v_1, 1}, E_{v_1, 2} = E_{v_2, 1}, E_{v_2, 2} = E_{v_3, 1}, \cdots$
 are aligned cyclically).
The following is clear.
\begin{lem}\label{lem:prod1}
The condition of Lemma \ref{lem:loopcomp} is equivalent to the condition that we can take
 the pre-log curves $\varphi_{0, v_i}$ so that 
 the product of the ratios
\[
\prod_{i=1}^k\frac{Z_j(p_{v_i, 2})}{Z_j(p_{v_i, 1})}
\]
 is equal to one for all $j = r+1, \dots, n$.\qed
\end{lem}

Now consider the ratio
\[
\frac{Z_j(p_{v_i, 2})}{Z_j(p_{v_i, 1})} 
%= 
%\left(\frac{b_3-b_2}{b_3}\right)^{w_3u_{3, j}}\left(\frac{b_4-b_2}{b_4}\right)^{w_4u_{4, j}}
%\cdots  \left(\frac{b_{l+1}-b_2}{b_{l+1}}\right)^{w_{l+1}u_{l+1, j}}
 = \left(1-\frac{b_2}{b_3}\right)^{w_{v_i, 3}u_{v_i, 3, j}}
 \left(1-\frac{b_2}{b_4}\right)^{w_{v_i, 4}u_{v_i, 4, j}}
 \cdots  \left(1-\frac{b_2}{b_{l_i+1}}\right)^{w_{v_i, l_i+1}u_{v_i, l_i+1, j}},
\]
 here the valence of the vertex $v_i$ is $l_i+2$.
%For the curve $\varphi_{0, v_i}$ to be torically transverse, the non-zero constants
% $b_2, b_3, \dots, b_{l_i+1}$ need to be all different.
Thus, by changing the configuration of the curve $\psi_0$, 
 the ratio $\frac{Z_j(p_{v_i, 2})}{Z_j(p_{v_i, 1})}$ varies over the set
\[
 \{C_{v_i, j} := \zeta_{v_i, 3}^{w_{v_i, 3}u_{v_i, 3, j}}
 \zeta_{v_i, 4}^{w_{v_i, 4}u_{v_i, 4, j}}\cdots 
 \zeta_{v_i, l_i+1}^{w_{v_i, l_i+1}u_{v_i, l_i+1, j}}\;|\;
 \zeta_{v_i, q}\neq 1 \;\;\forall q,\;\; \zeta_{v_i, q_1}\neq \zeta_{v_i, q_2}\;\; 
  3\leq \forall q_1<\forall q_2\leq l_i+1\}.
\]
Let $\mathcal H_i$ be the set
\[
\mathcal H_i = \{(\zeta_{v_i, 3}, \zeta_{v_i, 4}, \dots, \zeta_{v_i, l_i+1})\in \Bbb C^{l_i-1}\;|\; 
 \zeta_{v_i, q}\neq 1 \;\;\forall q,\;\; \zeta_{v_i, q_1}\neq \zeta_{v_i, q_2}\;\; 
 3\leq \forall q_1<\forall q_2\leq l_i+1\}.
\]

Using this description, the condition of Lemma \ref{lem:prod1} can be written as follows.
\begin{lem}\label{lem:prodmap}
The condition of Lemma \ref{lem:prod1} is equivalent to the condition that
 the image of the product map
\[
\prod_{i=1}^k\mathcal H_i\to (\Bbb C^*)^{n-r},\;\;
 (\zeta_{v_i, q})\mapsto \left(\prod_{i=1}^k C_{v_i, j}\right)_{j=r+1, \dots, n}
\] 
 contains $(1, 1, \dots, 1)$.\qed
\end{lem}

Now note that if the direction of the edge $E_{v_i, q}$ is 
  contained in the linear subspace $\bar A$ spanned by the 
  directions of the edges in the loop, then 
  $u_{v_i, q, j} = 0$ for all $r+1\leq j\leq n$.
Therefore, we only need to consider those edges not contained in the affine
 subspace $A$.

We can assume that we have numbered the edges emanating from $v_i$ so that 
\[
E_{v_i, 3}, \dots, E_{v_i, m_i},\;\; 1\leq i\leq k,\;\; 4\leq m_i\leq l_i+2
\]
 are those not contained in $A$.
When all the edges emanating from $v_i$ are contained in $A$, 
 we set $m_i$ to be three.
Note that when there is an edge emanating from $v_i$ not contained in $A$, 
 then $m_i\geq 4$ by the balancing condition.
We choose $m_i-3$ edges $E_{v_i, 3}, \dots, E_{v_i, m_i-1}$
 for each $i = 1, \dots, k$.
Let $\mathcal E$ be the set of these edges:
\[
\mathcal E = \{E_{v_i, q_i}\}_{1\leq i\leq k, 3\leq q_i\leq m_i-1}
\]

\begin{rem}
There is freedom of choices of $m_i-3$ edges from 
 the $m_i-2$ edges $E_{v_i, 3}, \dots, E_{v_i, m_i}$
 for each vertex $v_i$, $1\leq i\leq k$.
However, the conclusion does not depend on these choices.
\end{rem}

Then we construct a flag of affine subspaces of $N_{\Bbb R}$
 containing $A$ in the following way.
First choose an edge $E_{v_i, q}$ from $\mathcal E$.   
Let $A_{r+1}$ be the minimal affine subspace containing both $A$ and $E_{v_i, q}$.
Let $n_{r+1}$ be the number of the edges in $\mathcal E$ contained in $A_{r+1}$.
Next, choose an edge $E_{v_j, s}$ from $\mathcal E$ which is not contained
 in $A_{r+1}$.
Let $A_{r+2}$ be the minimal affine subspace containing both $A_{r+1}$ and $E_{v_j, s}$.
Let $n_{r+2}$ be the number of the edges in $\mathcal E$ contained in $A_{r+2}$.

Continuing this process, we can construct a flag
\[
A = A_r\subset A_{r+1}\subset \cdots \subset A_{n-1}\subset A_{n} = N_{\Bbb R}.
\]
There is also an associated sequence of the numbers
\[
n_{r+1}<n_{r+2}<\cdots < n_{n-1}<n_{n} = \sum_{i=1}^k(m_i-3).
\]

Among these flags, we choose one for which the sequence of the positive integers
\[
n_{r+1}, n_{r+2}-n_{r+1}, \dots, n_{n}-n_{n-1}
\]
 is minimal with respect to the lexicographic order
 (in other words, we choose a flag so that $n_{i+1}-n_i$
 with smaller $i$ will be as small as possible). 
There may be several flags which give the same sequence of the integers, 
 and any of them will do.
We write such a flag by   
 $A_r\subset A_{r+1}\subset \cdots \subset A_{n-1}\subset A_{n}$. 
  
Now choose a basis $\{g_{r+1}, \dots, g_n\}$ of the sublattice 
 $N^{\vee}\cap \bar A^{\perp}$ 
 of the
 dual lattice $N^{\vee}$  
 in the following way.
First, let $g_{n}$ be one of the generators of the annihilating subspace
 of $\bar A_{n-1}$. 
Next, take $g_{n-1}$ so that $\{g_{n-1}, g_{n}\}$ will be a basis of the annihilating
 subspace of $\bar A_{n-2}$.
Similarly, take $\{g_{r+1}, \dots, g_{n}\}$ so that 
 $\{g_{n-i+1}, \dots, g_{n}\}$ is a basis of the annihilating subspace
 of $\bar A_{n-i}$.

The basis $\{g_{r+1}, \dots, g_{n}\}$ 
 can be obtained from the given basis 
 $\{e^{\vee}_{r+1}, \dots, e^{\vee}_n\}$ by an invertible integral linear transform.
Therefore, we can replace the functions $Z_{r+1}, \dots, Z_n$ in 
 Lemma \ref{lem:loopcomp} by those functions corresponding to 
 the vectors $\{g_{r+1}, \dots, g_{n}\}$.
We also write these functions by $Z_{r+1}, \dots, Z_n$.

Let us order the elements of the set of edges $\mathcal E$   
 as 
\[
E_1, E_2, \dots, E_{n_n}
\]
 so that the first $n_{r+j}$
 edges are contained in $A_{r+j}$, $j = 1, \dots, n-r$.
Let $v_m$ be the primitive integral generator of the 
 direction of the edge $E_m$ and $w_m$ be the weight of it.  
The values 
\[
w_mg_j(v_m), \;\; j = r+1, \dots, n
\]
 give an element of $\Bbb Z^{n-r}$ for 
 $m = 1, \dots, n_n$.   
Taking these vectors as column vectors, we make an
 $(n-r)\times n_n$-matrix $G'$.

Now among the positive integers
$n_{r+1}, n_{r+2}-n_{r+1}, \dots, n_n-n_{n-1}$, 
let $j\in \{0, 1, \dots, n-r-1\}$ be the smallest integer
such that $n_{r+j+1}-n_{r+j}$ is larger than 1
(here we set $n_r = 0$).
Such a $j$ may not exist, but when it exists, 
 then the following is easy to see.

\begin{lem}\label{lem:tuple}
For each $i\geq  j$, the inequality $n_{r+i+1}-n_{r+i}\geq 2$ holds. 
\end{lem}  
\proof
Suppose that $n_{r+i+1}-n_{r+i} =1$ holds for some $i\geq j$.
Recall we have the flag 
 $A_r\subset A_{r+1}\subset\cdots\subset A_{n-1}\subset A_n$
 of affine subspaces in $N_{\Bbb R}$.
By assumption, the affine subspace $A_{r+i+1}$ is the minimal affine subspace
 containing both $A_{r+i}$ and the edge $E_{n_{r+i}+1}$, and 
 $E_{n_{r+i}+1}$ is the unique edge of $\mathcal E$
 contained in $A_{r+i+1}\setminus A_{r+i}$.

Now take another flag $\{A'_{k}\}_{k = r, \dots, n}$ so that 
 $A'_k = A_k$ for $k = r, \dots, r+j-1$ and
 define $A'_{r+j}$ as the minimal affine subspace 
 containing $A_{r+j-1}$ and the edge $E_{n_{r+i}+1}$.
Then by construction $E_{n_{r+i}+1}$ is the unique edge of $\mathcal E$
 contained in $A'_{r+j}\setminus A_{r+j-1}$.
This contradicts to the minimality of $n_{r+1}, n_{r+2}-n_{r+1}, \dots, n_n-n_{n-1}$
 in the lexicographic order.\qed\\

Recall that each column of $G'$ corresponds to an edge in $\mathcal E$.
In turn, to each such an edge corresponds a complex number $\zeta_{v_i, q}$  
 for some $i$ and $q$.
Let $j\in\{0, 1, \dots, n-r-1\}$ be the minimal integer such that 
 $n_{r+j+1}-n_{r+j} \geq 2$ holds as above (assuming such an integer exists).
By Lemma \ref{lem:tuple}, any integer $a\in \{j, j+1, \dots, n-r-1\}$ satisfies
 this property.
Consider the product $\prod_{i=1}^k C_{v_i, n}$.
By assumption, it contains at least two factors of the form
 $\zeta_{v_i, q}^{\alpha}$, $\alpha\neq 0$.
Therefore, clearly there is a continuum of elements of $\prod_{i=1}^k\mathcal H_i$
 which satisfy $\prod_{i=1}^k C_{v_i, n} = 1$.

Now, if $j\leq n-r-2$, then the product $\prod_{i=1}^k C_{v_i, n-1}$ contains 
 at least two factors of the form $\zeta_{v_{i'}, q'}^{\alpha'}$ which do not appear in 
 $\prod_{i=1}^k C_{v_i, n}$.
Thus, there is a continuum of elements of $\prod_{i=1}^k\mathcal H_i$
 which satisfy both $\prod_{i=1}^k C_{v_i, n} = 1$ and
 $\prod_{i=1}^k C_{v_i, n-1} = 1$.

This can be continued so that we see there is a continuum of elements of $\prod_{i=1}^k\mathcal H_i$
 which satisfy
\[
\prod_{i=1}^k C_{v_i, n} = 
 \prod_{i=1}^k C_{v_i, n-1} = \cdots
 =\prod_{i=1}^k C_{v_i, r+j+1} = 1. 
\]
Let $\mathcal H'$ be the subset of $\prod_{i=1}^k\mathcal H_i$
 whose elements satisfy these equalities.

Now consider the product $\prod_{i=1}^k C_{v_i, l}$, $r\leq l\leq r+j$. 
If the $(l-r)$-th row of $G'$ contains a non-zero component at the $m$-th column, 
 where $j+1\leq m\leq n_n$, then it is easy to see that
 there is a continuum of elements in $\mathcal H'$ which satisfy
 $\prod_{i=1}^k C_{v_i, l} = 1$.

Then consider another product $\prod_{i=1}^k C_{v_i, l'}$.
If the $(l'-r)$-th row of $G'$ contains a non-zero component at the $m'$-th column, 
 where $j+1\leq m'\leq n_n$ or $m' = l-r$, 
 then again there is a continuum of elements in $\mathcal H'$ which satisfy
 $\prod_{i=1}^k C_{v_i, l} = \prod_{i=1}^k C_{v_i, l'} = 1$.
 
Repeating this and exchanging the columns of $G'$ (which corresponds to 
 reordering the edges in $\mathcal E$), and also rows of $G'$ 
 (which corresponds to reordering the basis $\{g_{r+1}, \dots, g_n\}$), 
 the existence of the pre-log curves reduces to the case where 
 all the integers
 $n_{r+1}, n_{r+2}-n_{r+1}, \dots, n_n-n_{n-1}$ are one.
In this case, we write the matrix $G'$ as $\overline G$.
In general, we obtain a matrix of the size smaller than $n-r$, 
 see Figure \ref{fig:matrix}.

Assuming $n_{r+1}, n_{r+2}-n_{r+1}, \dots, n_n-n_{n-1}$ are all one, 
 $n_n = n-r$ and $\overline G$ is an invertible square matrix
 with integer entries (the inverse matrix need not be defined over
 the integers).
Let us write the matrix $\overline G$ by entries, $\overline G = (g_{ij})$, 
 $1\leq i, j\leq n-r$.
In this case, we rearrange the order of the set of edges $\mathcal E$ so that
 the first $i_1 = m_1-3$ edges are attached to the vertex $v_1$, the next
 $i_2 = m_2-3$ edges are attached to $v_2$, and so on.
Note that $\sum_{a=1}^{k}i_a = n_n = n-r$.

\begin{figure}[h]
\[
\bordermatrix{
	&&&&&&i&&&j+1&&&&& \cr
	& * & * & * & \cdots & \cdots & \cdots & \cdots & \cdots &\cdots &\cdots& \cdots& \cdots& \cdots&
	 \cdots& \cdots \cr
	& 0 & * & * & \cdots & * & \beta & \cdots & \cdots & \cdots & \cdots& \cdots&\cdots &\cdots
	& \cdots& \cdots \cr
	& \vdots & \vdots & \vdots & \vdots & \vdots & \vdots & \vdots & \vdots & \vdots & \vdots& \vdots
	& \vdots& \vdots& \vdots & \vdots\cr
 i	& 0 & 0 & \cdots & \cdots &  0 & * &\cdots &\cdots &\cdots &\cdots& 0 & \alpha & 0 &\cdots &\cdots\cr
 	& \vdots & \vdots & \vdots & \vdots & \vdots & \vdots & \vdots & \vdots & \vdots & \vdots & \vdots
 	& \vdots& \vdots& \vdots  & \vdots\cr
 j+1 & 0 & 0  & \cdots & \cdots & \cdots & \cdots & \cdots & 0 &* &* & \cdots& \cdots & \cdots& \cdots
     & \cdots\cr
 	& \vdots & \vdots & \vdots & \vdots & \vdots & \vdots & \vdots & \vdots & \vdots & \vdots & \vdots
 	  & \vdots& \vdots& \vdots& \vdots\cr
 	& 0 & 0 & \cdots & \cdots & \cdots & \cdots & \cdots& \cdots &0 & 0 & \cdots 
 	 & * & * & \cdots& \cdots \cr
 	& 0 & 0 & \cdots & \cdots & \cdots & \cdots & \cdots& \cdots & 0 & 0& \cdots& 0 & 0 & * & *
}
\]
\caption{A possible form of the $(n-r)\times n_n$-matrix $G'$. 
The $k$-th row ($k\geq j+1$) contains at least two nonzero entries on the columns
 on which the entries of the $l$-th rows ($k+1\leq l\leq n-r$) are zero.
The upper left $j\times j$ submatrix is upper triangular with nonzero diagonals.
When the $i$-th row ($1\leq i\leq j$) contains a nonzero entry on the $m$-th column
 with $j+1\leq m\leq n_n$ (the entry $\alpha$ in the matrix), then the constant
 $\prod_{a=1}^k C_{v_a, r+i}$ can be set to one.
Then if the $i'$-th row with $i'<i$ contains a nonzero entry on the $i$-th column
 (the entry $\beta$ in the matrix), then the constant
 $\prod_{a=1}^k C_{v_a, r+i'}$ can be set to one.
By lowering these rows for which the constant $\prod_{a=1}^k C_{v_a, r+i}$
 can be set to one in the above way, we collect the remaining edges at the top part of the matrix.
When the number of these edges is $\rho$, then by reordering the columns, we obtain 
 a new matrix in which the upper left $\rho\times\rho$ matrix is upper triangular with
 nonzero diagonals, and the entries of these rows on the $k$-th column ($k\geq \rho+1$)
 are all zero.
The existence of a pre-log curve is reduced to the property of this matrix, which we write 
 by $\overline G$.
}\label{fig:matrix}
\end{figure}

Then the condition in Lemma \ref{lem:prodmap} can be written in the 
 following way.
Namely, for the existence of a pre-log curve of type
 $(\Gamma, h)$, it suffices to show that
 there are complex numbers $\zeta_1, \dots, \zeta_{n-r}$
 satisfying the following conditions:
\begin{enumerate}
\item $\zeta_i\in \Bbb C^*\setminus \{1\}$,
\item $\zeta_i\neq \zeta_j$ when $i\neq j$, 
 $\sum_{a=1}^ci_a+1\leq i, j\leq\sum_{a=1}^{c+1}i_a$ for some $0\leq c\leq k-1$, and
\item $\prod_{i=1}^{n-r}\zeta_i^{g_{ji}} = 1$ for all $j = 1, \dots, n-r$.
\end{enumerate}

Taking the logarithm, this can be further rephrased in the following form.
That is, consider the system of linear polynomials
\[
\overline G\begin{pmatrix}
x_1\\ \vdots \\ x_{n-r}
\end{pmatrix},
\]
 where $x_1, \dots, x_{n-r}$ are variables taking values in complex numbers.
Then the condition in Lemma \ref{lem:prodmap} is equivalent to the requirement
 that the equation
\[
\overline G\begin{pmatrix}
x_1\\ \vdots \\ x_{n-r}
\end{pmatrix}
 = 2\pi i \begin{pmatrix}
 a_1\\ \vdots\\ a_{n-r}
 \end{pmatrix}
\]
 has a solution for some integers $a_1, \dots, a_{n-r}$, where 
 $x_i\notin 2\pi i\Bbb Z$ and 
 $x_i-x_j\notin 2\pi i\Bbb Z$ for $i\neq j$,
 $\sum_{a=1}^ci_a+1\leq i, j\leq\sum_{a=1}^{c+1}i_a$ for some $0\leq c\leq k-1$.

Moreover, the solutions $\begin{pmatrix}
x_1\\ \vdots \\ x_{n-r}
\end{pmatrix}$ and $\begin{pmatrix}
x_1'\\ \vdots \\ x'_{n-r}
\end{pmatrix}$ such that 
 $\begin{pmatrix}
 x_1\\ \vdots \\ x_{n-r}
 \end{pmatrix}-\begin{pmatrix}
 x'_1\\ \vdots \\ x'_{n-r}
 \end{pmatrix}\in 2\pi i\Bbb Z^{n-r}$ 
 give the same solution for the original equation.
Thus, replacing $x_i$ by $X_i = \frac{x_i}{2\pi i}$, 
 we only need to look for solutions in the higher dimensional 
 open unit cube
 $I^{n-r} = (0, 1)^{n-r}$.

Since the matrix $\overline G$ is invertible, the image of $I^{n-r}$ by 
 $\overline G$ is a 
 higher dimensional open parallelotope $P\subset\Bbb R^{n-r}$
 spanned by the columns of 
 $\overline G$.
For $i\neq j$,
$\sum_{a=1}^ci_a+1\leq i, j\leq\sum_{a=1}^{c+1}i_a$ for some $0\leq c\leq k-1$, 
 let $H_{ij}$ be the hyperplanes in $\Bbb R^{n-r}$ defined by
\[
H_{ij} = \left\{\begin{pmatrix} X_1\\ \vdots \\ X_{n-r}\end{pmatrix} \; \bigg| \; X_i = X_j\right\}
\]
 and $H$ be their union
\[
H = \coprod H_{ij},
\]
 where the pair $i, j$ runs through those which satisfy the above condition.
Then, 
 summarizing the argument so far, we have the following.
\begin{thm}\label{thm:pre-logexist}
When $n_n = n-r$, the number of elements in the
 set $\prod_{i=1}^k \mathcal H_i$
 satisfying the condition of Lemma \ref{lem:prodmap}
 is equal to the number of integral points in $P\setminus \overline G(H)$.\qed
\end{thm}

Translating to the language of pre-log curves, we have the following.
We use the same notation as in the above argument.

\begin{cor}\label{cor:main}
Given a tropical curve $(\Gamma, h)$ of genus one such that the 
 number of the edges emanating from the vertices on the loop satisfies 
 $n_n = n-r$, there are pre-log curves of type $(\Gamma, h)$ 
 if and only if the set of integral points in $P\setminus \overline G(H)$ is non empty.
Moreover, the number of families of such pre-log curves is 
 given by the number of integral points in $P\setminus \overline G(H)$.\qed
\end{cor}

\begin{rem}\label{rem:family}
	The family of pre-log curves corresponding to an integral point in $P\setminus \overline G(H)$
	consists of curves obtained by gluing marked rational curves with fixed moduli associated 
	to the vertices of $\Gamma$ (the marked points correspond to nodes or intersection
	with toric divisors) as we explained in Subsection \ref{subsec:corrthm}. 
\end{rem}

Combined with Theorem \ref{thm:loop}, we have the following
 characterization of smoothable tropical curves of genus one.
 
\begin{cor}\label{cor:genusonesm}
Let $(\Gamma, h)$ be a tropical curve of genus one satisfying Assumption A.
Assume that the directions of the edges emanating from the vertices on the loop
 spans the whole $\Bbb R^n$.
Then $(\Gamma, h)$ is smoothable if and only if the associated 
 set $P\setminus \overline G(H)$ contains an integral point.\qed
\end{cor}

\begin{rem}\label{rem:abundant}
In this corollary, we do not assume $n_n = n-r$, and the size of 
 the matrix $\overline G$ will be in general smaller than $(n-r)\times (n-r)$.
When this happens, the family of pre-log curves mentioned in Corollary \ref{cor:main}
 and Remark \ref{rem:family} will acquire additional freedom coming from 
 deforming some of the pointed rational curves in the components of $C_0$.
\end{rem}

\begin{rem}
The number of families obtained by smoothing pre-log curves of 
 type $(\Gamma, h)$ is in general larger than the number of 
 integral points in $P\setminus G(H)$ since the same pre-log curve can carry
 several log structures relevant to our purpose.
For a given pre-log curve of type $(\Gamma, h)$, 
 the number of different families of smoothings can be calculated from 
 the edge weights of $\Gamma$ using log deformation theory,
 as in \cite{NS}.
\end{rem}

\begin{rem}
In \cite{Tor}, K. Torchiani defined a combinatorial counting number of tropical curves
 of genus one with nice invariance properties, using combinatorial moduli space of
 well-spaced tropical curves.
There she gave positive multiplicity to any loop of the kind we considered in this section.
In particular, the tropical curves like the one given in Example \ref{ex:1} also contribute to 
 her counting number.
The result in this section shows that the moduli of classical curves is subtler than the
 combinatorial moduli of tropical curves.
\end{rem}

\begin{rem}
The argument in this section can be used to study the existence of 
 pre-log curves of type $(\Gamma, h)$ with higher genus.
In the cases of higher genus, such a study becomes more important since
 as we saw in Example \ref{ex:1.5}, even in 3-valent cases, 
 the existence of pre-log curves will be a nontrivial problem.
 
Also, it turns out that when a suitable transversality 
 condition holds (in the  case of genus one, this is the condition that 
 the directions of the edges emanating from 
 the vertices on the loop span $\Bbb R^n$),
 the existence of pre-log curves immediately implies the smoothability
 of the tropical curve in the sense of Definition \ref{def:smoothable}.
These issues will be treated elsewhere.
\end{rem}

\subsection{Examples}
In this subsection, we give a few examples concerning 
 Theorems \ref{thm:loop} and \ref{thm:pre-logexist}.

\begin{example}\label{ex:4}
First, let us look at Examples \ref{ex:1} and \ref{ex:2}
 from the viewpoint of Theorems \ref{thm:loop} and \ref{thm:pre-logexist}.
In these cases, $n = 3, r = 2$ and $n_n = n-r = 1$.
Therefore, the matrix $\overline G$ is just an integer and 
 it is 1 in the case of Example \ref{ex:1} and 2 in the case of 
 Example \ref{ex:2}.
The open parallelotope $P$ is $(0, 1)$ and $(0, 2)$, respectively, 
 and the space $H$ is empty in both cases.
Thus, $P\setminus G(H)$ contains no integral point 
 in the case of Example \ref{ex:1}
 and one integral point in the case of Example \ref{ex:2},
 recovering the results discussed in these examples.
 
Similar situation happens in Example \ref{ex:super1}.
When the edges $A$ and $B$ are both contracted (see the picture
 on the left of Figure \ref{fig:0H}), 
 then there are infinitely many isomorphism classes of corresponding 
 pre-log curves by Corollary \ref{cor:genusonesm} and Remark \ref{rem:abundant}, 
 and all of them are smoothable.

\begin{figure}[h]
\includegraphics{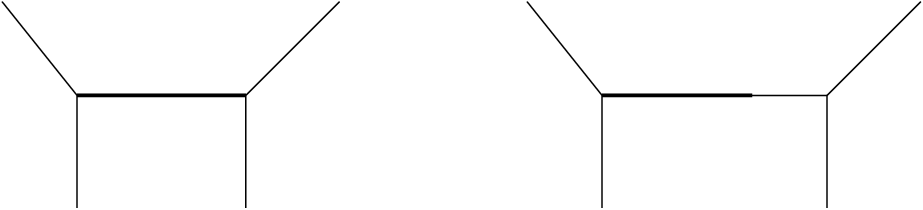}
\caption{}\label{fig:0H}	
\end{figure}

Now consider the case when only the edge  
 $A$  is contracted
 (see the picture on the right of Figure \ref{fig:0H}).
In this case, since the weights of the edges outside the loop are two, 
 it follows that the open parallelotope is $(0, 2)$.
Therefore, there is one
 pre-log curve of the corresponding type  and it is smoothable.
Contrary to the case of Example \ref{ex:super1},
 the length of $B$ does not matter to the smoothability.

\end{example}

\begin{example}
Consider the tropical curve in $\Bbb R^4$ whose image is 
 given by Figure \ref{fig:0E}.

\begin{figure}[h]
\includegraphics{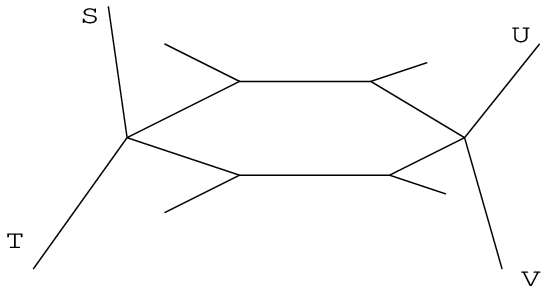}
\caption{}\label{fig:0E}
\end{figure}

Here the edges other than $S, T, U$ and $V$ are contained in 
 the two dimensional affine subspace $A$ spanned by the directions of the 
 edges in the loop.
Let $e_1, e_2, e_3, e_4$ be a basis of $\Bbb Z^4\subset \Bbb R^4$
 and assume that the two dimensional plane $\bar A$ parallel to $A$ is 
 spanned by $e_1$ and $e_2$.

Let $g_1, g_2$ be the primitive integral generators of the edges $S$ and $U$,
 and assume the weights of these edges are one.
Assume these generators are given by
\[
g_1 = \begin{pmatrix}
a_1\\ b_1\\ 1 \\ 0
\end{pmatrix},\;\;
g_2 = \begin{pmatrix}
a_2\\ b_2\\ 3 \\ 4
\end{pmatrix},
\]
 where $a_i, b_i$ are some integers.
In this case, $P$ is the open parallelogram spanned by the edges
 $\begin{pmatrix}
 1\\ 0
 \end{pmatrix}, 
 \begin{pmatrix}
 3\\ 4
 \end{pmatrix},
 $
 and $H$ is the empty set.
Thus, the relevant integral points are three  points 
 $\begin{pmatrix}
  1 \\ 1 
 \end{pmatrix},
 \begin{pmatrix}
 2 \\ 2
 \end{pmatrix},
 \begin{pmatrix}
 3 \\ 3
 \end{pmatrix}$.
Therefore, there are three families of pre-log curves of the given type, and
 all of them are smoothable by Theorem \ref{thm:loop}.

On the other hand, consider the tropical curve in $\Bbb R^4$ 
 whose image is given by Figure \ref{fig:0F}.

\begin{figure}[h]
	\includegraphics{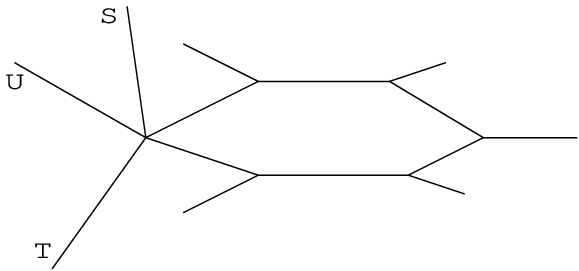}
	\caption{}\label{fig:0F}
\end{figure}

This time the edges other than $S, T$ and $U$ are contained in 
 the two dimensional affine subspace $A$ spanned by the directions of the 
 edges in the loop.

Assume that the edges $S$ and $U$ are generated by the same vectors
 $g_1$ and $g_2$ as before,
 and also assume that the weights of these edges are one.

In this case, $P$ is the same parallelogram as the above example, 
 but the space $H$ is the diagonal of $\Bbb R^2$.
The space $G(H)$ contains all the integral points
 $\begin{pmatrix}
 1 \\ 1 
 \end{pmatrix},
 \begin{pmatrix}
 2 \\ 2
 \end{pmatrix},
 \begin{pmatrix}
 3 \\ 3
 \end{pmatrix}$,
 so that in this case there are no pre-log curves of the type given by 
 the tropical curve in Figure \ref{fig:0F}.
In particular, this tropical curve is not smoothable. 
\end{example}

\end{document}